\documentclass[12pt]{article}
\usepackage{geometry}
\def\subclassname{{\bfseries Mathematics Subject Classification
(2010)}\enspace}
\def\subclass#1{\par\addvspace\medskipamount{\rightskip=0pt plus1cm
\def\and{\ifhmode\unskip\nobreak\fi\ $\cdot$
}\noindent\subclassname\ignorespaces#1\par}}
\def\PACSname{\textbf{PACS}\enspace}
\def\PACS#1{\par\addvspace\medskipamount{\rightskip=0pt plus1cm
\def\and{\ifhmode\unskip\nobreak\fi\ $\cdot$
}\noindent\PACSname\ignorespaces#1\par}}
\def\CRclassname{{\bfseries CR Subject Classification}\enspace}
\def\CRclass#1{\par\addvspace\medskipamount{\rightskip=0pt plus1cm
\def\and{\ifhmode\unskip\nobreak\fi\ $\cdot$
}\noindent\CRclassname\ignorespaces#1\par}}
\usepackage{titlesec}

\titleformat{\section}
  {\normalfont\fontsize{14}{15}\bfseries}{\thesection}{1em}{}
\titleformat{\subsection}
{\normalfont\fontsize{14}{15}\bfseries}{\thesubsection}{1em}{}

\newcommand{\footremember}[2]{%
   \footnote{#2}
    \newcounter{#1}
    \setcounter{#1}{\value{footnote}}%
}
\providecommand{\keywords}[1]
{
  \textbf{\textit{Keywords---}} #1
}

\usepackage{graphicx}
\usepackage{float}
\usepackage{caption}
\usepackage{subcaption}%
\usepackage{amsmath,amsthm}      
\usepackage{color}
\usepackage{latexsym}
\usepackage{hyperref}
\numberwithin{equation}{section}
\newtheorem{defn}{Definition}[section]
\newtheorem{theorem}{Theorem}[section]
\newtheorem{remark}{Remark}[section]

\usepackage[svgnames]{xcolor}
\usepackage{soul}
\usepackage{ulem}

\usepackage{lineno} 

\begin{document}

\title{Entropy stable flux correction for hydrostatic reconstruction scheme for shallow water flows 
}


\author{Sergii Kivva\footremember{alley}{Institute of Mathematical Machines and System Problems, National Academy of Sciences, Ukraine}
}



\date{\Large}

\maketitle

\begin{abstract}
First-order hydrostatic reconstruction (HR) schemes for shallow water equations are highly diffusive whereas high-order schemes can produce entropy-violating solutions. Our goal is to develop a flux correction with maximum antidiffusive fluxes to obtain entropy solutions of shallow water equations with variable bottom topography. 
For this purpose, we consider a hybrid explicit HR scheme whose flux is a convex combination of first-order Rusanov flux and high-order flux. The conditions under which the explicit first-order HR scheme for shallow water equations satisfies the fully discrete entropy inequality have been studied.
The flux limiters for the hybrid scheme are calculated from a corresponding optimization problem. Constraints for the optimization problem consist of inequalities that are valid for the first-order HR scheme and applied to the hybrid scheme. We apply the discrete cell entropy inequality with the proper numerical entropy flux to single out a physically relevant solution to the shallow water equations. A nontrivial approximate solution of the optimization problem yields expressions to compute the required flux limiters. 
Numerical results of testing various HR schemes on different benchmarks are presented.

 \subclass{MSC 65M08 \and MSC 65M22 \and 76M12 }
\end{abstract}

\section{Introduction}
\label{intro}

In this paper, we consider a design of entropy stable flux correction for a hydrostatic reconstruction scheme for  shallow water equations with variable bottom topography. For simplicity, without loss of generality, we focus on the Saint-Venant system of one-dimensional shallow water equations, given by
\begin{equation}
\label{eq:1}
\begin{split}
& \partial_t h + \partial_x Q = 0,  \\
& \partial_t Q +\partial_x \left(\frac{Q^2}{h} +g \frac{h^2}{2}\right) = 
-gh \partial_x z,
\end{split}
\end{equation}
subject to the initial conditions
\begin{equation}
\label{eq:2}
h(x,0) = h_0(x), \quad Q(x,0) = Q_0(x),
\end{equation}
where $h(x,t)$ is the water depth, $Q(x,t)$ is the water discharge, $g$ is the gravitational constant, and $z(x)$ is the bottom topography. The system \eqref{eq:1} is considered in a certain spatial domain $D$, and if $D \neq {R}$, then on the boundary of $D$ a corresponding boundary conditions should be specified.

In vector form, the system \eqref{eq:1} can be written as
\begin{equation} 
\label{eq:3} 
\frac{\partial \boldsymbol u}{\partial t}+\frac{\partial }{\partial  x} \boldsymbol f( \boldsymbol u) = \boldsymbol s,   
\end{equation}  
where $\boldsymbol u= (h,Q)^T$ is the vector of conserved variables, $\boldsymbol f=(Q,Q^2/h+gh^2/2)^T $ is the flux vector, and $\boldsymbol s=(0, ghz_x)^T$ is the source vector.

It is well known \cite{Lax} that solutions of \eqref{eq:1}-\eqref{eq:2} may develop singularities in finite time even for a smooth initial condition. Hence, we should interpret \eqref{eq:1} in the sense of distribution and search for weak solutions. However, such weak solutions are not unique. To single out a unique physically relevant weak solution, the latter should satisfy
\begin{equation}
\label{eq:4} 
\frac{\partial U(\boldsymbol u)}{\partial t} + \frac{\partial F(\boldsymbol u)}{\partial x} \le 0	
\end{equation}
in the sense of distribution for every entropy pair ($U,F$). Here $U$ is a convex function of $\boldsymbol u$, the so-called entropy function, and $F$ is its entropy flux that satisfies
\begin{equation}
\label{eq:5}
F^T_{\boldsymbol u}(\boldsymbol u) = U_{\boldsymbol u}(\boldsymbol u) \boldsymbol f_{\boldsymbol u}(\boldsymbol u).
\end{equation}

For shallow water equations \eqref{eq:1} with bottom topography $z(x)$, the total energy
\begin{equation}
\label{eq:6}
U(\boldsymbol u) = \frac{1}{2} \left(\frac{Q^2}{h} +g h^2\right) +ghz,
\end{equation}
serves as an entropy function with entropy flux
\begin{equation}
\label{eq:7}
F(\boldsymbol u) = \frac{Q^3}{2h^2} +ghQ +gQz.
\end{equation}

We discretize \eqref{eq:3} by the difference scheme 
\begin{equation}
\label{eq:8} 
\frac{1}{\Delta t}(\boldsymbol{\hat v_i} - \boldsymbol{v_i}) + \frac{1}{\Delta x}\left[ \boldsymbol g_{i+1/2} - 
\boldsymbol g_{i-1/2} \right] = \boldsymbol s_i, 								
\end{equation}
where the numerical flux $\boldsymbol g_{i+1/2}$ is calculated as
\begin{equation}
\label{eq:9} 
\boldsymbol g_{i+1/2} = \boldsymbol g_{i+1/2}^L + {\alpha _{i+1/2}}\left[ {\boldsymbol g_{i+1/2}^H - \boldsymbol g_{i+1/2}^L} \right]. 		
\end{equation}
Here, $\boldsymbol v_i = \boldsymbol v(x_i,t)= (y(x_i,t),q(x_i,t))^T$ is the discrete solution at the grid point (${x_i} = i\Delta x,t)$; $\boldsymbol {\hat v_i} = \boldsymbol v({x_i},t + \Delta t)$; $\Delta x $ and $\Delta t $ are the spatial and temporal computational grid size, respectively. $\boldsymbol g_{i+1/2}^H$ and $\boldsymbol g_{i+1/2}^L$ are a high-order and low-order numerical fluxes such that $\boldsymbol g_{i+1/2} = \boldsymbol g( \boldsymbol v_{i-l+1},...,\boldsymbol v_{i+r})$ is the Lipschitz continuous numerical flux consistent with the differential flux, that is $\boldsymbol g(\boldsymbol u,...,\boldsymbol u) = \boldsymbol f(\boldsymbol u)$ for all flux-limiters ${\alpha _{i+1/2}} \in \left[ {0,1} \right]$.

The expression in square brackets on the right-hand side of \eqref{eq:5} can be considered as an antidiffusive flux. For flux-correction we compute the flux limiters ${\alpha _{i+1/2}}$ as an approximate solution to the corresponding optimization problem. 
The classical two-step Flux-Corrected Transport (FCT) was firstly developed by Oran and Boris~\cite{b2} to solve the transient continuity equation. 
The procedure of two-step flux correction consists of computing the time-advanced low-order solution in the first step and correcting the solution by adding antidiffusive fluxes in the second step to produce accurate and monotone results. The antidiffusive fluxes, which define as the difference between the high and low-order fluxes, are limited in such a way that neither new extrema are created nor existing extrema are increased. 
Later, Zalesak~\cite{b4,b3} extended FCT to multidimensional explicit difference schemes. In~\cite{b3}, using characteristic variables, Zalesak proposed FCT algorithms for nonlinear systems of conservation laws. 
Several implicit FEM-FCT schemes for unstructured grids were proposed by Kuzmin and coworkers~\cite{b5,b8}. However, the known FCT algorithms do not guarantee entropy solutions for hyperbolic conservation laws.
 
We discretize the entropy inequality \eqref{eq:4} as follows 
\begin{equation} 
\label{eq:10} 
 U(\boldsymbol {\hat v_i}) - U({\boldsymbol v_i}) + \frac{{\Delta t}}{{\Delta x}}\left[ {G_{i+1/2}} - {G_{i-1/2}} \right] \le 0 
\end{equation}
where ${G_{i+1/2}} = G({\boldsymbol v_{i-l+1}},...,\boldsymbol {v_{i+r}})$ is the numerical entropy flux consistent with the differential one $ G(\boldsymbol u,...,\boldsymbol u) = F(\boldsymbol u)$. 

A difference scheme \eqref{eq:8} is called {\it entropy stable} if computed solutions satisfy the discrete cell entropy inequality \eqref{eq:6}. We mention here the pioneering studies of entropy stable schemes by Lax~\cite{Lax}. 
Entropy stable schemes were developed by several authors \cite{b35,b38,b43,b39,b42,b16,b12,b41}. 
To single out a physically relevant solution, we use the so-called proper numerical entropy flux, the concept of which was formulated by Merriam~\cite{b9} and Sonar~\cite{b10}. Zhao and Wu~\cite{b11} proved that three-point monotone semi-discrete schemes in conservative form satisfy the corresponding semi-discrete entropy inequality with the proper numerical entropy flux. Fully discrete entropy stable schemes with the proper numerical entropy flux for scalar conservation laws were obtained in~\cite{kivva2021,kivva2022}. The numerical entropy flux $G({\boldsymbol v_{i-l+1}},...,{\boldsymbol v_{i+r}})$ for $F$ is not unique. The distinguishing feature of the proper numerical entropy flux among others is that it satisfies property \eqref{eq:5} of the differential entropy flux. 
 
In this paper,  we apply a first-order hydrostatic reconstruction (HR) scheme as a low-order scheme to design flux correction. 
A first-order HR scheme was originally developed by Audusse et al.~\cite{Audusse2004}, and it does not properly account for the acceleration due to a sloped bottom~\cite{Delestre2012} for shallow downhill flow. Morales de Luna et al.~\cite{Morales2013} improved the original first-order HR scheme for partially wet interfaces. Using a technique of subcell reconstructions, Chen and Noelle~\cite{Chen2017} proposed a new reconstruction with a better approximation of the source term for shallow downhill flows. The main properties of the original HR scheme or its modifications are positivity preserving, well-balanced, and satisfying a semi-discrete in-cell entropy inequality. Unfortunately, it is well known that semi-discrete entropy inequalities are insufficient to obtain a suitable convergence to the entropy weak solution or to get relevant energy estimates. Audusse et al.~\cite{Audusse2016} showed that the HR scheme combined with a kinetic solver satisfies a fully discrete entropy inequality with an error term coming from the topography. Thus, we can expect the convergence of this scheme for Lipschitz continuous bathymetry.  
Berthon et al.~\cite{Berthon2019} suggested to introduce artificial viscosity into the HR scheme to get fully discrete entropy inequalities.
 
Using the approach proposed in~\cite{kivva2021,kivva2022},  we construct a flux correction for 1D shallow water equations \eqref{eq:1} to obtain numerical entropy solutions for which the antidiffusive fluxes are maximal. For this, the flux limiters for the hybrid scheme \eqref{eq:8}-\eqref{eq:9} are computed from the optimization problem with constraints that are valid for the low-order scheme. An approximate solution to the optimization problem yields the desired flux correction formulas.
Moreover, considering the flux limiters as functions of the numerical solution, we prove the unique solvability of the hybrid scheme \eqref{eq:8}-\eqref{eq:9} under general assumptions on them. We show that the approximate limiters satisfy the assumptions under which the hybrid scheme has a unique solution.
The developed approach is a novel view on the known FCT method.

The paper is organized as follows. In Section \ref{sec:1}, we present estimates that are valid for an explicit first-order HR scheme with the Rusanov numerical flux. 
Section~\ref{sec:2} defines the proper numerical entropy flux and studies the conditions under which the explicit HR scheme for homogeneous and inhomogeneous shallow water equations satisfies the fully discrete entropy inequality.
The unique solvability of flux correction for the HR scheme, the optimization problem for finding flux limiters, and the algorithm for its solution are described in Section \ref{sec:3}.
An approximate solution of the optimization problem is derived in Section \ref{sec:4}. The results of numerical experiments with different HR schemes are given in Section \ref{sec:5}. Concluding remarks are drawn in Section \ref{sec:6}.

\section{First-Order Hydrostatic Reconstruction Scheme}
\label{sec:1}

We consider an explicit first-order HR scheme of Chen and Noelle~\cite{Chen2017} in the form
\begin{equation}
\label{eq:1.1}
\boldsymbol {\hat v_i} - \boldsymbol {v_i} + \frac{{\Delta t}}{{\Delta x}}\left[ \boldsymbol{g_{i+1/2}^L(\boldsymbol{v^-_{i+1/2}},\boldsymbol{v^+_{i+1/2}}) - \boldsymbol g_{i-1/2}^L(\boldsymbol{v^-_{i-1/2}},\boldsymbol{v^+_{i+1/2}})} \right] = \Delta t \; \boldsymbol s_i,	
\end{equation}
where $\boldsymbol g_{i+1/2}^L$ is the Rusanov numerical flux~\cite{b14} consistent with the differential flux $\boldsymbol f$ and given by
\begin{equation}
\label{eq:1.2}
\boldsymbol g_{i+1/2}^L(\boldsymbol{v^-_{i+1/2}},\boldsymbol{v^+_{i+1/2}}) = \frac{1}{2} \left( \boldsymbol f(\boldsymbol v^-_{i+1/2}) + \boldsymbol f(\boldsymbol v^+_{i+1/2}) -c_{i+1/2} (\boldsymbol v^+_{i+1/2} -\boldsymbol v^-_{i+1/2}) \right).
\end{equation}

The vectors of conservative variables $\boldsymbol v^\pm_{i+1/2}$ are given by
\begin{equation}
\label{eq:1.3}
\boldsymbol v^-_{i+1/2} = \left( \begin{array}{c}
y^-_{i+1/2} \\ y^-_{i+1/2} u_i
\end{array} \right), \quad 
\boldsymbol v^+_{i+1/2} = \left( \begin{array}{c}
y^+_{i+1/2} \\ y^+_{i+1/2} u_{i+1}
\end{array} \right), \quad
u_i = \frac{ \sqrt{2} y_i q_i}{\sqrt{ y_i^4 +\max \left( y_i^4,\epsilon  \right)}},
\end{equation}
where $\epsilon $ is a small a-priori chosen positive number. The water depths are calculated as
\begin{equation}
\label{eq:1.4}
y^-_{i+1/2} = \min \left( w_i - z_{i+1/2},y_i \right), \quad
y^+_{i+1/2} = \min \left( w_{i+1} - z_{i+1/2},y_{i+1} \right)
\end{equation}
with water levels $w_i = z_i +y_i$, and the cell interface bottom
\begin{equation}
\label{eq:1.5}
z_{i+1/2} = \min \left( \max \left( z_i, z_{i+1} \right), \min \left( w_i, w_{i+1} \right) \right).
\end{equation}

The source term $\boldsymbol s_i = - \boldsymbol s^+_{i-1/2} + \boldsymbol s^-_{i+1/2} =( 0, -s^+_{i-1/2})^T +(0,s^-_{i+1/2} )^T$ is discretized as
\begin{equation}
\label{eq:1.6}
\begin{split}
 & s^-_{i+1/2} = -g \frac{y_i + y^-_{i+1/2}}{2} \frac{z_i -z_{i+1/2}}{\Delta x } , \\
 & s^+_{i+1/2} = -g \frac{y^+_{i+1/2} +y_{i+1} }{2} \frac{z_{i+1} -z_{i+1/2}}{\Delta x } .
\end{split}
\end{equation}

Finally, the local speed $c_{i+1/2}$ in \eqref{eq:1.2} is calculated using the eigenvalues of the Jacobian $\boldsymbol f_{\boldsymbol u} (\boldsymbol u)$ as follows
\begin{equation}
\label{eq:1.7}
c_{i+1/2} = \max \left( \vert u_i \vert +\sqrt{g \, y^-_{i+1/2}},  \; \vert u_{i+1} \vert +\sqrt{g \, y^+_{i+1/2}} \right) .
\end{equation}

The following theorem gives estimates for the numerical solution of the HR scheme \eqref{eq:1.1}.

\begin{theorem}
\label{th:1.1}
 Assume that $c_{i+1/2} \geq \max \left( \vert u_i \vert ,\, \vert u_{i+1} \vert \right)$ for all $i$. Then for
\begin{equation}
\label{eq:1.8}
\Delta t\; \le \frac{2 \Delta x} {\max \limits_{i} ( c_{i+1/2} - u_{i+1} + c_{i-1/2} + u_{i-1} )},	
\end{equation}
the following inequalities hold for the numerical solution of the system of equations \eqref{eq:1.1}-\eqref{eq:1.2}
\begin{equation}
\label{eq:1.9}
\begin{split}
& \frac{\Delta x}{\Delta t} \min(w_i,w^-_{i-1/2},w^+_{i+1/2}) \leq \frac{\Delta x}{\Delta t} \hat w_i - \frac{u_{i} +u_{i+1}}{2} z_{i+1/2} -\left( \frac{c_{i+1/2} -u_{i+1}}{2} w_i -\frac{c_{i+1/2} +u_{i}}{2} w^-_{i+1/2} \right)\\
& +\frac{u_{i} +u_{i-1}}{2} z_{i-1/2} -\left( \frac{c_{i-1/2} +u_{i-1}}{2} w_i -\frac{c_{i-1/2} -u_{i}}{2} w^+_{i-1/2} \right)\leq \frac{\Delta x}{\Delta t} \max(w_i,w^-_{i-1/2},w^+_{i+1/2}),   
\end{split} 
\end{equation} 	
\begin{equation}
\label{eq:1.10}
\begin{split}
\frac{\Delta x}{\Delta t} & \min(q_i,q^-_{i-1/2},q^+_{i+1/2}) \leq \frac{\Delta x}{\Delta t} \hat q_i -\left( \frac{c_{i+1/2} -u_{i+1}}{2} q_i -\frac{c_{i+1/2} +u_{i}}{2} q^-_{i+1/2} \right)\\
& + \frac{g}{2} \left[ \frac{1}{2} \left(y^{-^2}_{i+1/2} +y^{+^2}_{i+1/2} \right) +(y_i +y^-_{i+1/2})(z_{i+1/2} -z_i) \right] \\
& - \frac{g}{2} \left[ \frac{1}{2} \left(y^{-^2}_{i-1/2} +y^{+^2}_{i-1/2} \right) +(y_i +y^+_{i-1/2})(z_{i-1/2} -z_i) \right] \\
&  -\left( \frac{c_{i-1/2} +u_{i-1}}{2} q_i -\frac{c_{i-1/2} -u_{i}}{2} q^+_{i-1/2} \right)\leq \frac{\Delta x}{\Delta t} \max(q_i,q^-_{i-1/2},q^+_{i+1/2}),   
\end{split} 
\end{equation} 		
where ${w^{\pm}_{i+1/2} = y^{\pm}_{i+1/2} +z_{i+1/2}}$. 
\end{theorem}
\begin{proof}
Let us prove inequalities \eqref{eq:1.9}. Inequalities \eqref{eq:1.10} are proved similarly.

We rewrite the equation \eqref{eq:1.1} for the conservative variable $y_i$ in the form
\begin{equation}
\label{eq:1.11}
\begin{split}
 \frac{\Delta x}{\Delta t} \hat y_i = & \frac{\Delta x}{\Delta t}y_i  - \frac{c_{i+1/2} + u_i}{2} y^-_{i+1/2} - \frac{c_{i-1/2}-u_i}{2} y^+_{i-1/2} \\
+ & \frac{c_{i+1/2} -u_{i+1}}{2} y^+_{i+1/2}  + \frac{c_{i-1/2} +u_{i-1}}{2} y^-_{i-1/2}.
\end{split}
\end{equation}

Substituting the water level $w_i$ in \eqref{eq:1.11} instead of the water depth $y_i$, we obtain
\begin{equation}
\label{eq:1.12}
\begin{split}
 & \frac{\Delta x}{\Delta t} \hat w_i = \left( \frac{\Delta x}{\Delta t} -\frac{c_{i+1/2} -u_{i+1}}{2} - \frac{c_{i-1/2} +u_{i-1}}{2} \right) w_i  +\frac{c_{i+1/2} -u_{i+1}}{2} w^+_{i+1/2} \\
& +\frac{c_{i-1/2} +u_{i-1}}{2}  w^-_{i-1/2} +\frac{u_{i} +u_{i+1}}{2} z_{i+1/2} -\frac{u_{i} +u_{i-1}}{2} z_{i-1/2} \\
& +\left( \frac{c_{i+1/2} -u_{i+1}}{2} w_i -\frac{c_{i+1/2} +u_{i}}{2} w^-_{i+1/2} \right) +\left( \frac{c_{i-1/2} +u_{i-1}}{2} w_i -\frac{c_{i-1/2} -u_{i}}{2} w^+_{i-1/2} \right). 
\end{split}
\end{equation}

Note that under the condition \eqref{eq:1.8}, the first three terms in the right-hand side of \eqref{eq:1.12} are a convex combination of  $w_i,w^+_{i+1/2}$, and $w^-_{i-1/2}$, which proves the theorem.
\end{proof}
\begin{remark}
We note that if $c_{i+1/2}$ satisfies the following inequalities
\begin{equation}
\label{eq:1.13}
\begin{split}
& \frac{\Delta x}{\Delta t} y_i -\frac{c_{i+1/2} +u_{i}}{2} y^-_{i+1/2} - \frac{c_{i-1/2} -u_{i}}{2}  y^+_{i-1/2} \ge 0, \\
& c_{i+1/2} -u_{i+1} \ge 0, \quad
 c_{i-1/2} +u_{i-1} \ge 0, 
\end{split}
\end{equation}
then the difference scheme \eqref{eq:1.11} preserves the non-negativity of the water depth $y$. 

\end{remark}

\section{Cell Entropy Inequality for Fully Discrete HR Scheme }
\label{sec:2}

In this section we study the cell entropy inequality for the fully discrete HR scheme \eqref{eq:1.1}-\eqref{eq:1.2}.

We consider a homogeneous three-point low-order scheme in the form
\begin{equation}
\label{eq:2.1}
\boldsymbol {\hat v_i} - \boldsymbol {v_i} + \frac{{\Delta t}}{{\Delta x}}\left[ {\boldsymbol g_{i+1/2}^L({\boldsymbol v_i},{\boldsymbol v_{i+1}}) - \boldsymbol  g_{i-1/2}^L({\boldsymbol v_{i-1}},{\boldsymbol v_i})} \right] = 0,
\end{equation}
where the low-order numerical flux $\boldsymbol g_{i+1/2}^L = \boldsymbol g(\boldsymbol v,\boldsymbol w)$ is consistent with the smooth differential flux $\boldsymbol f(\boldsymbol u):R^m\rightarrow R^m$ of the conservative variables $\boldsymbol u=(u^1,\ldots,u^m)^T$.

We define the numerical entropy flux as follows.
\begin{defn}
 Numerical entropy flux $G(\boldsymbol v_{i-l+1},...,\boldsymbol v_{i+r})$ of the difference scheme \eqref{eq:2.1} is called proper if for any $\boldsymbol v_{i-l+1},...,\boldsymbol v_{i+r}\in R^m$ we have
\begin{equation}
\label{eq:2.2}  
\begin{split} 
& \frac{\partial }{\partial v^j_p} G(\boldsymbol v_{i-l+1},...,\boldsymbol v_{i+r}) = 
\sum\limits_k \frac{\partial U(\boldsymbol v_p)}{\partial v^k_p} 
\frac{\partial }{\partial v^j_p} 
g^k(\boldsymbol v_{i-l+1},...,\boldsymbol v_{i+r}), \quad p = i - l + 1,...,i + r.
\end{split} 		 
\end{equation} 
\end{defn}

Then the proper numerical entropy flux for the difference scheme \eqref{eq:8} and \eqref{eq:9} can be written in the form

\begin{equation}
\label{eq:2.3}  
  {G_{i+1/2}} = G_{i+1/2}^L + {\alpha _{i+1/2}}\left( {G_{i+1/2}^H - G_{i+1/2}^L} \right), 							\end{equation}
where $G_{i+1/2}^L$ and $G_{i+1/2}^H$ are the low-order and high-order proper numerical entropy fluxes corresponding to the numerical fluxes $\boldsymbol g_{i+1/2}^L$ and $\boldsymbol g_{i+1/2}^H$.

\begin{theorem}
\label{th:2.1}
Suppose that $\boldsymbol f:R^m \rightarrow R^m$ is  hemicontinuosly Gateaux differentiable, $U:~R^m \rightarrow R$ is a strictly convex function with a hemicontinuos second Gateaux derivative. If matrices $U''(\boldsymbol w)g'_{\boldsymbol u}(\boldsymbol u,\boldsymbol v_i)$ and $U''(\boldsymbol w)g'_{\boldsymbol u}(\boldsymbol v_i,\boldsymbol u)$ are positive and negative definite, respectively, for any $\boldsymbol u, \boldsymbol w \in R^m$, $\Delta t$ satisfies the inequality  
\begin{equation}
\label{eq:2.4} 
\begin{split} 
\Delta t & \mathop {\max }\limits_{\boldsymbol s \in ({\boldsymbol  v_i},\boldsymbol {\hat v_i})} \lambda \left( U''(\boldsymbol s) \right) \;\langle {\left( {\boldsymbol g_{i+1/2}^L - \boldsymbol g_{i-1/2}^L} \right)}, {\left( {\boldsymbol g_{i+1/2}^L - \boldsymbol g_{i-1/2}^L} \right)} \rangle \\
 & \le 2\Delta x \left[ \langle U'(\boldsymbol {v_i}), \; \left( \boldsymbol g_{i+1/2}^L - \boldsymbol g_{i-1/2}^L  \right) \rangle - G_{i+1/2}^L + G_{i-1/2}^L \right], 
\end{split}		  		
\end{equation}
then the fully discrete scheme \eqref{eq:2.1} satisfies the discrete cell entropy inequality
\begin{equation}
\label{eq:2.5}  
U(\boldsymbol {\hat v_i}) - U(\boldsymbol {v_i}) + \frac{{\Delta t}}{{\Delta x}}\left[ {G_{i+1/2}^L(\boldsymbol {v_i},\boldsymbol v_{i+1}) - G_{i-1/2}^L( \boldsymbol {v_{i-1}},\boldsymbol {v_i})} \right] \le 0.		
\end{equation}
where $G_{i+1/2}^L$ is the proper numerical entropy flux corresponding to the numerical flux $\boldsymbol g_{i+1/2}^L$,
  $\langle \cdot, \cdot \rangle$ denotes the Euclidean inner product.
\end{theorem}

\begin{proof}
Multiplying  \eqref{eq:2.1} by $U'(\boldsymbol {v_i})$ and subtracting it from the left-hand side of \eqref{eq:10}, we get
\begin{equation}
\label{eq:2.6}  
\begin{split}  
 & \quad U(\boldsymbol {\hat v_i}) - U(\boldsymbol {v_i}) \; + \; \frac{{\Delta t}}{{\Delta x}}\left[ {G_{i+1/2}^L - G_{i-1/2}^L} \right] \\
  =  U(\boldsymbol {\hat v_i}) - U(\boldsymbol {v_i}) - & \langle U'(\boldsymbol {v_i}), \left( {{\boldsymbol {\hat v}_i} - \boldsymbol {v_i}} \right) \rangle  \; + \; \frac{{\Delta t}}{{\Delta x}}\left[ G_{i+1/2}^L - G_{i-1/2}^L - \left\langle U'(\boldsymbol {v_i}), \left( \boldsymbol g_{i+1/2}^L - \boldsymbol g_{i-1/2}^L \right) \right\rangle \right]  \\
& = \frac{1}{2}{\left( {\frac{\Delta t}{\Delta x}} \right)^2}  \left\langle U''(\boldsymbol s )  \left({\boldsymbol g_{i+1/2}^L - \boldsymbol  g_{i-1/2}^L} \right),\left({\boldsymbol g_{i+1/2}^L - \boldsymbol  g_{i-1/2}^L} \right) \right\rangle \;  \\
& + \; \frac{{\Delta t}}{{\Delta x}}\left[ G_{i+1/2}^L - F(\boldsymbol v_i) - \left\langle U'(\boldsymbol {v_i}), \left( \boldsymbol g_{i+1/2}^L - \boldsymbol f(\boldsymbol v_i) \right) \right\rangle \right] \\
& + \; \frac{{\Delta t}}{{\Delta x}}\left[ F(\boldsymbol v_i) - G_{i-1/2}^L - \left\langle U'(\boldsymbol {v_i}), \left( \boldsymbol f(\boldsymbol v_i) - \boldsymbol g_{i-1/2}^L \right) \right\rangle \right]   
\end{split}
\end{equation}
where $\boldsymbol s  = \theta \boldsymbol  {\hat v_i} + (1 - \theta ) \boldsymbol {v_i}$, $0 < \theta  < 1$. 

It is easy to see that the first term on the right-hand side of \eqref{eq:2.6} is non-negative. Now we show that the second and third terms in square brackets are non-positive. Indeed, we rewrite the second and third terms as follows
\begin{equation}
\label{eq:2.7}  
\begin{split} 
& G^L(\boldsymbol v_{i},\boldsymbol v_{i+1})-G^L(\boldsymbol v_{i},\boldsymbol v_i) 
- \sum\limits_j \frac{\partial U}{\partial v^j}(\boldsymbol v_i)  \left[ g^{L,j}(\boldsymbol v_{i},\boldsymbol v_{i+1}) -g^{L,j}(\boldsymbol v_{i},\boldsymbol v_i) \right] \\
= & {\int \limits_0^1 \sum\limits_{j,k} \left[ \frac{\partial U}{\partial v^j}(\boldsymbol v_i +\xi \Delta \boldsymbol v_{i+1/2}) -\frac{\partial U}{\partial v^j}(\boldsymbol v_i ) \right] \frac{\partial g^{L,j}}{\partial v^{k}}(\boldsymbol v_{i},\boldsymbol v_i +\xi \Delta \boldsymbol v_{i+1/2})  \Delta v^{k}_{i+1/2} d\xi} \\
= & {\int\limits_0^1 \int \limits_0^1 \sum\limits_{k,l} \left[ \sum\limits_j { \frac{\partial^2 U}{\partial v^j \partial v^l} (\boldsymbol v_i +\eta \xi \Delta \boldsymbol v_{i+1/2}) }  \frac{\partial g^{L,j}}{\partial v^{k}}(\boldsymbol v_{i},\boldsymbol v_i +\xi \Delta \boldsymbol v_{i+1/2}) \right] \Delta v^{k}_{i+1/2} \Delta v^{l}_{i+1/2} d\eta \, \xi d\xi} 
\end{split}
\end{equation}
\begin{equation}
\label{eq:2.8}  
\begin{split} 
& G^L(\boldsymbol v_{i},\boldsymbol v_i)- G^L(\boldsymbol v_{i-1},\boldsymbol v_{i})
- \sum\limits_j \frac{\partial U}{\partial v^j}(\boldsymbol v_i)  \left[ g^{L,j}(\boldsymbol v_{i},\boldsymbol v_{i}) -g^{L,j}(\boldsymbol v_{i-1},\boldsymbol v_i) \right] \\
= & {\int \limits_0^1 \sum\limits_{j,k} \left[ \frac{\partial U}{\partial v^j}(\boldsymbol v_i -\xi \Delta \boldsymbol v_{i-1/2}) -\frac{\partial U}{\partial v^j}(\boldsymbol v_i ) \right] \frac{\partial g^{L,j}}{\partial v^{k}}(\boldsymbol v_i -\xi \Delta \boldsymbol v_{i-1/2}, \boldsymbol v_{i})  \Delta v^{k}_{i-1/2} d\xi} \\
= - \int\limits_0^1 & \int \limits_0^1 \sum\limits_{k,l} \left[ \sum\limits_j { \frac{\partial^2 U}{\partial v^j \partial v^l} (\boldsymbol v_i -\eta \xi \Delta \boldsymbol v_{i-1/2}) }  \frac{\partial g^{L,j}}{\partial v^{k}}(\boldsymbol v_i -\xi \Delta \boldsymbol v_{i-1/2}, \boldsymbol v_{i}) \right] \Delta v^{k}_{i-1/2} \Delta v^{l}_{i-1/2} d\eta \, \xi d\xi 
\end{split}
\end{equation}
where $\Delta \boldsymbol v_{i+1/2} = \boldsymbol v_{i+1} -\boldsymbol v_{i}$. 

Thus, according to our assumption, the integrals in \eqref{eq:2.7}-\eqref{eq:2.8} do not change the sign over the integration interval, which means that the second and third terms are negative.

The second and third terms on the right side of \eqref{eq:2.6} are linear in $\Delta t$, and the first term is of second-order.
Therefore, we can choose the time step small enough that the second and third terms  dominate over the first term.
Consequently, the right-hand side of \eqref{eq:2.5} is non-positive if $\Delta t$ satisfies \eqref{eq:2.4}. This completes the proof of the theorem.
\end{proof}

The proper numerical entropy flux for the first-order HR scheme~\eqref{eq:1.1} can be written as follows 
\begin{equation}
\label{eq:2.9}
\begin{split}
  G^L_{i+1/2} = \frac{1}{2} \left( F(U^-) +F(U^+) -c_{i+1/2} \; (U^+ -U^-) \right) . 
\end{split}
\end{equation}

Multiplying  \eqref{eq:1.1} by $U'(\boldsymbol {v_i})$ and subtracting it from the left-hand side of \eqref{eq:10}, we obtain
\begin{equation}
\label{eq:2.10}  
\begin{split}  
  \quad U(\boldsymbol {\hat v_i})& \; - \; U(\boldsymbol {v_i}) \; + \; \frac{{\Delta t}}{{\Delta x}}\left[ {G_{i+1/2}^L - G_{i-1/2}^L} \right] 
  =  U(\boldsymbol {\hat v_i}) - U(\boldsymbol {v_i}) -  \langle U'(\boldsymbol {v_i}), \left( {{\boldsymbol {\hat v}_i} - \boldsymbol {v_i}} \right) \rangle  \; \\ 
&  + \; \frac{{\Delta t}}{{\Delta x}}\left[ G_{i+1/2}^L - G_{i-1/2}^L 
  -  \left\langle U'(\boldsymbol {v_i}), \left( \boldsymbol g_{i+1/2}^L - \boldsymbol g_{i-1/2}^L - \boldsymbol s^-_{i+1/2} +\boldsymbol  s^+_{i-1/2} \right) \right\rangle \right]  \\
& = \; U(\boldsymbol {\hat v_i}) - U(\boldsymbol {v_i}) - \langle U'(\boldsymbol {v_i}), \left( {{\boldsymbol {\hat v}_i} - \boldsymbol {v_i}} \right) \rangle  \; + \; \frac{{\Delta t}}{{\Delta x}}\left[\Delta G_{i+1/2}^{L,-} - \Delta G_{i-1/2}^{L,+} \right] , 
\end{split}
\end{equation}
where $\Delta G_{i+1/2}^{L,\pm}$ are defined as
\begin{equation}
\label{eq:2.11}
\begin{split} 
& \qquad \Delta G_{i+1/2}^{L,\pm} = G_{i+1/2}^L(U^-_{i+1/2},U^+_{i+1/2}) -F(U_i) \\
  - & \left\langle U'(\boldsymbol {v_i}), \left( \boldsymbol g_{i+1/2}^L(\boldsymbol v^-_{i+1/2},\boldsymbol v^+_{i+1/2}) -  \boldsymbol s^{\pm}_{i+1/2} -\boldsymbol  f(\boldsymbol u_i) \right) \right\rangle . 
\end{split}
\end{equation}
or substituting the values of the corresponding functions, $\Delta G_{i+1/2}^{L,\pm}$ can be represented as
\begin{equation}
\label{eq:2.12}
\begin{split} 
 \Delta G_{i+1/2}^{L,-} & =  \; \frac{1}{2} \left\lbrace \frac{u_{i+1}}{2} y^+_{i+1/2} (u_{i+1} -u_i)^2 + g \left( u_{i+1} y^+_{i+1/2} - u_i y_i \right) (z_{i+1/2} - z_i ) \right.  \\
& + g \left. \left[ \frac{u_i}{2} ( y^-_{i+1/2} -y_{i})^2 + ( y^+_{i+1/2} -y_{i}) \left( u_{i+1} y^+_{i+1/2} -\frac{u_i}{2} ( y^+_{i+1/2} +y_i ) \right) \right] \right. \\
& - c_{i+1/2} \left. \left[ \frac{1}{2} y^+_{i+1/2} (u_{i+1} -u_i)^2 +\frac{g}{2} ( y^+_{i+1/2} -y_{i})^2 - \frac{g}{2} ( y^-_{i+1/2} -y_{i})^2 \right. \right. \\
& \qquad \qquad + \left. \left. g ( y^+_{i+1/2} -y^-_{i+1/2}) (z_{i+1/2} -z_i) \right] \right\rbrace ,
\end{split}
\end{equation}

\begin{equation}
\label{eq:2.13}
\begin{split} 
 \Delta G_{i-1/2}^{L,+} & =  \; \frac{1}{2} \left\lbrace  \frac{u_{i-1}}{2} y^-_{i-1/2} (u_{i-1} -u_i)^2 - g \left( u_{i-1} y^-_{i-1/2} - u_i y_i \right) (z_i - z_{i-1/2} ) \right.  \\
& + g \left. \left[ \, \frac{u_i}{2} ( y^+_{i-1/2} -y_{i})^2 + ( y^-_{i-1/2} -y_{i}) \left( u_{i-1} y^-_{i-1/2} -\frac{u_i}{2} ( y^-_{i-1/2} +y_i ) \right) \right] \right. \\
& + c_{i-1/2} \left. \left[ \, \frac{1}{2} y^-_{i-1/2} (u_{i-1} -u_i)^2 -\frac{g}{2} ( y^+_{i-1/2} -y_{i})^2 + \frac{g}{2} ( y^-_{i-1/2} -y_{i})^2 \right. \right. \\
& \qquad \qquad + \left. \left. g ( y^+_{i-1/2} -y^-_{i-1/2}) (z_i -z_{i-1/2}) \right] \right\rbrace .
\end{split}
\end{equation}

We consider two cases: (1) homogeneous and (2) inhomogeneous shallow water equations.

\begin{theorem}
\label{th:2.2}
 Assume that $c_{i+1/2}$ and $\Delta t$  satisfy the inequalities   
\begin{equation}
\label{eq:2.14}
\begin{split} 
c_{i+1/2} \ge \frac{1}{2} \max \left( -u_{i+1} -u_i +\sqrt{(u_{i+1}-u_i)^2 + 4g y_i }, u_{i+1} +u_i +\sqrt{(u_{i+1}-u_i)^2 + 4g y_{i+1} } \right), 
\end{split}
\end{equation}
\begin{equation}
\label{eq:2.15} 
\begin{split} 
\Delta t & \mathop {\max }\limits_{\boldsymbol v \in ({\boldsymbol  v_i},\boldsymbol {\hat v_i})} \lambda \left( U''(\boldsymbol v) \right) \;\langle {\left( {\boldsymbol g_{i+1/2}^L - \boldsymbol g_{i-1/2}^L} \right)}, {\left( {\boldsymbol g_{i+1/2}^L - \boldsymbol g_{i-1/2}^L} \right)} \rangle \\
 & \le 2\Delta x \left[ \langle U'(\boldsymbol {v_i}), \; \left( \boldsymbol g_{i+1/2}^L - \boldsymbol g_{i-1/2}^L  \right) \rangle - G_{i+1/2}^L + G_{i-1/2}^L \right]. 
\end{split}		  		
\end{equation}
where $\lambda \left( U''(\boldsymbol v) \right) = \frac{1}{2y} (u^2 + g y + \sqrt{(u^2+g y)^2 -4 g y}$.

Then for homogeneous shallow water equations, the fully discrete HR scheme \eqref{eq:1.1}-\eqref{eq:1.2} satisfies the discrete cell entropy inequality 
\begin{equation}
\label{eq:2.16}  
U(\boldsymbol {\hat v_i}) - U(\boldsymbol {v_i}) + \frac{{\Delta t}}{{\Delta x}}\left[ {G_{i+1/2}^L(\boldsymbol {v_i},\boldsymbol v_{i+1}) - G_{i-1/2}^L( \boldsymbol {v_{i-1}},\boldsymbol {v_i})} \right] \le 0.		
\end{equation}
where $G_{i+1/2}^L$ is the proper numerical entropy flux~\eqref{eq:2.9}.
\end{theorem}

\begin{proof}
For homogeneous shallow water equations, we have that $\boldsymbol v^-_{i-1/2} = \boldsymbol v_{i-1}, \boldsymbol v^-_{i+1/2} = \boldsymbol v^+_{i-1/2} = \boldsymbol v_{i}$ and $\boldsymbol v^+_{i+1/2} = \boldsymbol v_{i+1}$. Then $\Delta G_{i+1/2}^{L,\pm}$ in \eqref{eq:2.12}-\eqref{eq:2.13} take the form
\begin{equation}
\label{eq:2.17}
\begin{split} 
 \Delta G_{i+1/2}^{L,-} & =  \; \frac{1}{2} \left\lbrace \frac{u_{i+1}}{2} y_{i+1} (u_{i+1} -u_i)^2  + g ( y_{i+1} -y_{i}) \left( u_{i+1} y_{i+1} -\frac{u_i}{2} ( y_{i+1} +y_i ) \right) \right. \\
& - c_{i+1/2} \left. \left[ \frac{y_{i+1}}{2}  (u_{i+1} -u_i)^2 +\frac{g}{2} ( y_{i+1} -y_{i})^2  \right] 
\right\rbrace ,
\end{split}
\end{equation}
\begin{equation}
\label{eq:2.18}
\begin{split} 
 \Delta G_{i-1/2}^{L,+} & =  \; \frac{1}{2} \left\lbrace  \frac{u_{i-1}}{2} y_{i-1} (u_{i-1} -u_i)^2 
 + g  ( y_{i-1} -y_{i}) \left( u_{i-1} y_{i-1} -\frac{u_i}{2} ( y_{i-1} +y_i ) \right) \right. \\
& + c_{i-1/2} \left. \left[ \, \frac{1}{2} y_{i-1} (u_{i-1} -u_i)^2 + \frac{g}{2} ( y_{i-1} -y_{i})^2 \right] 
\right\rbrace .
\end{split}
\end{equation}

Note that the multipliers in the square brackets for $c_{i\pm1/2}$ in \eqref{eq:2.17}-\eqref{eq:2.18} are non-negative. It is easy to check that when they are zero, then $\Delta G^{L,\pm}_{i\mp1/2}$ are also zero. Let us show that $c_{i\pm1/2}$ can be chosen so that $\Delta G_{i+1/2}^{L,-}$ and $\Delta G_{i-1/2}^{L,+}$ are non-positive and non-negative, respectively. Indeed, we rewrite \eqref{eq:2.17} as
\begin{equation}
\label{eq:2.19}
\begin{split} 
 \Delta G_{i+1/2}^{L,-}  = & \; \frac{1}{2} \left\lbrace \frac{u_{i+1}}{2} y_{i+1} (u_{i+1} -u_i)^2  + g \, y_{i+1}( y_{i+1} -y_{i}) \left( u_{i+1} -u_i \right)  \right. \\
 + & g \, \frac{u_i}{2} \left(  y_{i+1} -y_i  \right)^2 - c_{i+1/2} \left. \left[ \frac{y_{i+1}}{2}  (u_{i+1} -u_i)^2 +\frac{g}{2} ( y_{i+1} -y_{i})^2  \right] 
\right\rbrace \\
 = &  \; \frac{1}{2} \left\lbrace \frac{y_{i+1}}{2}  (u_{i+1} -c_{i+1/2} ) (u_{i+1} -u_i)^2  + g \, y_{i+1}( y_{i+1} -y_{i}) \left( u_{i+1} -u_i \right)  \right. \\
+ & \left. \frac{g}{2} (u_i - c_{i+1/2}) \left(  y_{i+1} -y_i  \right)^2   
\right\rbrace .
\end{split}
\end{equation}
$\Delta G_{i+1/2}^{L,-}$ is a quadratic form with respect to $(u_{i+1} -u_i)$ and $( y_{i+1} -y_{i})$. For it to be non-positive, it is sufficient for the leading coefficient and its discriminant to be non-positive. Then, $c_{i+1/2}$ should satisfy the following inequalities
\begin{equation}
\label{eq:2.20}
\begin{split}   
&  c_{i+1/2} -u_{i+1} > 0, \\
  c_{i+1/2}^2 -(u_{i+1} & +u_i) \,  c_{i+1/2} + u_{i+1} \, u_i -g y_{i+1} \ge 0 .  
\end{split}
\end{equation}

It is clear that inequalities \eqref{eq:2.20} hold for $c_{i+1/2} \ge \frac{1}{2} \left( u_{i+1} +u_i +\sqrt{(u_{i+1}-u_i)^2 + 4g y_{i+1}} \right)$. Similarly, it can be shown that $\Delta G_{i-1/2}^{L,+}$ is non-negative for $c_{i-1/2} \ge \frac{1}{2} \left( -u_{i-1} -u_i \right. \\ \left. +\sqrt{(u_{i-1}-u_i)^2 + 4g y_{i-1} }\right)$.

We rewrite the discrete cell entropy inequality \eqref{eq:2.10} in the form
\begin{equation}
\label{eq:2.21}  
\begin{split}  
&  \quad U(\boldsymbol {\hat v_i}) \; - \; U(\boldsymbol {v_i}) \; + \; \frac{{\Delta t}}{{\Delta x}}\left[ {G_{i+1/2}^L - G_{i-1/2}^L} \right] 
   \\ 
& = \; U(\boldsymbol {\hat v_i}) - U(\boldsymbol {v_i}) - \langle U'(\boldsymbol {v_i}), \left( {{\boldsymbol {\hat v}_i} - \boldsymbol {v_i}} \right) \rangle  \; + \; \frac{{\Delta t}}{{\Delta x}}\left[\Delta G_{i+1/2}^{L,-} - \Delta G_{i-1/2}^{L,+} \right]  \\
& = \frac{1}{2}{\left( {\frac{\Delta t}{\Delta x}} \right)^2}  \left\langle U''(\boldsymbol s )  \left({\boldsymbol g_{i+1/2}^L - \boldsymbol  g_{i-1/2}^L} \right),\left({\boldsymbol g_{i+1/2}^L - \boldsymbol  g_{i-1/2}^L} \right) \right\rangle   + \frac{{\Delta t}}{{\Delta x}}\left[\Delta G_{i+1/2}^{L,-} - \Delta G_{i-1/2}^{L,+} \right] .
\end{split}
\end{equation}

Thus, the non-positivity of \eqref{eq:2.21} can be achieved by choosing a sufficiently small $\Delta t$ so that the second term dominates over the first non-negative term on the right-hand side of \eqref{eq:2.21}. This completes the proof of the theorem.
\end{proof}

\begin{theorem}
\label{th:2.3}
 Suppose that $c_{i+1/2}$ and $\Delta t$  satisfy the inequalities   
\begin{equation}
\label{eq:2.22}
\begin{split}  
  c_{i+1/2} & \ge \max \left( \frac{a^-_{i} u_{i+1} +b^-_{i} u_i +\sqrt{(a^-_{i} u_{i+1} -b^-_{i} u_i)^2 +4 g \, a^-_{i} y^+_{i+1/2} (w^+_{i+1/2} -w_i)^2 }}{2a^-_{i}}, \right. \\
& \left. \frac{-a^+_{i+1} u_{i} -b^+_{i+1} u_{i+1} +\sqrt{(a^+_{i+1} u_{i} -b^+_{i+1} u_{i+1})^2 +4 g \, a^+_{i+1} y^-_{i+1/2} (w^-_{i+1/2} -w_{i+1})^2 }}{2a^+_{i+1}}  \right). 
\end{split}
\end{equation}
\begin{equation}
\label{eq:2.23} 
\begin{split} 
& \Delta t \mathop {\max }\limits_{\boldsymbol v \in ({\boldsymbol  v_i},\boldsymbol {\hat v_i})} \lambda \left( U''(\boldsymbol v) \right) \;\langle {\left( {\boldsymbol g_{i+1/2}^L - \boldsymbol g_{i-1/2}^L} - \boldsymbol s^-_{i+1/2} +\boldsymbol  s^+_{i-1/2}\right)}, {\left( {\boldsymbol g_{i+1/2}^L - \boldsymbol g_{i-1/2}^L} - \boldsymbol s^-_{i+1/2} +\boldsymbol  s^+_{i-1/2} \right)} \rangle \\
&  \le 2\Delta x \left[ \langle U'(\boldsymbol {v_i}), \; \left( \boldsymbol g_{i+1/2}^L - \boldsymbol g_{i-1/2}^L  - \boldsymbol s^-_{i+1/2} +\boldsymbol  s^+_{i-1/2} \right) \rangle - G_{i+1/2}^L + G_{i-1/2}^L \right]. 
\end{split}		  		
\end{equation}
where  
\begin{equation}
\label{eq:2.24}
\begin{split}  
& \lambda \left( U''(\boldsymbol v) \right) = \frac{1}{2y} (u^2 + g y + \sqrt{(u^2+g y)^2 -4 g y} \\
& a^{\mp}_{i} = \pm( y^+_{i\pm1/2} -y^-_{i\pm1/2}) \left(  w^\pm_{i+1/2}+w^-_{i\pm1/2}  -2 w_{i} \right), \\
& b^\mp_{i} =  \pm( y^-_{i\pm1/2} -y_{i})^2 \mp ( y^+_{i\pm1/2} -y_{i})^2 + 2 (  y^\pm_{i\pm1/2} - y_i)( w^\pm_{i\pm1/2} -w_i )  .
\end{split}
\end{equation}

Then for inhomogeneous shallow water equations, the fully discrete HR scheme \eqref{eq:1.1}-\eqref{eq:1.2} satisfies the discrete cell entropy inequality 
\begin{equation}
\label{eq:2.25}  
U(\boldsymbol {\hat v_i}) - U(\boldsymbol {v_i}) + \frac{{\Delta t}}{{\Delta x}}\left[ {G_{i+1/2}^L(\boldsymbol {v^-_{i+1/2}},\boldsymbol v^+_{i+1/2}) - G_{i-1/2}^L( \boldsymbol {v^-_{i-1/2}},\boldsymbol {v^+_{i-1/2}})} \right] \le 0.		
\end{equation}
where $G_{i+1/2}^L$ is the proper numerical entropy flux~\eqref{eq:2.9}.
\end{theorem}

\begin{proof}
We rewrite $\Delta G_{i+1/2}^{L,\pm}$ in \eqref{eq:2.12}-\eqref{eq:2.13} as follows
\begin{equation}
\label{eq:2.26}
\begin{split} 
 \Delta & G_{i+1/2}^{L,-} =  \; \frac{1}{2} \left\lbrace \frac{u_{i+1}}{2} y^+_{i+1/2} (u_{i+1} -u_i)^2 + g \, y^+_{i+1/2} ( u_{i+1}  - u_i ) (w^+_{i+1/2} - w_i ) \right.  \\
& + g \left. \left[ \frac{u_i}{2} ( y^-_{i+1/2} -y_{i})^2 - \frac{u_i}{2} ( y^+_{i+1/2} -y_{i})^2 + u_{i} (  y^+_{i+1/2} - y_i)( w^+_{i+1/2} -w_i )  \right] \right. \\
& - c_{i+1/2} \left. \left[ \frac{1}{2} y^+_{i+1/2} (u_{i+1} -u_i)^2 + g ( y^+_{i+1/2} -y^-_{i+1/2}) \left( \frac{1}{2} (w^+_{i+1/2}+w^-_{i+1/2})  -w_{i} \right) \right] \right\rbrace ,
\end{split}
\end{equation}
\begin{equation}
\label{eq:2.27}
\begin{split} 
 \Delta  & G_{i-1/2}^{L,+}=  \; \frac{1}{2} \left\lbrace  \frac{u_{i-1}}{2} y^-_{i-1/2} (u_{i-1} -u_i)^2 - g \, y^-_{i-1/2}( u_{i-1}  - u_i ) (w_i - w^-_{i-1/2} ) \right.  \\
& + g \left. \left[ \, \frac{u_i}{2} ( y^+_{i-1/2} -y_{i})^2 - \frac{u_i}{2}  ( y^-_{i-1/2} -y_{i})^2 + u_i( y^-_{i-1/2} -y_i) ( w^-_{i-1/2} -w_i ) \right] \right. \\
& + c_{i-1/2} \left. \left[ \, \frac{1}{2} y^-_{i-1/2} (u_{i-1} -u_i)^2 + g ( y^+_{i-1/2} -y^-_{i-1/2}) \left(w_i -\frac{1}{2} (w^+_{i-1/2} + w^-_{i-1/2} )  \right) \right] \right\rbrace .
\end{split}
\end{equation}

Let us show that the coefficients in the square brackets at $c_{i\pm1/2}$ in \eqref{eq:2.26}-\eqref{eq:2.27} are non-negative. Consider the following cases:

(i) In the fully wet case, $\min(w_i,w_{i+1}) > \max(z_i,z_{i+1})$. According to \eqref{eq:1.4}-\eqref{eq:1.5}, we have $w^-_{i+1/2}=w_i$ and $w^+_{i+1/2}=w_{i+1}$. 

Hence, if $y^+_{i+1/2} \ge y^-_{i+1/2}$, then $w^+_{i+1/2} \ge w^-_{i+1/2}$, and $( y^+_{i+1/2} -y^-_{i+1/2}) \left( \frac{1}{2} (w^+_{i+1/2}+w^-_{i+1/2})  -w_{i} \right)$ $ \ge 0$. Otherwise, if $y^+_{i+1/2} < y^-_{i+1/2}$, then also $w^+_{i+1/2} < w^-_{i+1/2}$, and the required inequality holds. 

(ii) In the partially wet case $ \min(w_i, w_{i+1}) \le \max(z_i, z_{i+1})$. Depending on which bottom is higher, right or left, we consider two subcases.

Let $z_i \ge z_{i+1}$. Then $z_{i+1/2} = w_{i+1}, \; y^-_{i+1/2} = y_i, \; y^+_{i+1/2}= 0$, and $w^-_{i+1/2},w^+_{i+1/2} \le w_i$. Therefore, $( y^+_{i+1/2} -y^-_{i+1/2}) \left( \frac{1}{2} (w^+_{i+1/2}+w^-_{i+1/2})  -w_{i} \right) \ge 0$.

If $z_i < z_{i+1}$, then
$z_{i+1/2} = w_i, \, y^-_{i+1/2} = 0, \, y^+_{i+1/2}= y_{i+1}$, and $w^-_{i+1/2} = w_i, \, w^-_{i+1/2} > w_i$. Thus, we again obtain the required inequality.

The non-negativity of the terms in square brackets at $c_{i-1/2}$ is proved similarly.

It is easy to check that when they are zero, then $\Delta G^{L,\pm}_{i\mp1/2}$ are also zero. Let us show that $c_{i\pm1/2}$ can be chosen so that $\Delta G_{i+1/2}^{L,-}$ and $\Delta G_{i-1/2}^{L,+}$ are non-positive and non-negative, respectively. Indeed, we consider $\Delta G_{i+1/2}^{L,-}$ and $\Delta G_{i-1/2}^{L,+}$ 
as quadratic equations with respect to $(u_{i+1} -u_i)$ and $(u_{i-1} -u_i)$, respectively.

$\Delta G_{i+1/2}^{L,-}$ will be non-positive if its leading coefficient and discriminant are non-positive, i.e.
\begin{equation}
\label{eq:2.28}
\begin{split} 
&  c_{i+1/2} > u_{i+1}, \\
-a^-_{i} c_{i+1/2}^2 +(a^-_{i} u_{i+1} +b^-_{i} u_i) c_{i+1/2} & - u_{i+1} u_i b^-_{i} +g \, y^+_{i+1/2} (w^+_{i+1/2} -w_i)^2 \le 0,   
\end{split}
\end{equation}
where 
\begin{equation}
\label{eq:2.29}
\begin{split}  
& a^-_{i} = ( y^+_{i+1/2} -y^-_{i+1/2}) \left( w^+_{i+1/2}+w^-_{i+1/2}  - 2 w_{i} \right), \\
& b^-_{i} =  ( y^-_{i+1/2} -y_{i})^2 - ( y^+_{i+1/2} -y_{i})^2 + 2 (  y^+_{i+1/2} - y_i)( w^+_{i+1/2} -w_i )  .
\end{split}
\end{equation}

It is clear that the inequalities \eqref{eq:2.28} hold for
\begin{equation}
\label{eq:2.30}
  c_{i+1/2} \ge \frac{a^-_{i} u_{i+1} +b^-_{i} u_i +\sqrt{(a^-_{i} u_{i+1} -b^-_{i} u_i)^2 +4 g \, a^-_{i} y^+_{i+1/2} (w^+_{i+1/2} -w_i)^2 }}{2a^-_{i}}.
\end{equation}

Similarly, it is proved that $\Delta G_{i-1/2}^{L,+}$ is non-negative for 
\begin{equation}
\label{eq:2.31}
  c_{i-1/2} \ge \frac{-a^+_{i} u_{i-1} -b^+_{i} u_i +\sqrt{(a^+_{i} u_{i-1} -b^+_{i} u_i)^2 +4 g \, a^+_{i} y^-_{i-1/2} (w^-_{i-1/2} -w_i)^2 }}{2a^+_{i}},
\end{equation}
where
\begin{equation}
\label{eq:2.32}
\begin{split}  
& a^+_{i} = ( y^+_{i-1/2} -y^-_{i-1/2}) \left( 2 w_{i} -  w^+_{i-1/2} - w^-_{i-1/2}  \right), \\
& b^+_{i} =  ( y^+_{i-1/2} -y_{i})^2 - ( y^-_{i-1/2} -y_{i})^2 + 2 (  y^-_{i-1/2} - y_i)( w^-_{i-1/2} -w_i )  .
\end{split}
\end{equation}

Thus, returning to the discrete entropy inequality \eqref{eq:2.21}, we can choose $\Delta t$ so that the non-positive terms in square brackets dominate over the non-negative first term on the right-hand side of \eqref{eq:2.21}.
This concludes the proof of the theorem.
\end{proof}


\section{Finding Flux Limiters}
\label{sec:3}
 
The system of equations \eqref{eq:8}-\eqref{eq:9} is nonlinear if we consider $\boldsymbol {\alpha}$ as a function of $\boldsymbol {\hat v}$, and it can be written in the form 
\begin{equation}
\label{eq:3.1} 
  \boldsymbol {\hat v} -\Delta t \, P \boldsymbol {\hat v} = \boldsymbol {\hat v^L},  
\end{equation}  
where mapping $P: R^{N}\times R^N \rightarrow R^{N} \times R^N$ is defined by  $P_i \boldsymbol {\hat v} = {\left[ {\alpha _{i -1/2} (\boldsymbol {\hat v})} \boldsymbol g_{i-1/2}^{AD} - {\alpha _{i+1/2} (\boldsymbol {\hat v})} \,\boldsymbol g_{i+1/2}^{AD} \right]}$ $ / {\Delta x}$, $\boldsymbol g_{i+1/2}^{AD} = \boldsymbol  g_{i+1/2}^H - \boldsymbol g_{i+1/2}^L$, and $ \hat{\boldsymbol  v}^L_i = \boldsymbol v_i - {\Delta t} / {\Delta x} \left( {\boldsymbol g_{i+1/2}^L - \boldsymbol g_{i-1/2}^L} \right) + \boldsymbol s_i$. Let $ O_0 = \bar O(\boldsymbol{ \hat{v}^L},\delta) $ be a closed ball with center at $ \boldsymbol {\hat v}^L $ and radius $ \delta > 0 $. Furthermore, we define a mapping $ S\boldsymbol v  : O_0 \rightarrow R^{N} \times R^N $ as
\begin{equation} 
\label{eq:3.2}
  S \boldsymbol v = \boldsymbol v - \Delta t P \boldsymbol v.
\end{equation}
 
Let us show that for sufficiently small $ \Delta t $ the system of equations \eqref{eq:3.1} is uniquely solvable in a neighborhood of the first-order HR solution of \eqref{eq:1.1}.
 
\begin{theorem} 
\label{th:3.1} 
  Assume that
\begin{equation}
\label{eq:3.3}
  \Vert P(\boldsymbol w) -P(\boldsymbol v) \Vert \le M \Vert \boldsymbol w -\boldsymbol v \Vert, \qquad \forall \boldsymbol v,\boldsymbol w \in O_0.
\end{equation} 

If $\Delta t$ satisfies
\begin{equation} 
\label{eq:3.4} 
\Delta t < \delta (\Vert P \boldsymbol {\hat v}^L \Vert + \delta M)^{-1},
\end{equation}
then the system of equations \eqref{eq:3.1} has a unique solution in $O_0$.   
\end{theorem} 
\begin{proof} Our proof mimics the proof of theorem 5.1.6~\cite[p.122]{ortega}.
For fixed $\boldsymbol d \in O \left( S \boldsymbol {\hat v}^L,\varepsilon \right)$,  we define the mapping $T: O_0 \rightarrow R^{N} \times R^N $ by 
\[ T \boldsymbol y =  \Delta t P \boldsymbol v + \boldsymbol d = 
\boldsymbol v - \left[ S \boldsymbol v - \boldsymbol d \right]. \]  
 Then, $ S \boldsymbol v = \boldsymbol d$ has a unique solution in $O_0$ if and only if $T$ has a unique fixed point. 
For any $\boldsymbol v, \boldsymbol w \in O_0$ 
\begin{equation}
\label{eq:3.5} 
\Vert T \boldsymbol v -T \boldsymbol w \Vert = \Delta t\Vert   P \boldsymbol v -P \boldsymbol w  \Vert 
\leq  \Delta t M \left\| \boldsymbol v - \boldsymbol w \right\| 
\end{equation}
and $S$ is contractive on $O_0$ if $\Delta t M < 1$. Moreover, for any $ \boldsymbol v \in O_0$,
\begin{equation}
\label{eq:3.6} 
\begin{split}
\Vert T \boldsymbol v -\boldsymbol {\hat v}^L \Vert  \leq \Vert T \boldsymbol v -T \boldsymbol {\hat v}^L \Vert +\Vert T \boldsymbol {\hat v}^L -\boldsymbol {\hat v}^L \Vert \leq \Delta t M  \Vert \boldsymbol v -\boldsymbol {\hat v}^L \Vert +\Vert S \boldsymbol v - \boldsymbol d \Vert 
  \leq \Delta t M \delta + \varepsilon. 
\end{split} 
\end{equation} 
For $\varepsilon = \delta (1 - \Delta t M)$, the expression on the right-hand side of \eqref{eq:3.6} equal to $\delta$. 
Hence, $T$ maps $O_0$ into $O_0$, and for any $\boldsymbol d \in O \left( S \boldsymbol {\hat v}^L,\varepsilon \right)$ the equation $S \boldsymbol v = \boldsymbol d$ has a unique solution in $O_0$.

Finally, we have that $\boldsymbol {\hat v}^L \in O \left( S \boldsymbol {\hat v}^L,\varepsilon \right)$ if 
\begin{equation}
\label{eq:3.7} 
\Vert S \boldsymbol {\hat v}^L - \boldsymbol {\hat v}^L \Vert =  \Delta t \Vert P \boldsymbol {\hat v}^L \Vert < \varepsilon.
\end{equation}
Combining all restrictions on $\Delta t$, we obtain that the nonlinear system of equation \eqref{eq:3.1} has a unique solution if $\Delta t$ satisfies \eqref{eq:3.4}. 
\end{proof}

\begin{remark}
\label{rm:3.1}
  Note that the mapping $S$ in \eqref{eq:3.2} is contractive in the vicinity of a low-order solution with the HR scheme \eqref{eq:1.1}-\eqref{eq:1.7} if the mapping $P$ in \eqref{eq:3.1} is Lipschitz-continuous. 
\end{remark}

Our goal is to find the maximum values of the flux limiters $\boldsymbol \alpha \in U_{ad} = \left\{ \boldsymbol{\alpha} \vert  \; \rm{0} \le \alpha_{i+1/2} \le \rm{1} \right\}$, for which the numerical solution of the hybrid scheme \eqref{eq:8}-\eqref{eq:9} satisfies the constraints \eqref{eq:1.9}-\eqref{eq:1.10} and the discrete cell entropy inequality \eqref{eq:2.25}. Then finding the flux limiters can be considered as the following optimization problem
\begin{equation}
\label{eq:3.8} 
\Im (\boldsymbol{\alpha} ) = \sum\limits_i {\alpha _{i+1/2}} \to \mathop {\max }\limits_{\boldsymbol{\alpha}  \in U_{ad}}  			\end{equation}
subject to 
\begin{multline}
\label{eq:3.9}  
 \frac{{\Delta x}}{{\Delta t}}\left( \underline  w_i -  w_i \right) + \frac{c_{i+1/2} -u_{i+1}}{2} (w_i - w^+_{i+1/2}) + \frac{c_{i-1/2} +u_{i-1}}{2} (w_i - w^-_{i-1/2})		
  \le  - \alpha_{i+1/2}\; g_{i+1/2}^{AD,y} \\ + \alpha_{i-1/2}\; g_{i-1/2}^{AD,y} \le 			
     \frac{{\Delta x}}{{\Delta t}}\left( {\overline w_i - w_i} \right) + \frac{c_{i+1/2} -u_{i+1}}{2} (w_i - w^+_{i+1/2}) + \frac{c_{i-1/2} +u_{i-1}}{2} (w_i - w^-_{i-1/2}) , 
\end{multline} 

\begin{multline}
\label{eq:3.10}  
 \frac{{\Delta x}}{{\Delta t}}\left( \underline  q_i -  q_i \right) + \frac{c_{i+1/2} -u_{i+1}}{2} (q_i - q^+_{i+1/2}) + \frac{c_{i-1/2} +u_{i-1}}{2} (q_i - q^-_{i-1/2})		
  \le  - \alpha_{i+1/2}\; g_{i+1/2}^{AD,q} \\ + \alpha_{i-1/2}\; g_{i-1/2}^{AD,q} \le 			
     \frac{{\Delta x}}{{\Delta t}}\left( {\overline q_i - q_i} \right) + \frac{c_{i+1/2} -u_{i+1}}{2} (q_i - q^+_{i+1/2}) + \frac{c_{i-1/2} +u_{i-1}}{2} (q_i - q^-_{i-1/2}) , 
\end{multline} 

\begin{multline}
\label{eq:3.11}  
 \frac{\Delta x}{\Delta t}\left[ U(\boldsymbol {\hat v}_i) - U(\boldsymbol {v_i}) - \left\langle U'(\boldsymbol v_i), \left( {\boldsymbol {\hat v}_i - \boldsymbol  v_i} \right) \right\rangle \right] -  \left\langle U'(\boldsymbol {v_i}), ( \boldsymbol g_{i+1/2}^L - \boldsymbol g_{i-1/2}^L - \boldsymbol s^-_{i+1/2} +\boldsymbol  s^+_{i-1/2} ) \right\rangle \\
 + G_{i+1/2}^L  - G_{i-1/2}^{L} 				
  \le \alpha _{i+1/2} \left( \langle U'(\boldsymbol v_i),\boldsymbol  g_{i+1/2}^{AD} \rangle - G_{i+1/2}^{AD} \right) -  \alpha _{i-1/2} \;( \langle U'(\boldsymbol v_i), \boldsymbol  g_{i-1/2}^{AD} \rangle - G_{i-1/2}^{AD} ),
\end{multline}

\begin{equation}
\label{eq:3.12}
 \frac{\Delta x}{\Delta t}\left( {\boldsymbol {\hat v}_i - \boldsymbol  v_i} \right) + \boldsymbol g_{i+1/2}^L - \boldsymbol  g_{i-1/2}^L -\boldsymbol s_{i+1/2}^- + \boldsymbol  s_{i-1/2}^+ + \alpha_{i+1/2} \, \boldsymbol  g_{i+1/2}^{AD} - \alpha_{i-1/2} \boldsymbol \,g_{i-1/2}^{AD} = 0,					
\end{equation}
where $\underline w_i = \min \left(w_i,w^-_{i-1/2},w^+_{i+1/2} \right)$, $\overline w_i = \max \left(w_i,w^-_{i-1/2},w^+_{i+1/2} \right)$, and $G^{AD}_{i+1/2} = G^{H}_{i+1/2} -G^{L}_{i+1/2} $.

Due to constraints \eqref{eq:3.11} the optimization problem \eqref{eq:3.8}-\eqref{eq:3.12} is nonlinear. Consequently, finding a numerical entropy solution of shallow water equations with variable bottom topography \eqref{eq:1}-\eqref{eq:2}  in one time step can be represented as the following iterative process.

\begin{description}
\item[\it Step 1.]  Initialize positive numbers $\delta $, $\epsilon_1 $, and  $\epsilon_2 $. Set $p=0$, 
 $ {\boldsymbol{\hat v}}^0 = \boldsymbol{v}$, $ \boldsymbol{\alpha} ^0 = 0$.
\item[\it Step 2.] Find $ \boldsymbol {\alpha}^{p + 1}$ as a solution to the following linear programming problem
\begin{equation}
\label{eq:3.13}
\Im (\boldsymbol {\alpha} ) = \sum\limits_i \alpha _{i +1/2}^{p + 1} \to \mathop {\max }\limits_{\boldsymbol {\alpha} ^{p + 1} \in U_{ad}} 					\end{equation}
subject to
\begin{multline}
\label{eq:3.14}  
 \frac{{\Delta x}}{{\Delta t}}\left( \underline  w_i -  w_i \right) + \frac{c_{i+1/2} -u_{i+1}}{2} (w_i - w^+_{i+1/2}) + \frac{c_{i-1/2} +u_{i-1}}{2} (w_i - w^-_{i-1/2})		
  \le  - \alpha_{i+1/2}^{p+1}\; g_{i+1/2}^{AD,y} \\ + \alpha_{i-1/2}^{p+1}\; g_{i-1/2}^{AD,y} \le 			
     \frac{{\Delta x}}{{\Delta t}}\left( {\overline w_i - w_i} \right) + \frac{c_{i+1/2} -u_{i+1}}{2} (w_i - w^+_{i+1/2}) + \frac{c_{i-1/2} +u_{i-1}}{2} (w_i - w^-_{i-1/2}) , 
\end{multline} 
\begin{multline}
\label{eq:3.15}  
 \frac{{\Delta x}}{{\Delta t}}\left( \underline  q_i -  q_i \right) + \frac{c_{i+1/2} -u_{i+1}}{2} (q_i - q^+_{i+1/2}) + \frac{c_{i-1/2} +u_{i-1}}{2} (q_i - q^-_{i-1/2})		
  \le  - \alpha_{i+1/2}^{p+1}\; g_{i+1/2}^{AD,q} \\ + \alpha_{i-1/2}^{p+1}\; g_{i-1/2}^{AD,q} \le 			
     \frac{{\Delta x}}{{\Delta t}}\left( {\overline q_i - q_i} \right) + \frac{c_{i+1/2} -u_{i+1}}{2} (q_i - q^+_{i+1/2}) + \frac{c_{i-1/2} +u_{i-1}}{2} (q_i - q^-_{i-1/2}) , 
\end{multline} 
\begin{multline}
\label{eq:3.16}  
 \frac{\Delta x}{\Delta t}\left[ U(\boldsymbol {\hat v}_i^p) - U(\boldsymbol {v_i}) - \left\langle U'(\boldsymbol v_i), \left( {\boldsymbol {\hat v}_i^p - \boldsymbol  v_i} \right) \right\rangle \right] -  \left\langle U'(\boldsymbol {v_i}), ( \boldsymbol g_{i+1/2}^L - \boldsymbol g_{i-1/2}^L - \boldsymbol s^-_{i+1/2} +\boldsymbol  s^+_{i-1/2} ) \right\rangle \\
 + G_{i+1/2}^L  - G_{i-1/2}^{L} 				
  \le \alpha _{i+1/2}^{p+1} \left( \langle U'(\boldsymbol v_i),\boldsymbol  g_{i+1/2}^{AD} \rangle - G_{i+1/2}^{AD} \right) -  \alpha _{i-1/2}^{p+1} \;( \langle U'(\boldsymbol v_i), \boldsymbol  g_{i-1/2}^{AD} \rangle - G_{i-1/2}^{AD} ),
\end{multline}
\item[\it Step 3.] For ${\boldsymbol {\alpha} ^{p + 1}}$ we find ${\boldsymbol {\hat v}^{p + 1}}$ from the system of linear equations
\begin{equation}
\label{eq:3.17}
 \frac{\Delta x}{\Delta t}\left( {\boldsymbol {\hat v}_i^{p+1} - \boldsymbol  v_i} \right) + \boldsymbol g_{i+1/2}^L - \boldsymbol  g_{i-1/2}^L -\boldsymbol s_{i+1/2}^- + \boldsymbol  s_{i-1/2}^+ + \alpha_{i+1/2}^{p+1} \, \boldsymbol  g_{i+1/2}^{AD} - \alpha_{i-1/2}^{p+1} \, \boldsymbol g_{i-1/2}^{AD} = 0,				\end{equation}
\item[\it Step 4.] Algorithm stop criterion
\begin{equation}
\label{eq:3.18}
 \frac{{\left| {\hat y_i^{p + 1} - \hat y_i^p} \right|}}{{\max \left( {\delta ,\left| {\hat y_i^{p + 1}} \right|} \right)}} < {\varepsilon _1},  \quad \frac{{\left| {\hat q_i^{p + 1} - \hat q_i^p} \right|}}{{\max \left( {\delta ,\left| {\hat q_i^{p + 1}} \right|} \right)}} < {\varepsilon _1},  \quad   \left| {\alpha _{i +1/2}^{p + 1} - \alpha _{i+1/2}^p} \right| < {\varepsilon _2}. 				
\end{equation}
If conditions \eqref{eq:3.18} hold, then set $\boldsymbol {\hat v} = {\boldsymbol {\hat v}^{p + 1}}$. Otherwise, set $p = p + 1$ and go to \textit{Step 2}.
\end{description}

\begin{remark}  
  It is clear that the linear programming problem \eqref{eq:3.13}-\eqref{eq:3.16} is solvable if $\Delta t$ satisfies inequalities \eqref{eq:1.8} and \eqref{eq:2.23}. It follows from the non-emptiness of the feasible set and the boundedness on $U_{ad}$ of the objective function $\Im (\boldsymbol {\alpha} )$. 
\end{remark}


\section{Approximate Solution to the Optimization Problem }
\label{sec:4}

Solving a linear programming problem is computationally expensive. So, at \textit{Step 2}, instead of solving the linear programming problem, it is reasonable to use its  computationally less expensive approximate solution. In this section our goal is to look for an approximate solution of the linear programming problem \eqref{eq:3.13} - \eqref{eq:3.16}.

First, we find a nontrivial $\boldsymbol {\alpha}  \in {U_{ad}}$ satisfying inequalities \eqref{eq:3.14}, which are rewritten in the form
\begin{equation}
\label{eq:4.1}
 - {\alpha _{i+1/2}} \; g_{i+1/2}^{AD,y} + {\alpha_{i-1/2}} \; g_{i-1/2}^{AD,y} \ge Q_i^{-,y} ,								
\end{equation}

\begin{equation}
\label{eq:4.2}
 - {\alpha_{i+1/2}} \; g_{i+1/2}^{AD,y} + {\alpha_{i-1/2}} \; g_{i-1/2}^{AD,y} \le Q_i^{ +,y} ,								
\end{equation}
where 
\begin{eqnarray*}
  && Q_i^ {+,y}  = \frac{\Delta x}{\Delta t} \left( \overline w_i - w_i \right) + \frac{c_{i+1/2} -u_{i+1}}{2} (w_i - w^+_{i+1/2}) + \frac{c_{i-1/2} +u_{i-1}}{2} (w_i - w^-_{i-1/2}) , \\
  && Q_i^ {-,y}  = \frac{\Delta x}{\Delta t}\left( \underline w_i - w_i \right) + \frac{c_{i+1/2} -u_{i+1}}{2} (w_i - w^+_{i+1/2}) + \frac{c_{i-1/2} +u_{i-1}}{2} (w_i - w^-_{i-1/2}).
\end{eqnarray*} 
Denote by $\alpha _i^{-,y} $ and $\alpha _i^{+,y} $ the maximum values of the components $\alpha$ for the negative and positive terms on the left-hand side of \eqref{eq:4.1}-\eqref{eq:4.2}, respectively. Then
\begin{eqnarray}
\label{eq:4.3}
 - {\alpha _{i+1/2}} \; g_{i+1/2}^{AD,y} + {\alpha_{i-1/2}}  \; g_{i-1/2}^{AD,y} \ge \alpha_i^{-,y} P_i^{-,y} ,\\							
 - {\alpha_{i+1/2}} \; g_{i+1/2}^{AD,y} + {\alpha_{i-1/2}}  \; g_{i-1/2}^{AD,y} \le \alpha_i^{+,y} P_i^{+,y} ,							
\label{eq:4.4}
\end{eqnarray}
where
\begin{eqnarray*}
  && P_i^{-,y}  = \min \left( {0, -  g_{i+1/2}^{AD,y}} \right) + \min \left( {0, g_{i-1/2}^{AD,y}} \right),		\\				
  && P_i^{+,y}  = \max \left( {0, -  g_{i+1/2}^{AD,y}} \right) + \max \left( {0, g_{i-1/2}^{AD,y}} \right) .
\end{eqnarray*} 
						
Each flux limiter ${\alpha _{i+1/2}}$ appears twice in \eqref{eq:4.3} and twice in \eqref{eq:4.4} with coefficients that differ only in sign. Substituting \eqref{eq:4.3}-\eqref{eq:4.4} into \eqref{eq:4.1}-\eqref{eq:4.2}, we obtain that ${\alpha_{i+1/2}}$ should not exceed 
\begin{equation}
\label{eq:4.5}
{\bar \alpha _{i+1/2}}^y = \left\{ {\begin{array}{*{20}{c}}
{\min (\alpha _i^{+,y} ,\alpha _{i+1}^{-,y} ) = \min (R_i^{+,y} ,R_{i + 1}^{-,y} ), \quad {\rm{    }} g_{i+1/2}^{AD,y} < 0}  ,\\
{\min (\alpha _i^{-,y} ,\alpha _{i+1}^{+,y} ) = \min (R_i^{-,y} ,R_{i+1}^{+,y} ), \quad {\rm{    }} g_{i+1/2}^{AD,y} > 0}  ,
\end{array}} \right. 					
\end{equation} 
where $R_i^{-,y}  = \min \left( 1,{\min(0,Q_i^{-,y}) \mathord{\left/
 {\vphantom {Q_i^{\pm,y} P_i^{\pm,y})}} \right.
 \kern-\nulldelimiterspace} {P_i^{-,y} }} \right)$ and $R_i^{+,y}  = \min \left(1,{\max ( 0,Q_i^{+,y} ) \mathord{\left/
 {\vphantom {Q_i^{\pm,y} P_i^{\pm,y})}} \right.
 \kern-\nulldelimiterspace} {P_i^{+,y} }}\right)$. 

Similarly, it is proved that inequalities \eqref{eq:3.15} hold for 
\begin{equation}
\label{eq:4.6}
{\bar \alpha _{i+1/2}}^q = \left\{ {\begin{array}{*{20}{c}}
{\min (\alpha _i^{+,q} ,\alpha _{i+1}^{-,q} ) = \min (R_i^{+,q} ,R_{i + 1}^{-,q} ), \quad {\rm{    }} g_{i+1/2}^{AD,q} < 0}  ,\\
{\min (\alpha _i^{-,q} ,\alpha _{i+1}^{+,q} ) = \min (R_i^{-,q} ,R_{i+1}^{+,q} ), \quad {\rm{    }} g_{i+1/2}^{AD,q} > 0}  ,
\end{array}} \right. 					
\end{equation} 
where $R_i^{-,q}  = \min \left( 1,{\min(0,Q_i^{-,q}) \mathord{\left/
 {\vphantom {Q_i^{\pm,q} P_i^{\pm,q})}} \right.
 \kern-\nulldelimiterspace} {P_i^{-,q} }} \right)$ and $R_i^{+,q}  = \min \left(1,{\max ( 0,Q_i^{+,q} ) \mathord{\left/
 {\vphantom {Q_i^{\pm,q} P_i^{\pm,q})}} \right.
 \kern-\nulldelimiterspace} {P_i^{+,q} }}\right)$ , 
\begin{eqnarray*}
  && Q_i^ {+,q}  = \frac{\Delta x}{\Delta t} \left( \overline q_i - q_i \right) + \frac{c_{i+1/2} -u_{i+1}}{2} (q_i - q^+_{i+1/2}) + \frac{c_{i-1/2} +u_{i-1}}{2} (q_i - q^-_{i-1/2}) , \\
  && Q_i^ {-,q}  = \frac{\Delta x}{\Delta t}\left( \underline q_i - q_i \right) + \frac{c_{i+1/2} -u_{i+1}}{2} (q_i - q^+_{i+1/2}) + \frac{c_{i-1/2} +u_{i-1}}{2} (q_i - q^-_{i-1/2}).
\end{eqnarray*} 
\begin{eqnarray*}
  && P_i^{-,q}  = \min \left( {0, -  g_{i+1/2}^{AD,q}} \right) + \min \left( {0, g_{i-1/2}^{AD,q}} \right),		\\				
  && P_i^{+,q}  = \max \left( {0, -  g_{i+1/2}^{AD,q}} \right) + \max \left( {0, g_{i-1/2}^{AD,q}} \right) .
\end{eqnarray*}

Finally, we rewrite \eqref{eq:3.16} in the form
\begin{equation}
\label{eq:4.7}
 {A_i} \le {\alpha _{i+1/2}}{d_{ii+1}} + {\alpha _{i-1/2}}{d_{ii - 1}} ,			
\end{equation} 
where
\begin{eqnarray*}
 && {A_i} = \frac{\Delta x}{\Delta t}\left( { U(\boldsymbol {\hat v}_i) - U(\boldsymbol {v_i}) - \left\langle U'(\boldsymbol v_i), \left( {\boldsymbol {\hat v}_i - \boldsymbol  v_i} \right) \right\rangle } \right) + G_{i+1/2}^{L} - G_{i-1/2}^L  \\
 && \qquad \quad  - \left\langle U'(\boldsymbol {v_i}), ( \boldsymbol g_{i+1/2}^L - \boldsymbol g_{i-1/2}^L - \boldsymbol s^-_{i+1/2} +\boldsymbol  s^+_{i-1/2} ) \right\rangle ,\\
 && d_{ik} = \left( \langle U'(\boldsymbol v_i),\boldsymbol g_{(i+k)/2}^{AD}  \rangle - G_{(i+k)/2}^{AD} \right) {\rm sgn} (k-i) .
\end{eqnarray*} 

By reasoning similar to the above, we obtain from \eqref{eq:4.7}  that the upper bound of ${\alpha _{i+1/2}}$ is equal to 
\begin{equation}
\label{eq:4.8} 
\begin{split} 
  {\bar \alpha} _{i+1/2}^U =  \min \left\{ 1, \frac{-A_i}{B_i}\min \left( {0,{\mathop{\rm sgn}} \, d_{ii+1}} \right) +\max \left( {0, {\rm sgn} \, d_{ii+1}} \right), \right. \\
 \left. \frac{-A_{i+1}}{ B_{i+1}} \min \left( {0,{\mathop{\rm sgn}} \, d_{i+1i}} \right) +\max \left( {0, {\rm sgn} \, d_{i+1i}} \right) 
  \right\}	, 
\end{split}   				
\end{equation} 
where
${B_i} = \min (0,d_{ii+1}) + \min (0,d_{ii-1})$ 

Thus, a nontrivial feasible solution to the linear programming problem \eqref{eq:3.13}-\eqref{eq:3.16} on ${U_{ad}}$ is equal to
\begin{equation}
\label{eq:4.9}
{\alpha _{i+1/2}} = \min ({\bar \alpha_{i+1/2}^y}, \;{\bar \alpha} _{i+1/2}^q, \;{\bar \alpha} _{i+1/2}^U) .
\end{equation}

\begin{theorem} 
\label{th:4_1}
Let $ U(\boldsymbol {\hat v}): R^{N} \times R^N \rightarrow R^{N} \times R^N$ be a Lipschitz-continuous function on a closed ball $O_0 = \bar O(\boldsymbol {\hat v}^L,\delta)$, where $\boldsymbol {\hat v}^L$ is a solution of the system of equations \eqref{eq:1.1}. Then the flux limiters $\boldsymbol \alpha$, defined by \eqref{eq:4.9}, are Lipschitz-continuous on $ O_0$.
\end{theorem} 
\begin{proof} It is clear that ${\bar \alpha}_{i+1/2}^y, \bar \alpha_{i+1/2}^q, B_i, d_{ik}$ are constants, and $A_i$ are Lipschitz-continuous functions on $O_0$. Thus, $\alpha_{i+1/2}, {{\bar \alpha}}_{i+1/2}^U$ are  Lipschitz-continuous on $O_0$, since the minimum of Lipschitz-continuous functions is again a Lipschitz-continuous function.
\end{proof}

\begin{remark} 
The hypotheses of Theorems \ref{th:3.1} and \ref{th:4_1} are satisfied if $U(\boldsymbol v)$ is a strictly convex function,  and $\alpha_{i+1/2}$ are calculated using \eqref{eq:4.1}-\eqref{eq:4.9}. In this case, the system of equations \eqref{eq:3.1} has a unique solution.  
\end{remark}

\begin{remark}
The approach presented in this paper can be extended to multidimensional and implicit HR schemes. For details we refer the reader to~\cite{kivva2021}.
\end{remark}

\section{Numerical Examples}
\label{sec:5}

In this section, we demonstrate the benefits of the proposed approach and compare numerical results with analytical and previous numerical studies. 
We also compare numerical results obtained with  flux limiters, which are  approximate and exact solutions to the corresponding optimization problems.

Applying the centered space flux as a high-order flux, we use the following hybrid HR scheme in our calculations
\begin{equation}
\label{eq:5.1}
\begin{split}
  & \boldsymbol {\hat v_i} - \boldsymbol {v_i} + \frac{\Delta t} {{\Delta x}}\left[ \boldsymbol{g_{i+1/2}^L(\boldsymbol{v^-_{i+1/2}},\boldsymbol{v^+_{i+1/2}}) - \boldsymbol g_{i-1/2}^L(\boldsymbol{v^-_{i-1/2}},\boldsymbol{v^+_{i+1/2}})} \right] \\
 + \frac{1}{2} \frac{\Delta t} {{\Delta x}} & \left[ \alpha_{i+1/2} \; c_{i+1/2} \; (\boldsymbol v^+_{i+1/2} -\boldsymbol v^-_{i+1/2}) - \alpha_{i-1/2} \; c_{i-1/2} \; (\boldsymbol v^+_{i-1/2} -\boldsymbol v^-_{i-1/2}) \right] = \Delta t \; \boldsymbol s_i,	
\end{split}
\end{equation}
where $\boldsymbol g_{i+1/2}^L$ is the Rusanov numerical flux \eqref{eq:1.2}.

Then the discrete cell entropy inequality \eqref{eq:10} can be written in the form
\begin{equation} 
\label{eq:5.2} 
\begin{split}
& U(\boldsymbol {\hat v_i}) - U({\boldsymbol v_i}) + \frac{{\Delta t}}{{\Delta x}}\left[ {G^L_{i+1/2}(\boldsymbol v^-_{i+1/2},\boldsymbol v^+_{i+1/2})} - {G^L_{i-1/2}(\boldsymbol v^-_{i-1/2},\boldsymbol v^+_{i-1/2})} \right. \\
+  \alpha_{i+1/2} \; & \frac{c_{i+1/2}}{2} \left.  ( U(\boldsymbol v^+_{i+1/2}) -U (\boldsymbol v^-_{i+1/2})) -  \alpha_{i-1/2} \; \frac{c_{i-1/2}}{2} \; ( U(\boldsymbol v^+_{i-1/2}) - U(\boldsymbol v^-_{i-1/2})) \right] \le 0 
\end{split}
\end{equation}
with the proper numerical entropy flux \eqref{eq:2.9}.

Below we will mark the numerical solutions of scheme \eqref{eq:5.1}-\eqref{eq:5.2} with a label indicating how the flux limiters are calculated. The letters $L$ and $A$ denote the applying linear programming or approximate solution to the optimization problem, respectively. The letters $H$, $Q$ and $E$ mean that the flux limiters were calculated using inequalities \eqref{eq:3.14}, \eqref{eq:3.15} and \eqref{eq:5.2}, respectively. Numerical solutions with flux limiters satisfying inequalities \eqref{eq:1.13} are denoted as $PP$. In the latter case, flux limiters are defined as follows
\begin{equation}
\label{eq:5.3}
  \alpha_{i+1/2} = \max \left( 0, \min \left( 1, 1-\frac{u_{i+1}}{c_{i+1/2}}, 1+\frac{u_{i}}{c_{i+1/2}} \right)  \right).  
\end{equation}

In addition, we use the following labels:
\begin{description}
\item{HR1} is a first-order hydrostatic reconstruction scheme with HLL numerical flux given in \cite{Chen2017};
\item{HR2} is a hydrostatic reconstruction scheme of second-order   spatial accuracy with explicit Euler time integration proposed in \cite{Buttinger2019};
\item{ZL} is a characteristic variable implementation of the Boris-Book flux limiter described in \cite{b3} and applied to the HR scheme \eqref{eq:1.1}-\eqref{eq:1.2}.
\end{description}
 
To solve linear programming problems we apply GLPK package v.4.65 ($\href{url}{https://www.gnu.}$ $\href{url}{org/software/glpk/}$).

\subsection{One-Dimensional Dam Break Over a Wet Flat Bed}
\label{sec:51}

In this section, we consider a dam break on a wet flat bed in a frictionless, horizontal, rectangular channel. The channel is 1000 m long. The dam is located in the middle of the channel. The water depth at the left and right hand sides of the dam is $100$ m and $1$ m, respectively. The dam instantly collapses across its entire width and the resulting flow consists 
of a shock wave traveling downstream and a rarefaction wave traveling upstream. In this problem there is a transition from subcritical upstream to supercritical downstream flows. The simulation is performed up to time t=10 s.

The analytical solution of this problem was given by Stoker (1957) \cite{stoker}.
The 1D dam-break on a wet flat bed is a classical test to verify the shock-capturing ability of numerical schemes.

\begin{figure*}[!t]
  \centering 
  \includegraphics{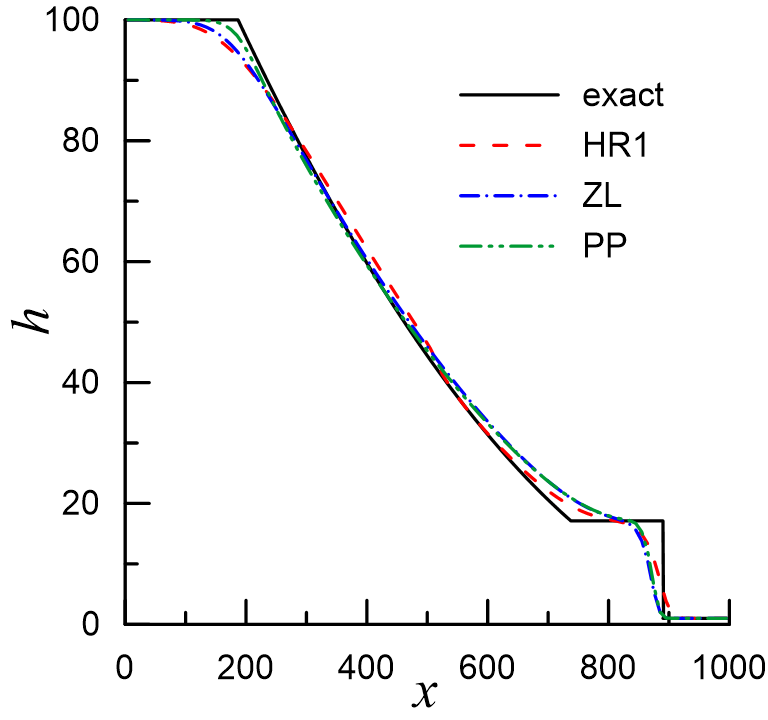}
  \includegraphics{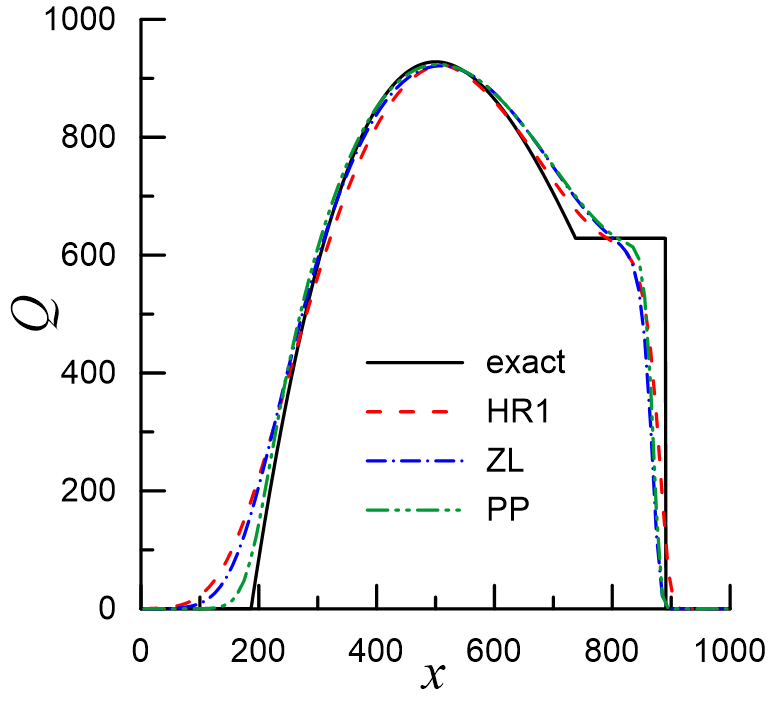}
\caption{One-dimensional dam break over a wet flat bed. Comparisons of exact solutions with simulated water depths (left) and discharges (right) using HR1, ZL, and PP at t=10 s. The number of cells is N=100.}
\label{fig:1}       
\end{figure*}

\begin{figure*}[!t]
  \centering 
  \includegraphics{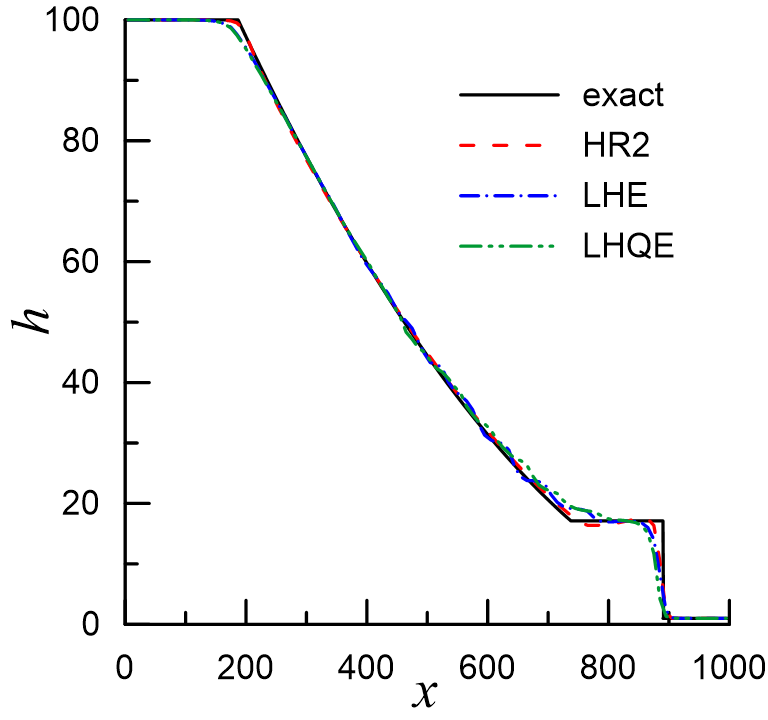}
  \includegraphics{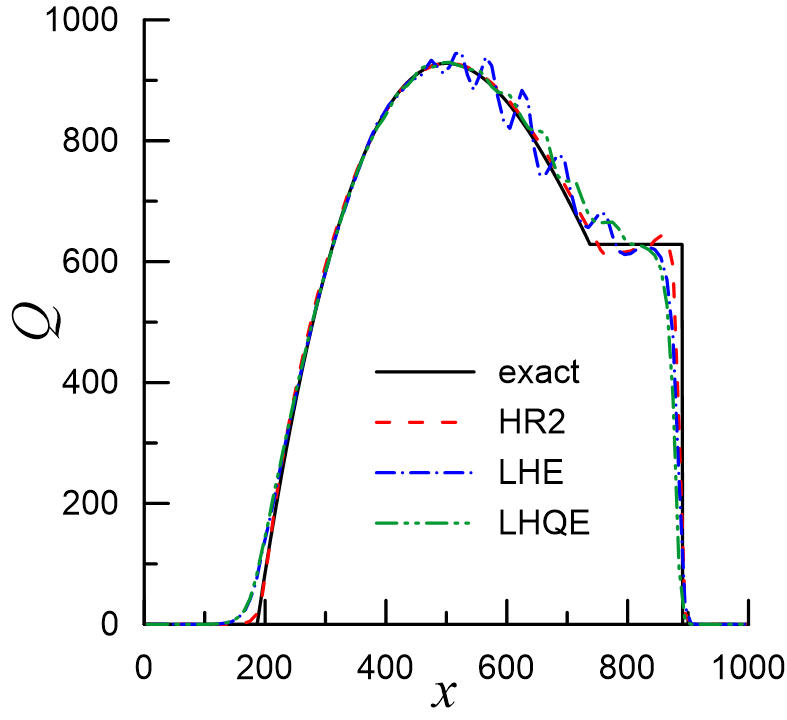}
  \includegraphics{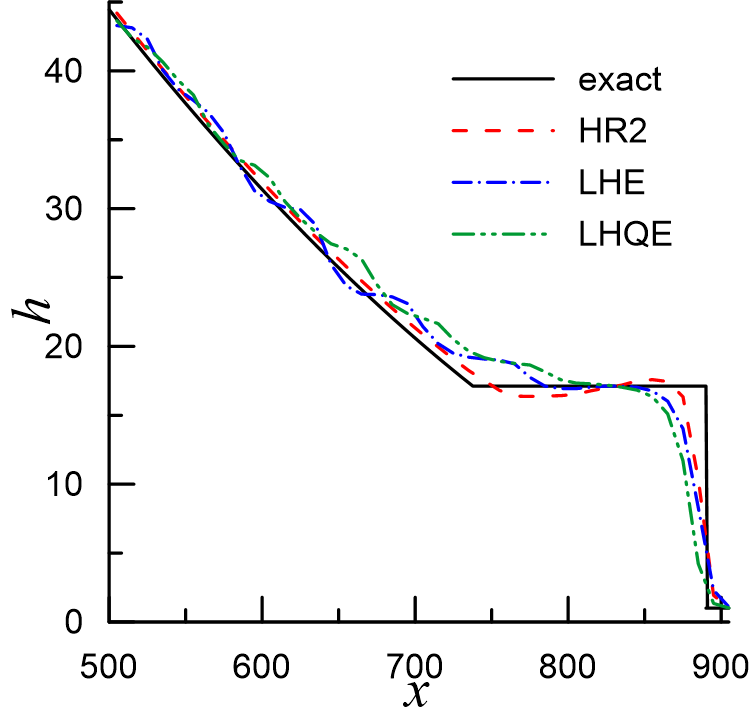}
  \includegraphics{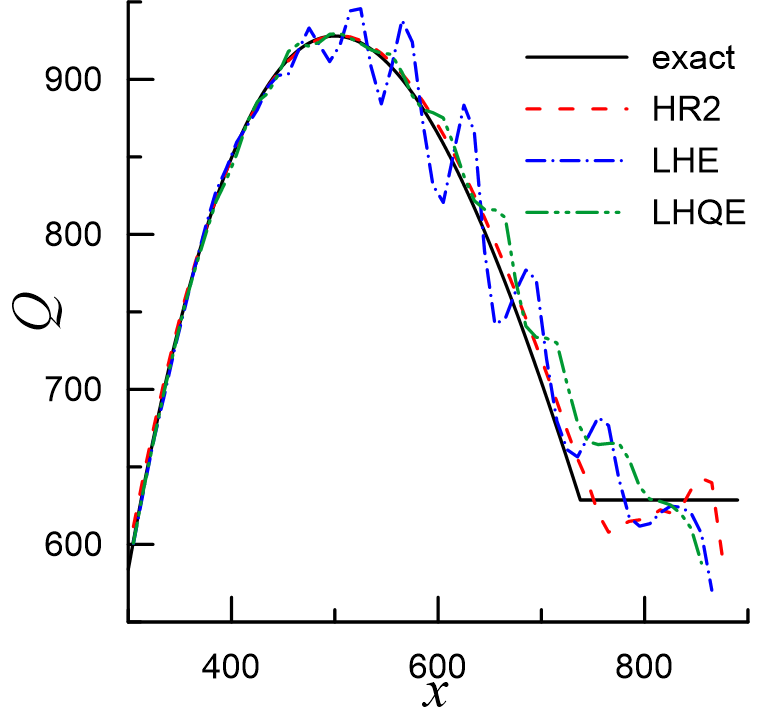}
\caption{One-dimensional dam break over a wet flat bed. Comparisons of exact solutions with simulated water depths (left) and discharges (right) using HR2, LHE, and LHQE at t=10 s. The second row is a zoom in the area behind the shock. The number of cells is N=100.}
\label{fig:2}       
\end{figure*}

\begin{figure*}[!t]
  \includegraphics[width=0.32\textwidth]{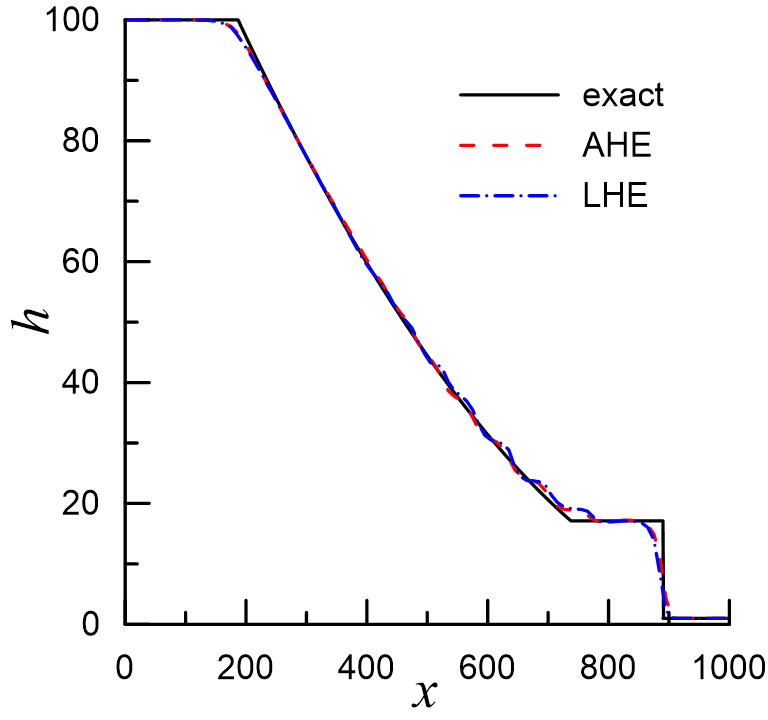}
  \includegraphics[width=0.32\textwidth]{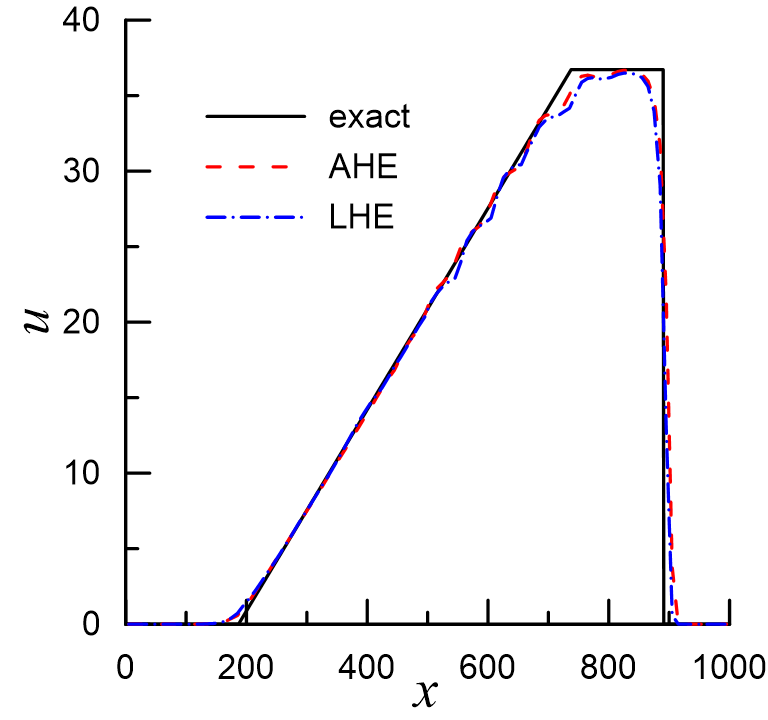}
  \includegraphics[width=0.32\textwidth]{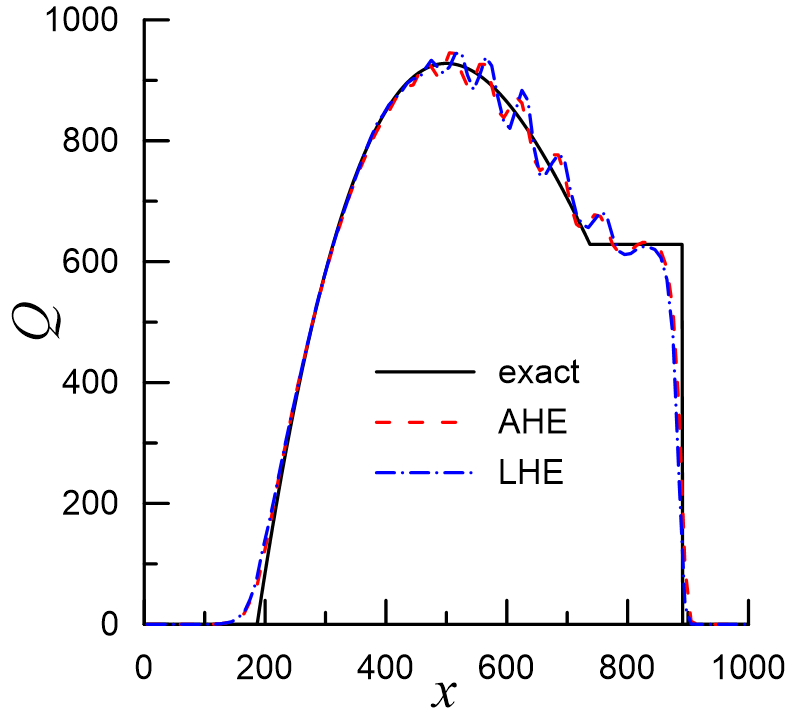}
  \includegraphics[width=0.32\textwidth]{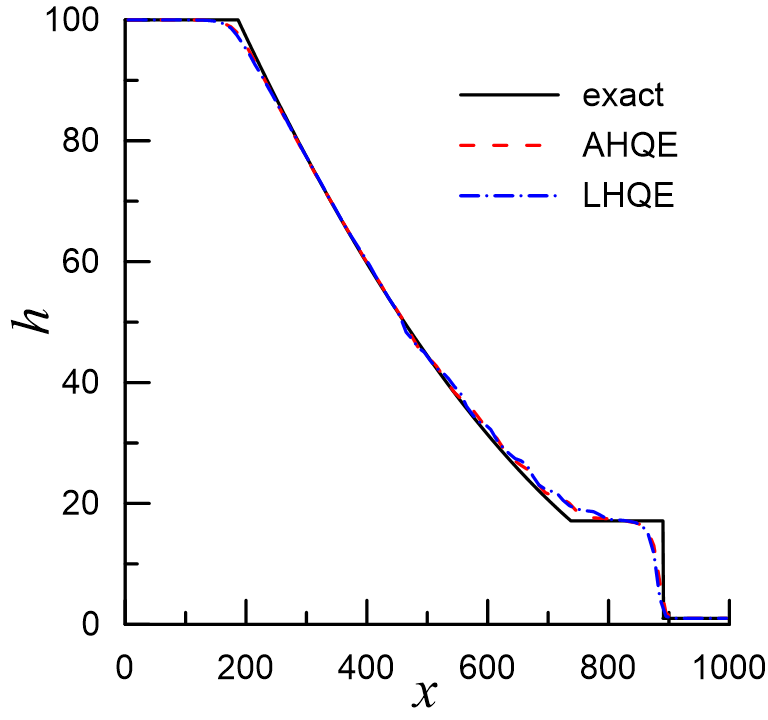}
  \includegraphics[width=0.32\textwidth]{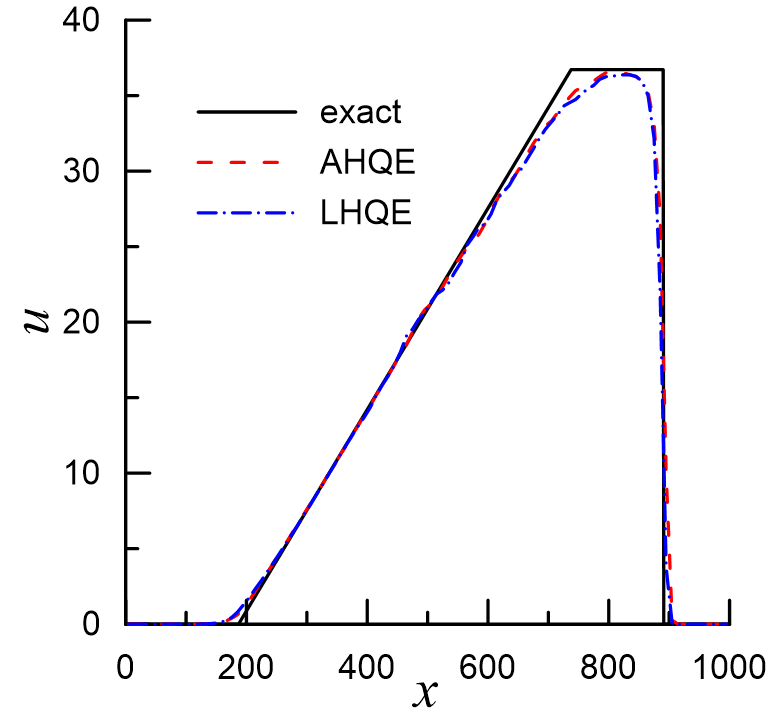}
  \includegraphics[width=0.32\textwidth]{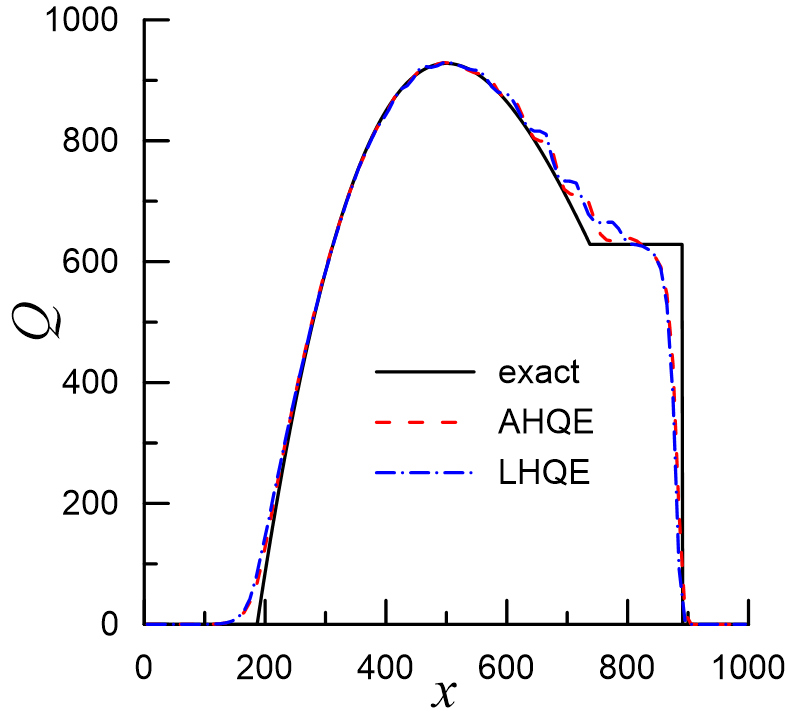}
\caption{One-dimensional dam break over a wet flat bed. Comparisons of numerical results obtained with FCT schemes  whose flux limiters are computed using exact and approximate solutions to a linear programming problem with discrete entropy inequality  and different constraints. The number of cells is N=100.}
\label{fig:3}       
\end{figure*}

Numerical results obtained with different schemes at time t=10 s on a uniform grid of N=100 cells are shown in Fig.~\ref{fig:1}-Fig.~\ref{fig:3}. As shown in Fig.~\ref{fig:2}, the shock wave resolutions using HR2, LHE, and LHQE 
are less dissipative (sharper) and better than with the other schemes shown in Fig.~\ref{fig:1}. 

The simulated results with PP are close to those obtained with ZL but require much less calculations. In the numerical results with LHE, AHE, LHQE, and AHQE, we observe the so-called "terracing" phenomenon characteristic of FCT methods. Numerical results for water discharge obtained using the LHE and AHE schemes have oscillations that are absent in the velocities (Fig~\ref{fig:3}).

The analytical and numerical solutions were compared quantitatively by the L1 error. The error is defined as

\begin{equation}
\label{eq:6.4}
L^1= \frac{1}{N} \sum\limits_{i=1}^N \vert y_i -y^a(x_i) \vert
\end{equation}
where $y_i$ is the numerical and $y^a$ is the analytical solution at point $x_i$, $N$ means the number of these points.
Table~\ref{tab1} shows the $L^1$-norms of errors of the numerical solutions obtained with different schemes.

%
  
%
  
%

\begin{table*}[!h]
\caption{\label{tab1}
$L^1$-norms of errors for the numerical solutions of the 1D dam break over a wet flat bed at t=10 s with N=100.
}
\centering
\begin{tabular}{lcccc@{} p{0.5cm} *{4}{p{4cm} @{}}}
\hline 
  & $HR1$  &  $ZL$  & $PP$ & $HR2$  &  \\[2pt] \hline
    H &	1.468$\times 10^{0}$ & 1.679$\times 10^{0}$ & 1.365$\times 10^{0}$ & 4.052$\times 10^{-1}$ &  \\	
	Q &	3.596$\times 10^{1}$ & 3.882$\times 10^{1}$ & 3.060$\times 10^{1}$ & 9.180$\times 10^{0}$ &  \\
\hline
 & $LHE$ & $AHE$ & $LHQE$ & $AHQE$ & \\[2pt] \hline 	
H & 6.153$\times 10^{-1}$ & 5.114$\times 10^{-1}$ & 7.898$\times 10^{-1}$ & 6.306$\times 10^{-1}$ & \\
Q &	1.912$\times 10^{1}$ & 1.619$\times 10^{1}$ & 2.003$\times 10^{1}$ & 1.647$\times 10^{1}$ & \\
\hline	    
\end{tabular}
\end{table*}

A comparison of analytical solutions with computed depths, as well as velocities and discharges at t=10 s using LHE(LHQE) and AHE(AHQE) are given in Fig.~\ref{fig:3}. The flux limiters for LHE(LHQE) and AHE(AHQE) are calculated using exact and approximate solutions to linear programming problems. We note good agreement between these numerical solutions, and the addition of  constraints on water discharges to calculate flux limiters leads to suppression of oscillations in the numerical solutions.

\subsection{One-Dimensional Dam Break Over a Dry Bed} 
\label{sec:52}

The dry bed dam-break test is usually applied to verify the ability of a difference scheme to propagate a wet/dry front at the correct speed and to keep water depth positive. The analytical solution of this problem was given by Stoker (1957) \cite{stoker}.

We consider a rectangular channel with $1000$ m length and a flat bed. The dam is located in the middle of the channel. The water depth at the left and right hand sides of the dam is 100 m and 0 m, respectively. 
The dam break is instantaneous and there is no friction. The solution consists of a single rarefaction wave with a wet/dry front at its lower end. 

The flow domain is discretized into 100 uniform cells. The simulation time is t=7 s.

  
\begin{figure*}[!t]
  \centering 
  \includegraphics{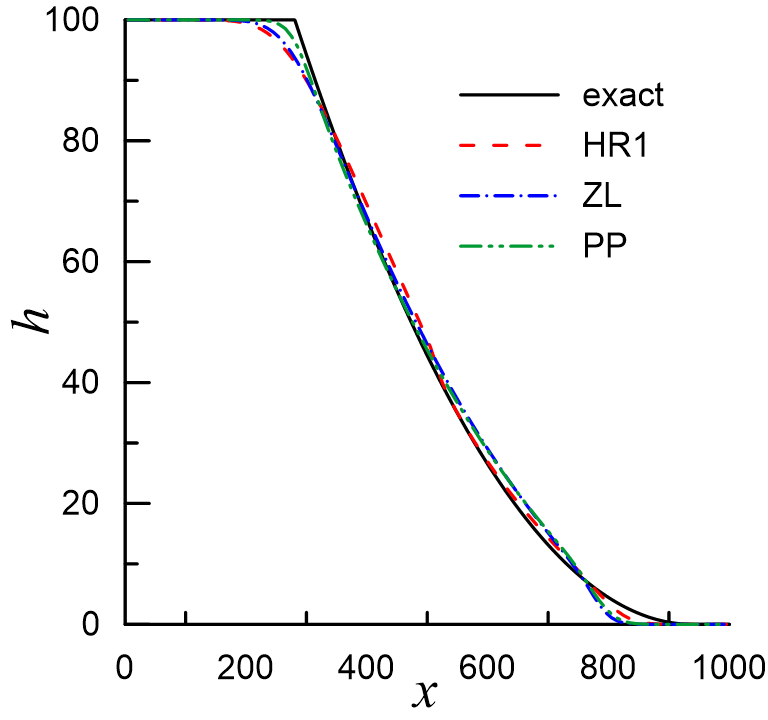}
  \includegraphics{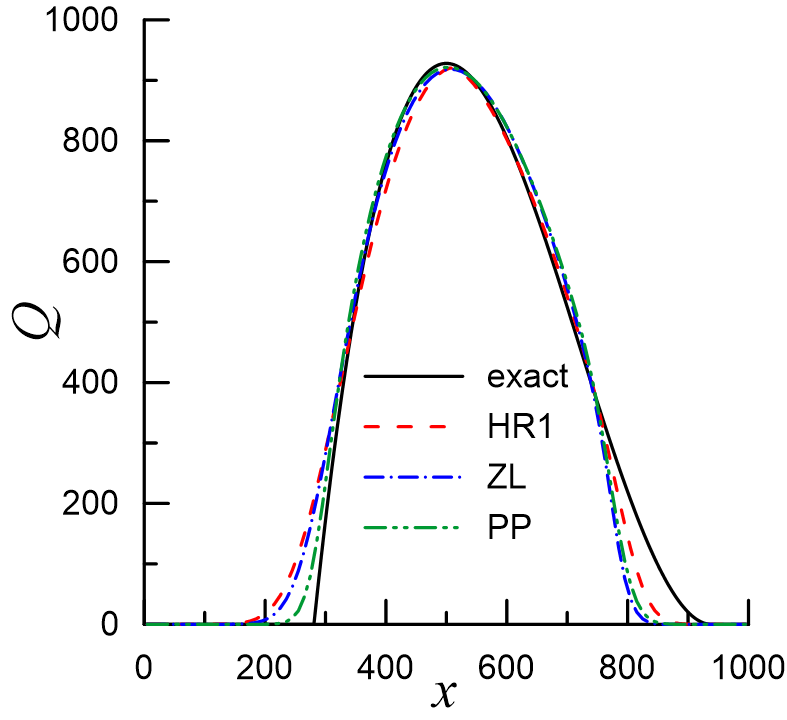}
\caption{One-dimensional dam break over a dry bed. Comparisons of exact solutions with simulated water depths and discharges using HR1, ZL, and PP at time t=7 s. The number of cells is N=100.}
\label{fig:4}       
\end{figure*}
%
  
\begin{figure*}[!t]
  \centering 
  \includegraphics{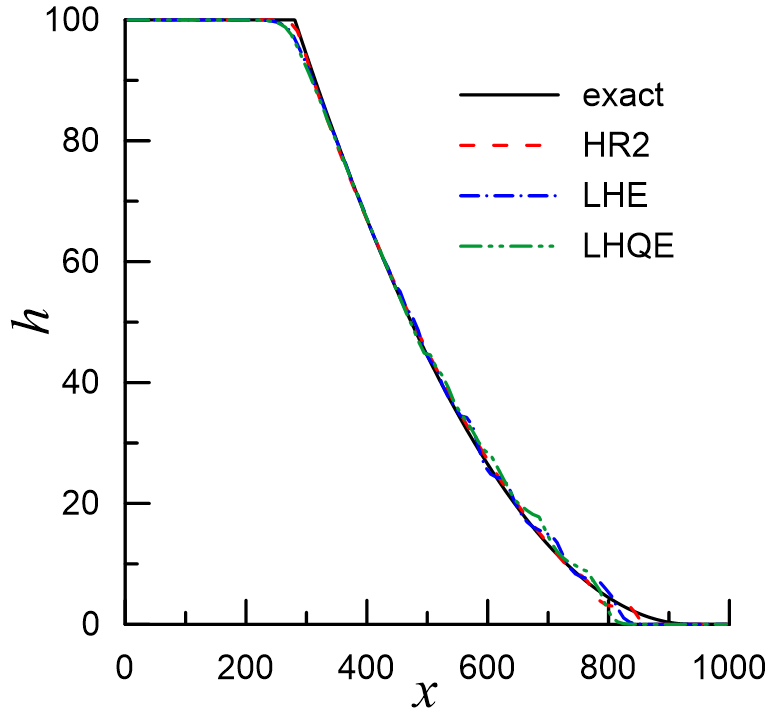}
  \includegraphics{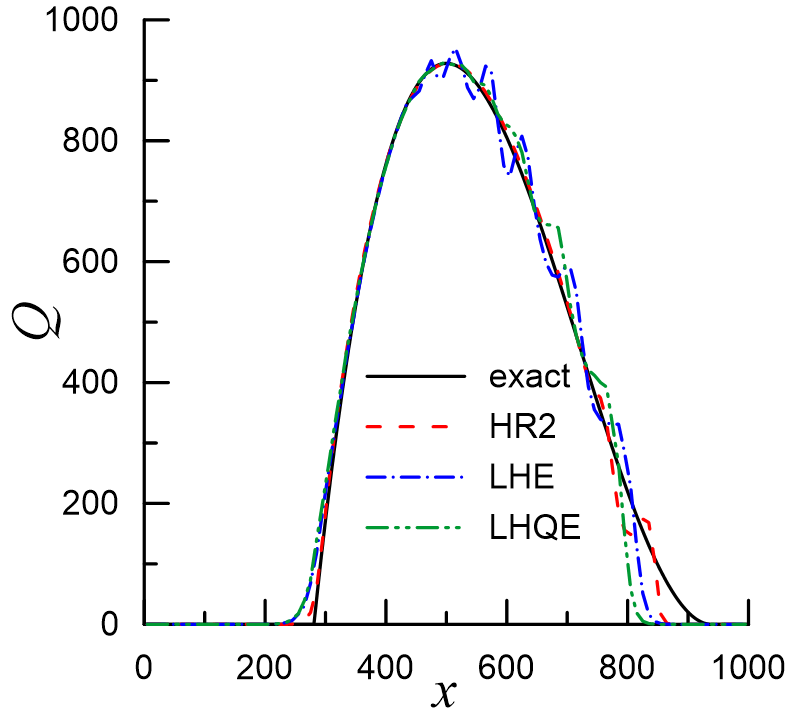}
  \includegraphics{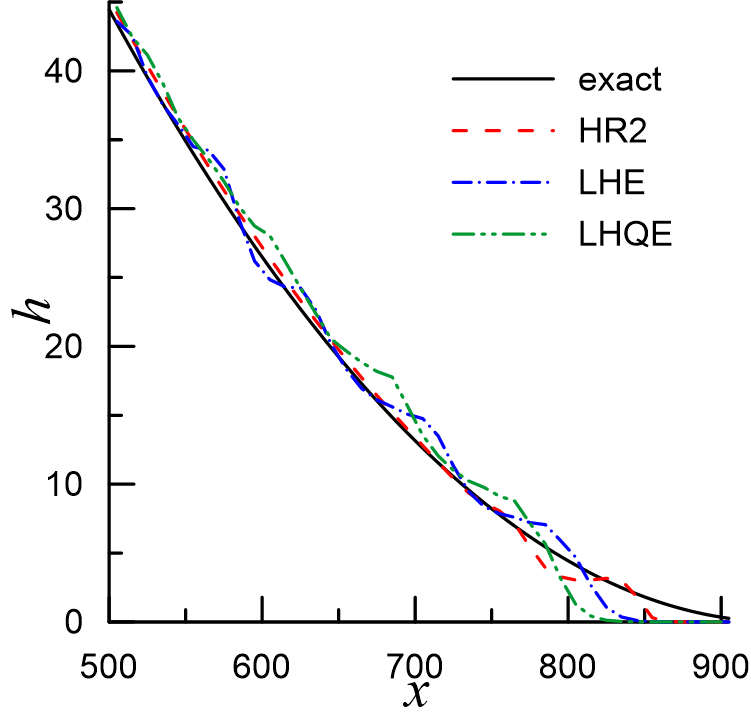}
  \includegraphics{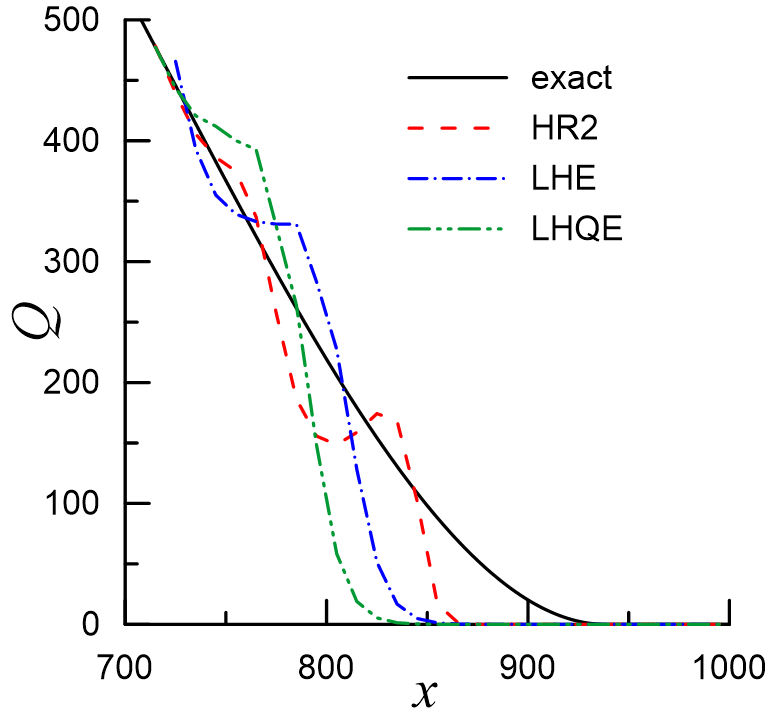}
\caption{One-dimensional dam break over a dry bed. Comparisons of exact solutions with simulated water depths and discharges using HR2, LHE, and LHQE at time t=7 s. The second row is a zoom in the area of the front of the moving water. The number of cells is N=100.}
\label{fig:5}       
\end{figure*}
%
  
\begin{figure*}[!t]
  \includegraphics[width=0.32\textwidth]{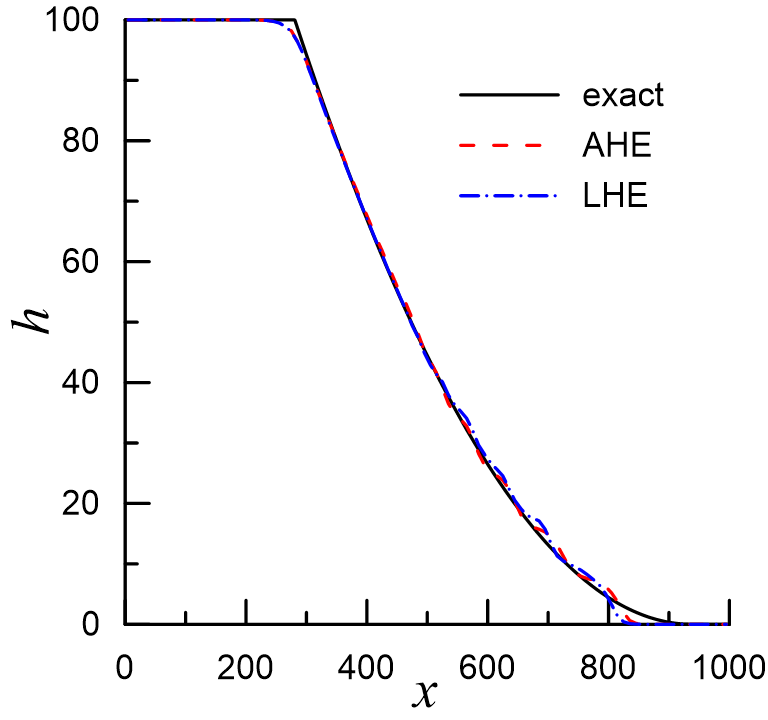}
  \includegraphics[width=0.32\textwidth]{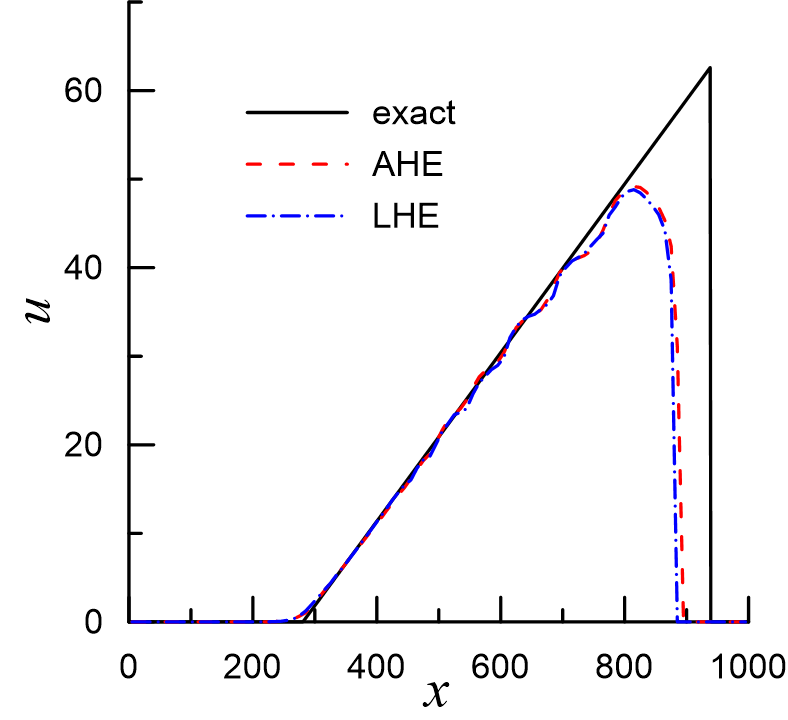}
  \includegraphics[width=0.32\textwidth]{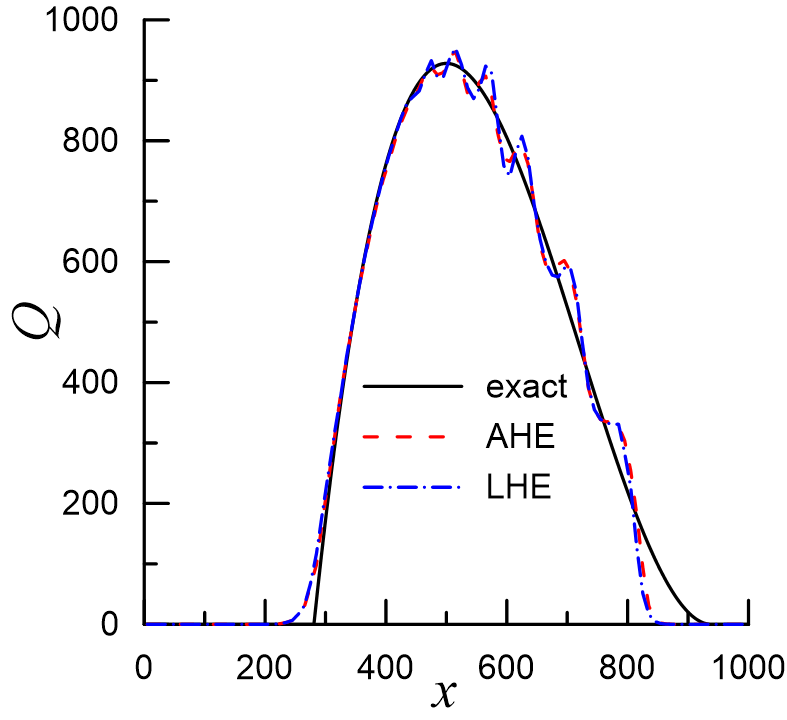}
  \includegraphics[width=0.32\textwidth]{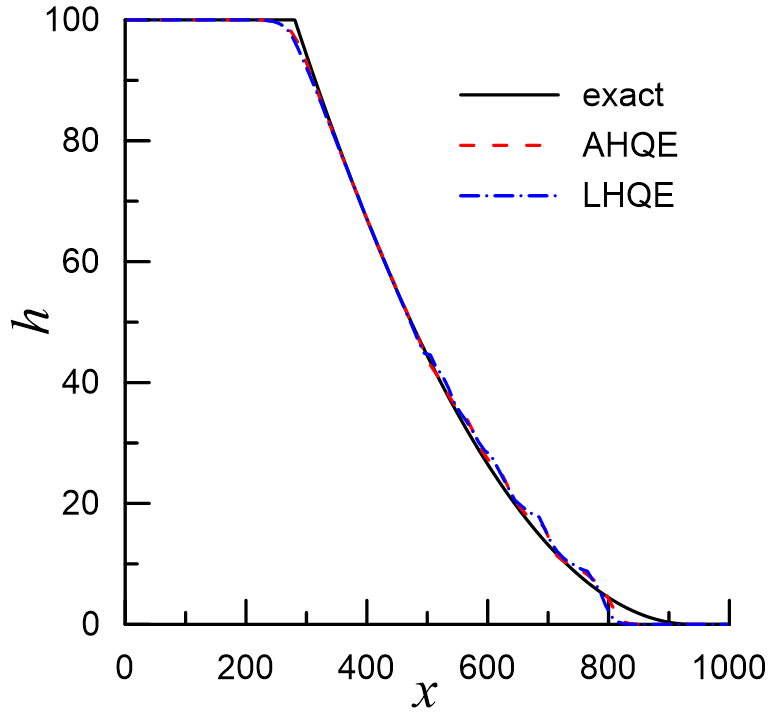}
  \includegraphics[width=0.32\textwidth]{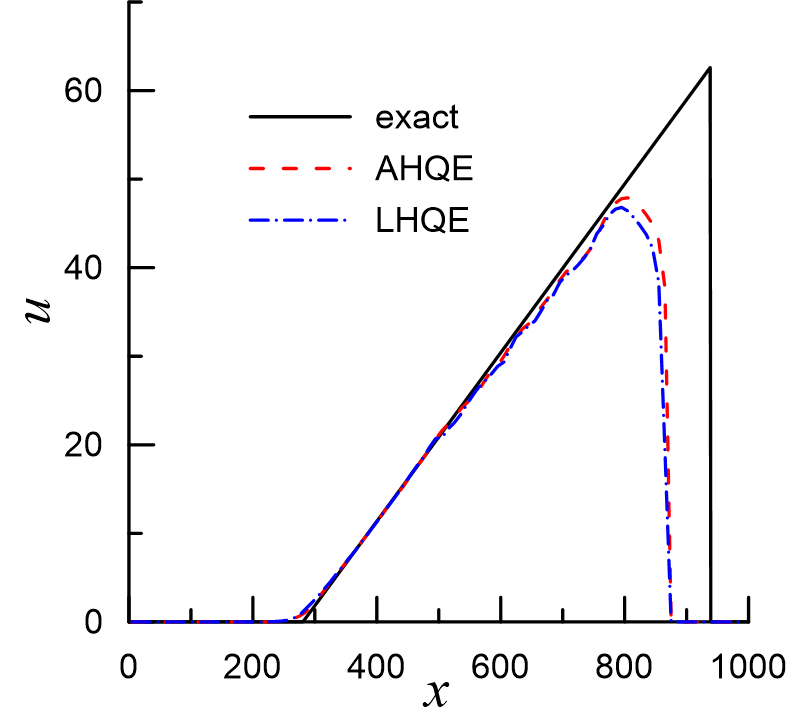}
  \includegraphics[width=0.32\textwidth]{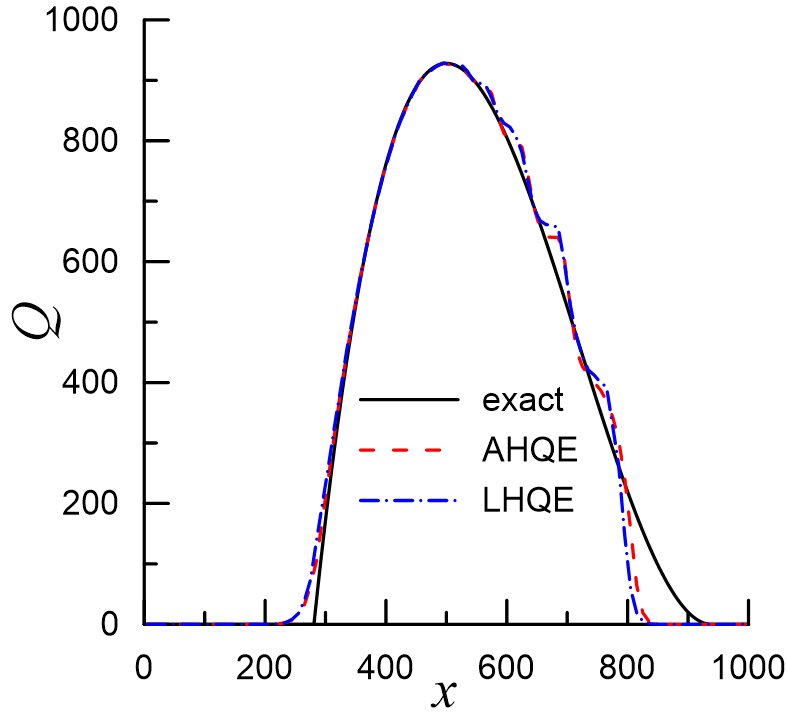}
\caption{One-dimensional dam break over a dry bed. Comparisons of numerical results obtained with FCT schemes  whose flux limiters are computed using exact and approximate solutions to a linear programming problem with discrete entropy inequality  and different constraints at time t=7 s.}
\label{fig:6}       
\end{figure*}
%

\begin{table*}[!h]
\caption{\label{tab2}
$L^1$-norms of errors for the numerical solutions of the 1D dam break over a dry bed at t=7 s with N=100.
}
\centering
\begin{tabular}{lcccc@{} p{0.5cm} *{4}{p{4cm} @{}}}
\hline 
  & $HR1$  &  $ZL$  & $PP$ & $HR2$  &  \\[2pt] \hline
    H &	1.145$\times 10^{0}$ & 1.320$\times 10^{0}$ & 1.050$\times 10^{0}$ & 3.684$\times 10^{-1}$ &  \\	
	Q &	2.838$\times 10^{1}$ & 3.218$\times 10^{1}$ & 2.590$\times 10^{1}$ & 1.038$\times 10^{1}$ &  \\
\hline
 & $LHE$ & $AHE$ & $LHQE$ & $AHQE$ & \\[2pt] \hline 	
H & 5.463$\times 10^{-1}$ & 5.216$\times 10^{-1}$ & 7.690$\times 10^{-1}$ & 5.786$\times 10^{-1}$ & \\
Q &	2.129$\times 10^{1}$ & 1.978$\times 10^{1}$ & 2.207$\times 10^{1}$ & 1.820$\times 10^{1}$ & \\
\hline	    
\end{tabular}
\end{table*}  
  
Comparisons of exact solutions with simulated depths as well as discharges at t=7 s using the six schemes are presented in Fig.~\ref{fig:4}-\ref{fig:5}. Among the proposed schemes, the HR1, ZL, and PP schemes present more dissipative results than the HR2, LHE, and LHQE schemes. The simulated results with PP are close to those obtained with ZL but require much less calculations. In the numerical results obtained with LHE, AHE, LHQE, and AHQE, we observe the so-called "terracing" phenomenon, which is characteristic of FCT methods. The LHE and AHE schemes produce oscillations in the water discharges that are absent in the velocities (Fig~\ref{fig:6}). For all the considered schemes, the 
largest error is observed at the front of the moving water.

Adding constraints on water discharges to the LHQE scheme to calculate flux limiters eliminates oscillations in numerical solutions. The numerical results obtained with LHE(LHQE) and AHE(AHQE) are in a good agreement (Fig.~\ref{fig:6}). The flux limiters for LHE(LHQE) and AHE(AHQE) are calculated using exact and approximate solutions to linear programming problems.

\subsection{Dam Break Over a Step.}
\label{sec:53}

In this test \cite{Buttinger2019}, a dam break over a downward bottom step is considered. The bottom topography and the initial data are given as follows
\begin{eqnarray}
\label{eq:6.5}
  z(x) = 
  \begin{cases} 
   1 &    \text{if  } x \leq 0, \\
   0 &   \text{otherwise  } ,
  \end{cases} \quad
 h(x,0) =
  \begin{cases} 
   0.75 &    \text{if  } x \leq 0, \\
   1.0 &   \text{otherwise  } ,
  \end{cases} \quad
  Q(x,0) = 0,
\end{eqnarray}

After a dam break, the solution consists of a left rarefaction wave, a stationary shock wave at an intermediate height of the bottom step between two stationary contact waves located at the bottom discontinuity, and a right shock wave~\cite{han2014,minatti2023}.

\begin{figure*}[!t]
  \centering 
  \includegraphics{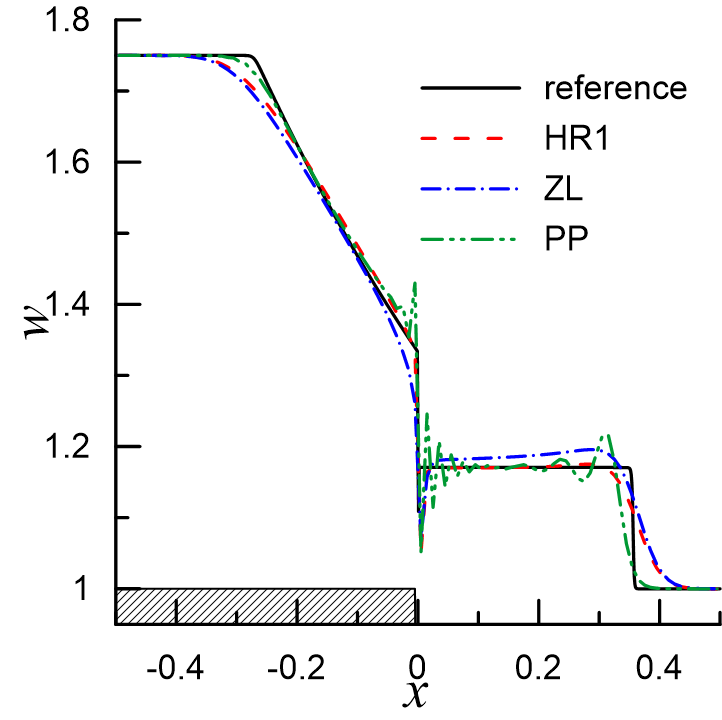}
  \includegraphics{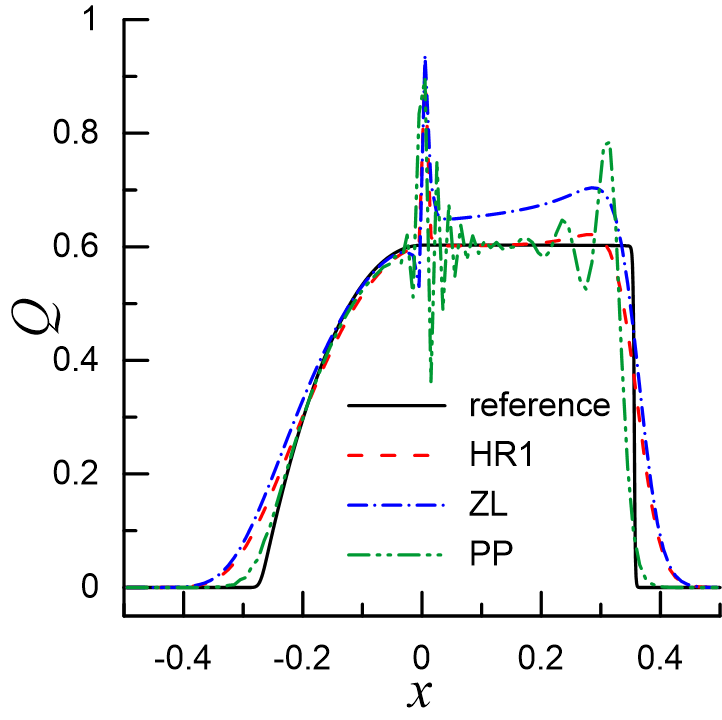}
\caption{Dam break over a step. Comparisons of reference solutions with simulated water depths and discharges using the HR1, ZL, and PP schemes at time t=0.1 s with N=200 cells.}
\label{fig:7}       
\end{figure*}
%
  
\begin{figure*}[!t]
  \centering 
  \includegraphics{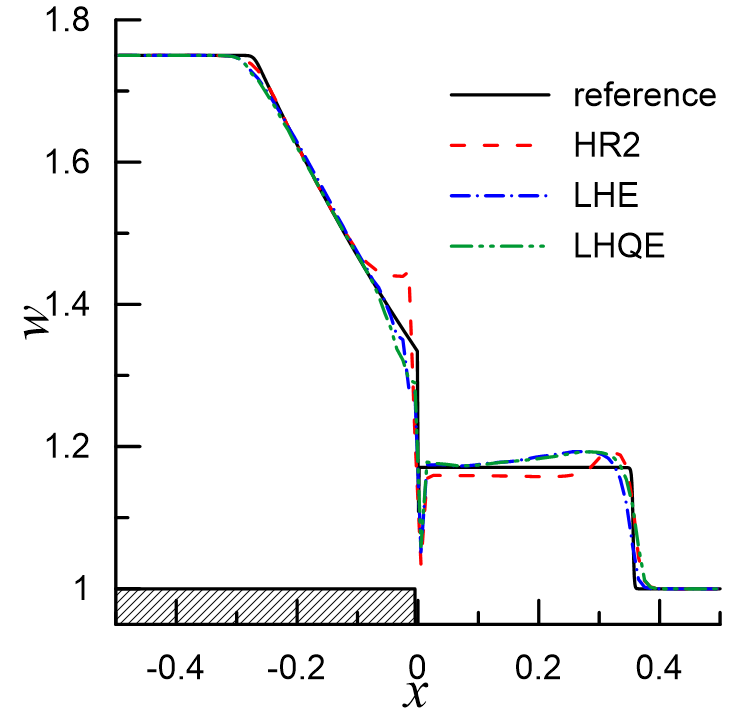}
  \includegraphics{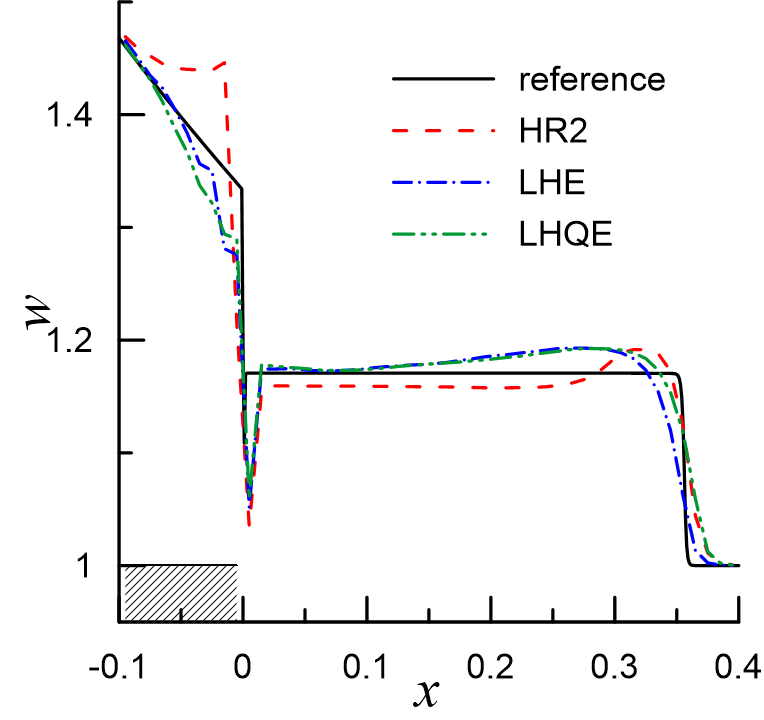}
  \includegraphics{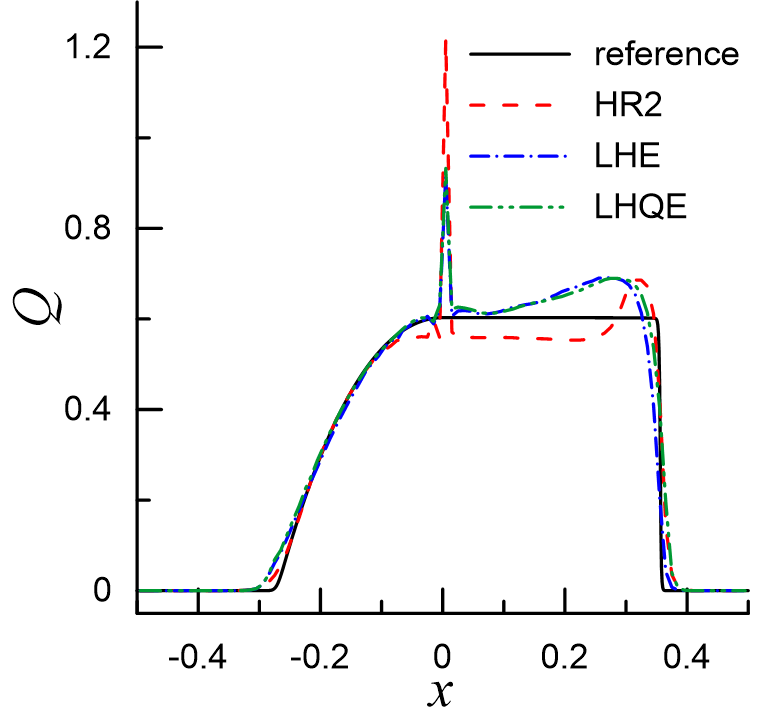}
  \includegraphics{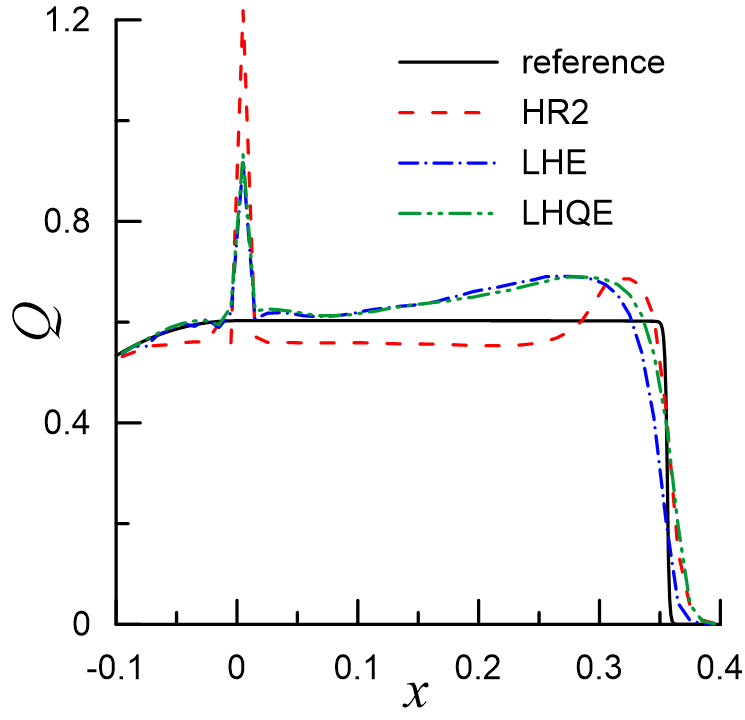}
\caption{Dam break over a step. Comparisons of reference solutions with simulated water depths and discharges using HR2, LHE, and LHQE at time t=0.1 s with N=200 cells. On the right is a zoom of the area of the bottom discontinuity and the right shock wave.}
\label{fig:8}       
\end{figure*}
%
  
\begin{figure*}[!t]
  \includegraphics[width=0.32\textwidth]{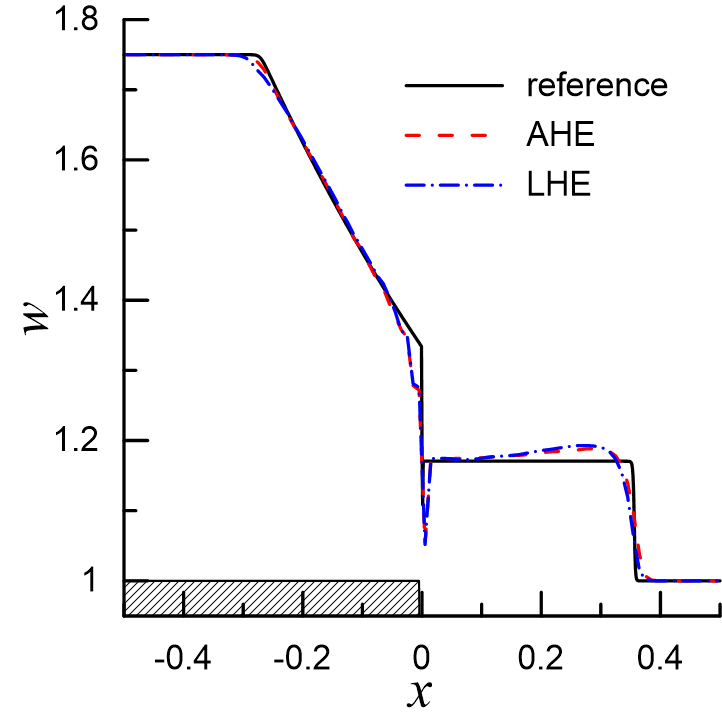}
  \includegraphics[width=0.32\textwidth]{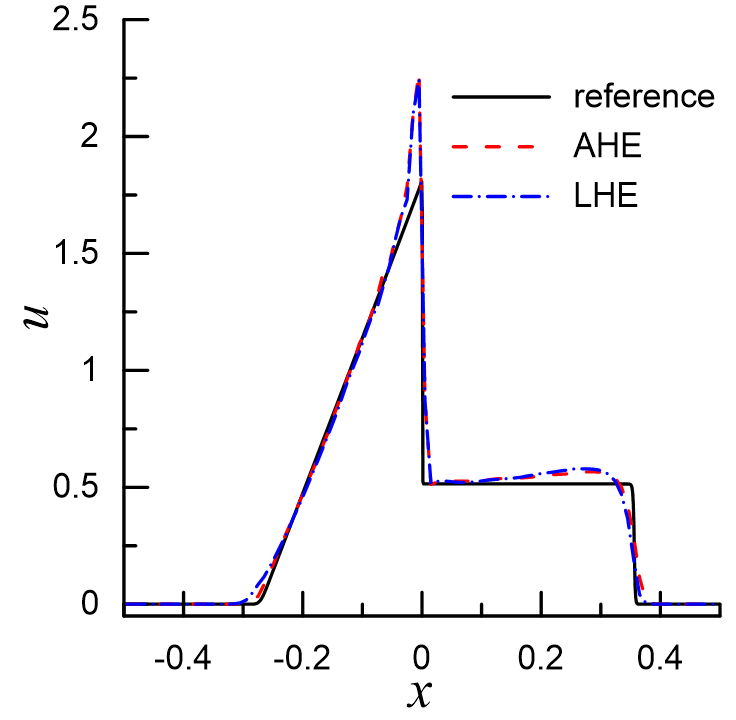}
  \includegraphics[width=0.32\textwidth]{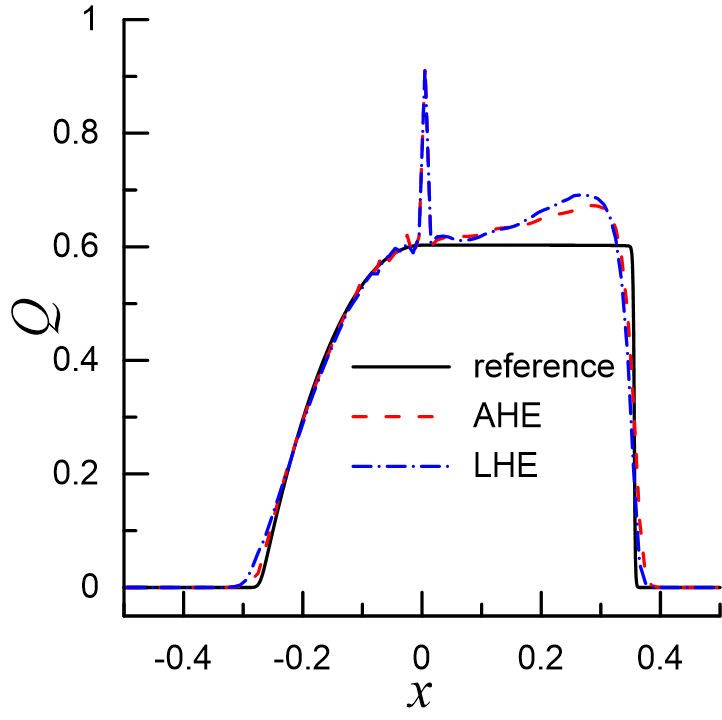}
  \includegraphics[width=0.32\textwidth]{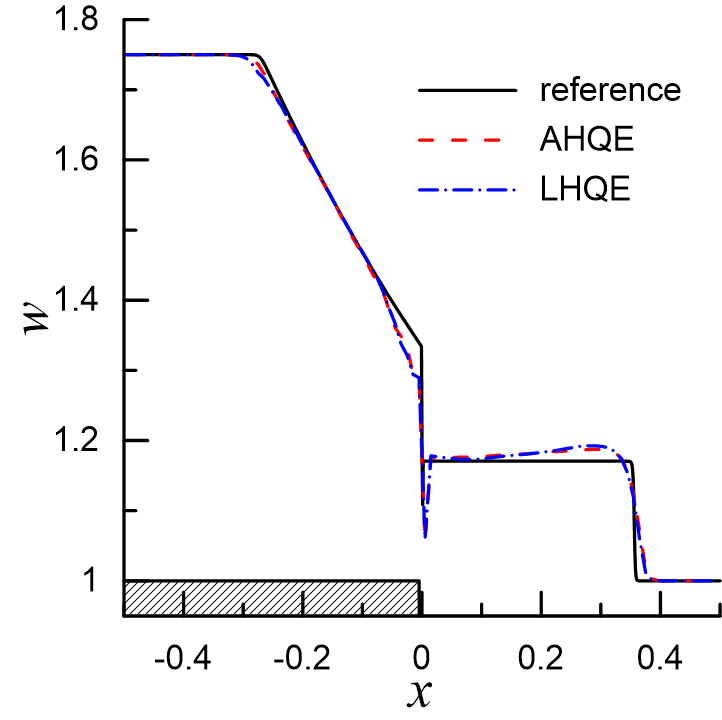}
  \includegraphics[width=0.32\textwidth]{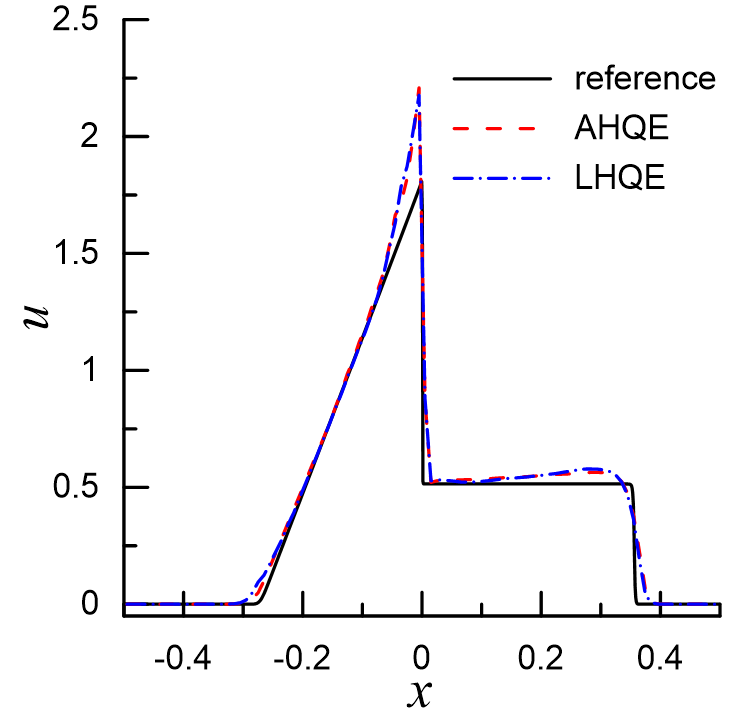}
  \includegraphics[width=0.32\textwidth]{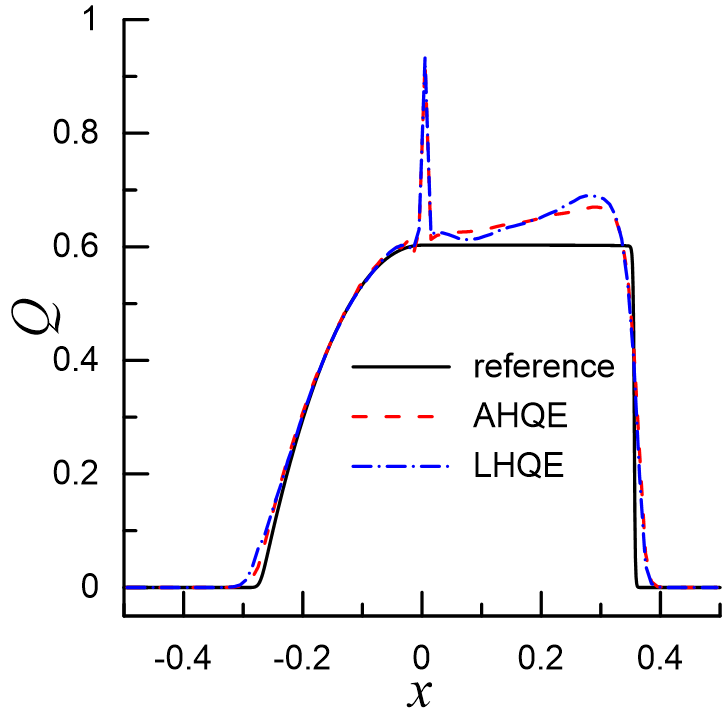}
\caption{Dam break over a step. Comparisons of numerical results obtained with FCT schemes whose flux limiters are computed using exact and approximate solutions to a linear programming problem with discrete entropy inequality  and different constraints. (t=0.1 s, N=200).}
\label{fig:9}       
\end{figure*}

Comparisons of the numerical results obtained on a uniform grid of 200 cells with a reference solution at $t$=0.1 after the dam  break are shown in Fig.~\ref{fig:7}-\ref{fig:9}. 
The reference solution was calculated using a central-upwind scheme of second-order spatial accuracy \cite{kivva2020} on a uniform grid with 2000 cells. In Fig.~\ref{fig:7}, the PP scheme generates oscillations in the numerical results in the area of the bottom discontinuity.
In the numerical results obtained with the ZL scheme, we see an overshoot of the water depth and discharge for the right shock wave. 
The second-order HR2 scheme does not reproduce the
left rarefaction wave in its whole entirety, as well as the shock wave (Fig.~\ref{fig:8}). 
In Fig.~\ref{fig:8}-\ref{fig:9}, the right side of the shock wave for the LHE, LHQE, AHE, and AHQE schemes shows an overshoot of the simulated water depth and discharge.
We note that none of the considered schemes reproduces the exact solution, especially in the bottom discontinuity.  
Table~\ref{tab3} shows the L1-norm error between the reference solution and the numerical solutions at time t = 0.1 for different difference schemes.

\begin{table*}[!h]
\caption{\label{tab3}
$L^1$-norms of errors for the numerical solutions of the 1D dam break over a step at $t=0.1 s$ with $N$=200.
}
\centering
\begin{tabular}{lcccc@{} p{0.5cm} *{4}{p{4cm} @{}}}
\hline 
  & $HR1$  &  $ZL$  & $PP$ & $HR2$  &  \\[2pt] \hline
    H &	5.219$\times 10^{-3}$ & 8.647$\times 10^{-3}$ & 6.039$\times 10^{-3}$ & 6.840$\times 10^{-3}$ &  \\	
	Q &	1.390$\times 10^{-2}$ & 2.489$\times 10^{-2}$ & 1.806$\times 10^{-2}$ & 1.703$\times 10^{-2}$ &  \\
\hline
 & $LHE$ & $AHE$ & $LHQE$ & $AHQE$ & \\[2pt] \hline 	
H & 4.830$\times 10^{-3}$ & 4.100$\times 10^{-3}$ & 4.695$\times 10^{-3}$ & 4.656$\times 10^{-3}$ & \\
Q &	1.326$\times 10^{-2}$ & 1.096$\times 10^{-2}$ & 1.243$\times 10^{-2}$ & 1.209$\times 10^{-2}$ & \\
\hline	    
\end{tabular}
\end{table*}  
  
We also note that the numerical results obtained with LHE(LHQE) and AHE(AHQE) agree well (Fig.~\ref{fig:9}). The flux limiters for LHE(LHQE) and AHE(AHQE) are calculated using exact and approximate solutions to linear programming problems.

\subsection{Steady Transcritical Flow With a Shock Over a Bump.}
\label{sec:54}

We consider a test taken from~\cite{delestre2012-1} consisting of a transcritical flow with a shock over a bump. The bed topography of a rectangular channel 25 m long is given as follows
\begin{eqnarray}
\label{eq:6.6}
  z(x) = 
  \begin{cases} 
   0.2-0.05(x-10)^2 &    \text{if  } 8 < x < 12, \\
   0 &   \text{otherwise  } .
  \end{cases} \quad
\end{eqnarray}

Initial conditions satisfy the hydrostatic equilibrium
\begin{equation}
\label{eq:6.7}
h+z=0.33 \quad \text{and} \quad Q=0.
\end{equation}

Discharge $Q$ = 0.18 $\text{m}^2$/s and water level $h+z$ = 0.33 m were set as upstream and downstream boundary conditions.
In the steady-state solution, the flow to the left of the bump is subcritical, then closer to the end of the bump it becomes supercritical, and after a hydraulic jump it is subcritical again.

\begin{figure*}[!h]
  \centering 
  \includegraphics[scale=0.92]{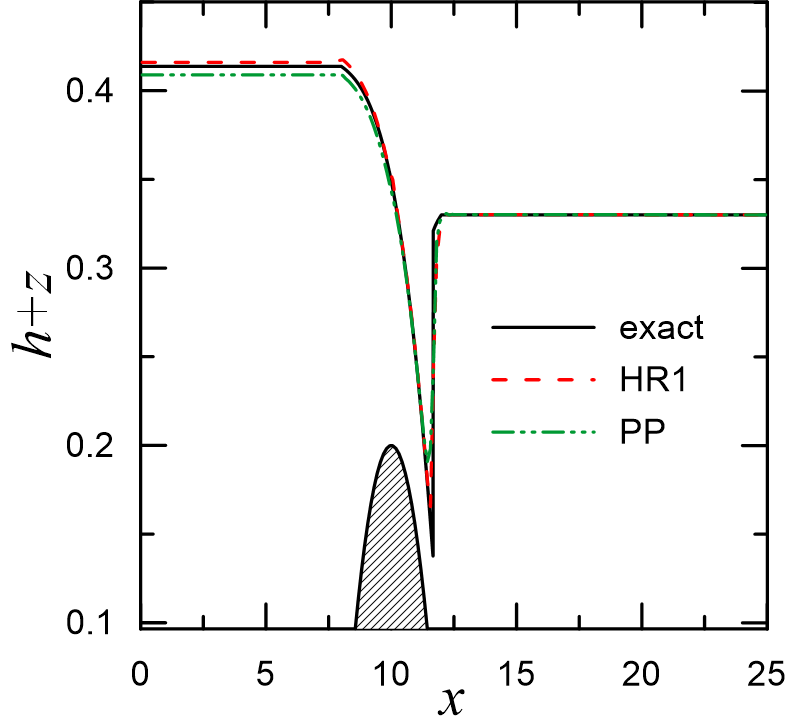}
  \includegraphics[scale=0.92]{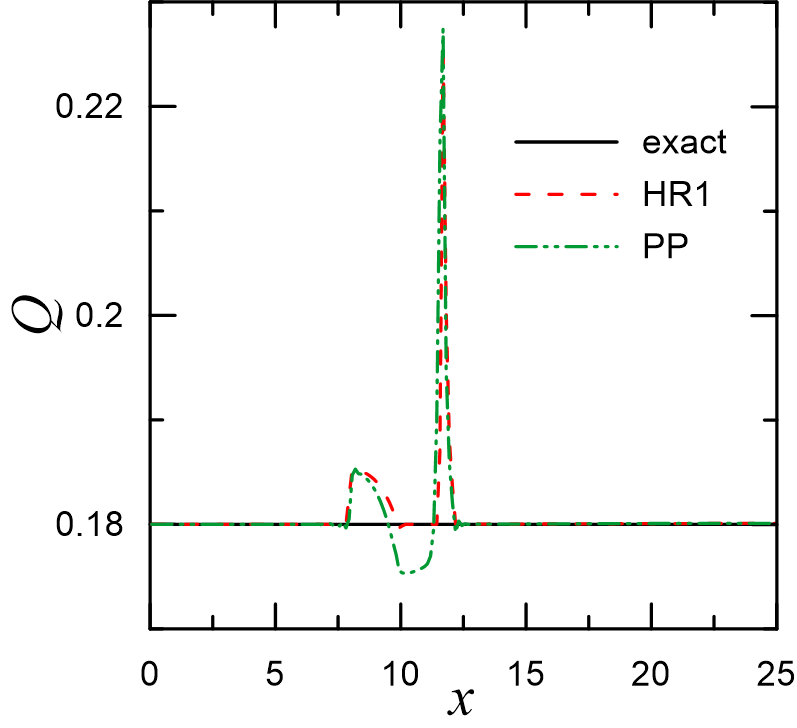}
  \includegraphics[scale=0.92]{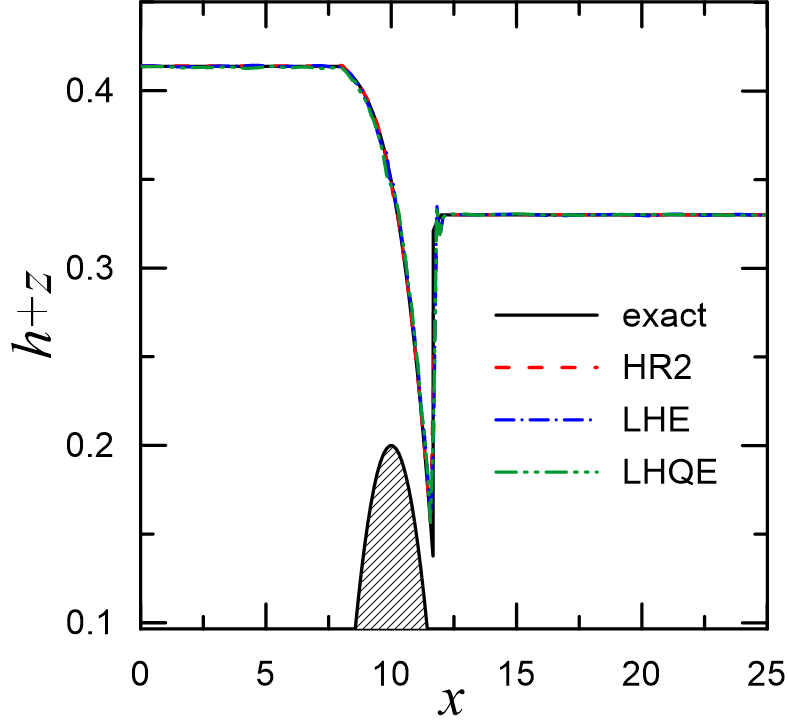}
  \includegraphics[scale=0.92]{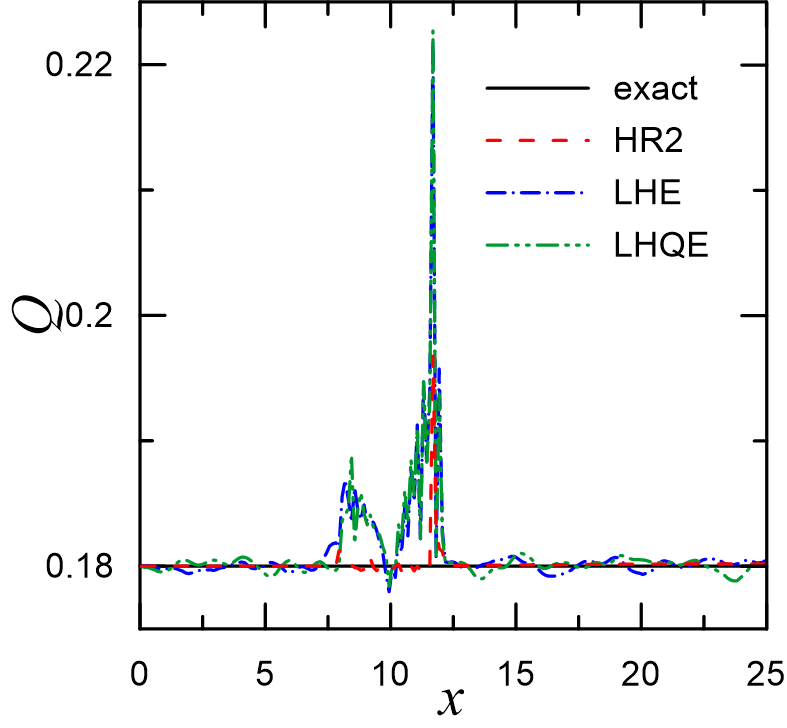}
\caption{Steady transcritical flow with a shock over a bump. Comparison of exact solutions with computed water depths and discharges obtained by HR1, PP, HR2, LHE, and LHQE with N=100.}
\label{fig:10}       
\end{figure*}
%
  
\begin{figure*}[!h]
  \includegraphics[width=0.32\textwidth]{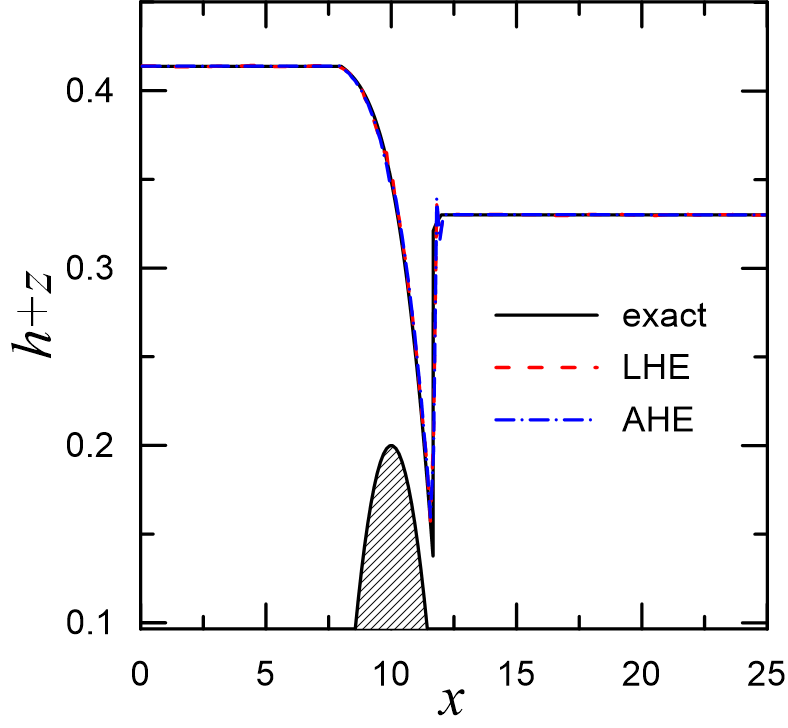}
  \includegraphics[width=0.32\textwidth]{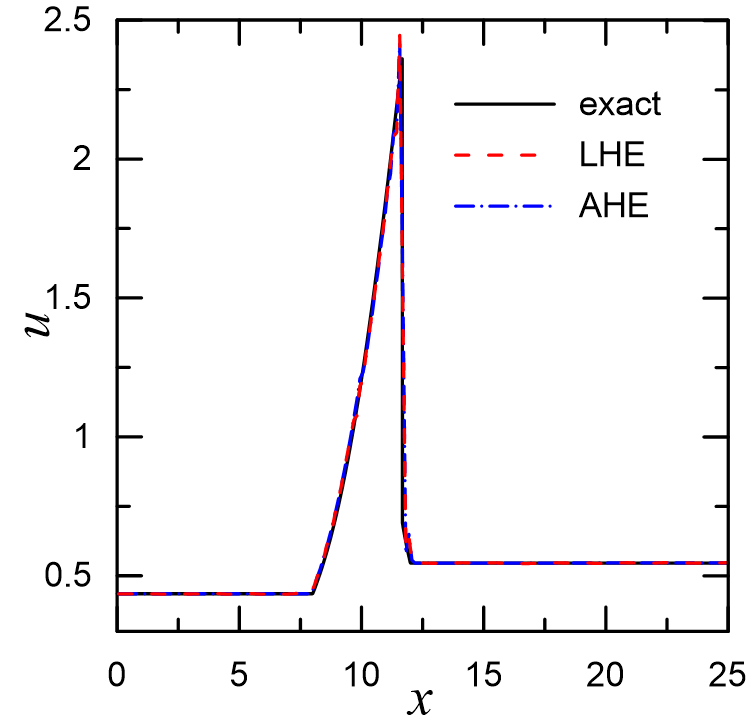}
  \includegraphics[width=0.32\textwidth]{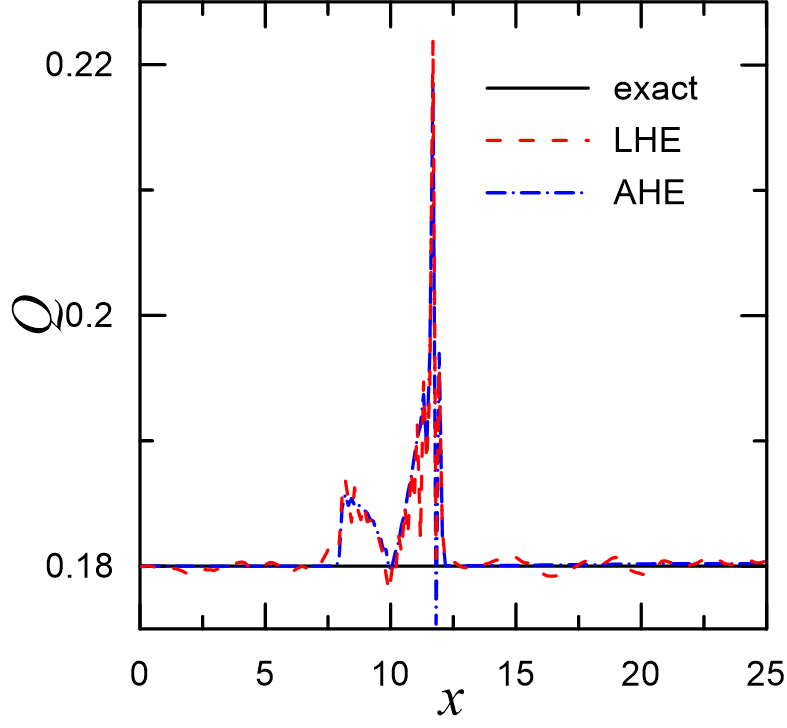}
  \includegraphics[width=0.32\textwidth]{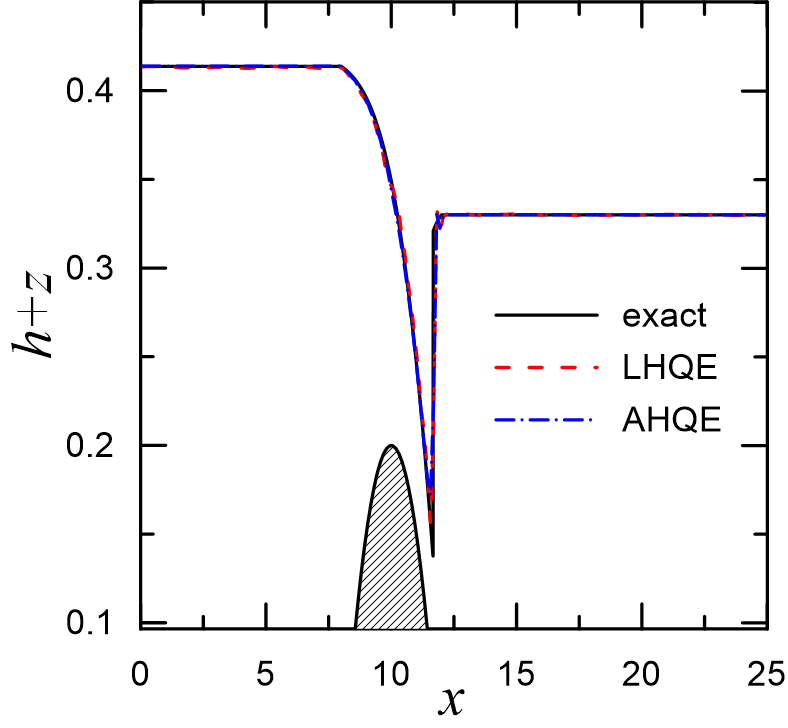}
  \includegraphics[width=0.32\textwidth]{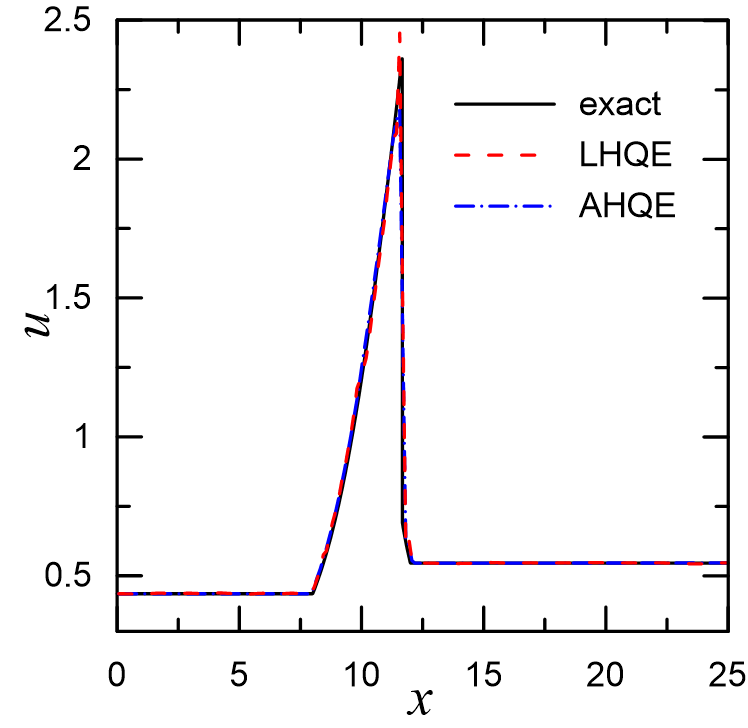}
  \includegraphics[width=0.32\textwidth]{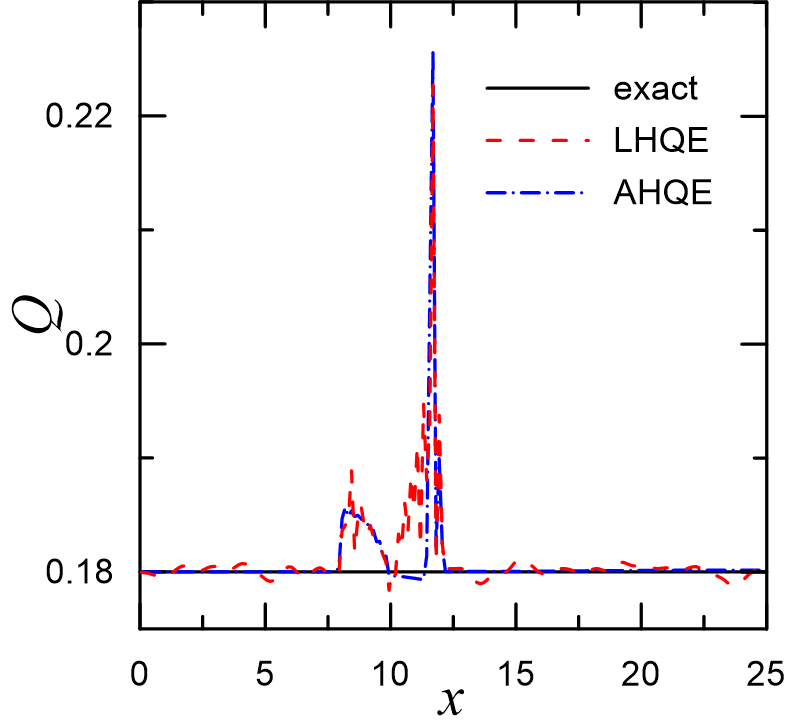}
\caption{Steady transcritical flow with a shock over a bump. Comparisons of numerical results obtained with FCT schemes  whose flux limiters are computed using exact and approximate solutions to a linear programming problem with  discrete entropy inequality  and different constraints.}
\label{fig:11}       
\end{figure*}
%

\begin{table*}[!h]
\caption{\label{tab4}
$L^1$-norms of errors of the transcritical
steady state flow with a shock over a bump.
}
\centering
\begin{tabular}{lcccc@{} p{0.5cm} *{4}{p{4cm} @{}}}
\hline 
  & $HR1$  & $PP$ & $HR2$  &  \\[2pt] \hline
    H &	1.633$\times 10^{-3}$ & 2.858$\times 10^{-3}$ & 6.258$\times 10^{-4}$ &  \\	
	Q &	7.534$\times 10^{-4}$ & 1.106$\times 10^{-3}$ & 2.201$\times 10^{-4}$ &  \\
\hline
 & $LHE$ & $AHE$ & $LHQE$ & $AHQE$ & \\[2pt] \hline 	
H & 1.298$\times 10^{-3}$ & 1.298$\times 10^{-3}$ & 1.501$\times 10^{-3}$ & 9.556$\times 10^{-4}$ & \\
Q &	1.270$\times 10^{-3}$ & 1.087$\times 10^{-3}$ & 1.282$\times 10^{-3}$ & 8.034$\times 10^{-4}$ & \\
\hline	    
\end{tabular}
\end{table*}  
  
Numerical results for the steady state, obtained on a uniform grid of 100 cells, are shown in Fig.~\ref{fig:10}-\ref{fig:11}.
In the numerical results obtained with the HR1 scheme, we observe an overshoot of the free surface before the bump and an undershoot of the free surface for the PP scheme.
The free water surfaces calculated with HR2, LHE, and LHQE agree fairly well with the analytical solution, with slight deviations around the hydraulic jump.
Numerical oscillations for water discharges near the hydraulic jump are present for all compared schemes. 
Small oscillations are also present in the calculated discharges with the LHE and LHQE schemes in the whole modeling area. The $L^1$-norm error between the exact and numerical solutions are shown in Table~\ref{tab4}.

Note that none of the considered schemes is well-balanced for moving water steady states with non-zero discharges. 

We also note that the numerical results obtained with LHE(LHQE) and AHE(AHQE) agree well (Fig.~\ref{fig:11}). The flux limiters for LHE(LHQE) and AHE(AHQE) are calculated using exact and approximate solutions to linear programming problems.

\begin{figure*}[!htbp]
  \centering 
  \includegraphics[scale=0.92]{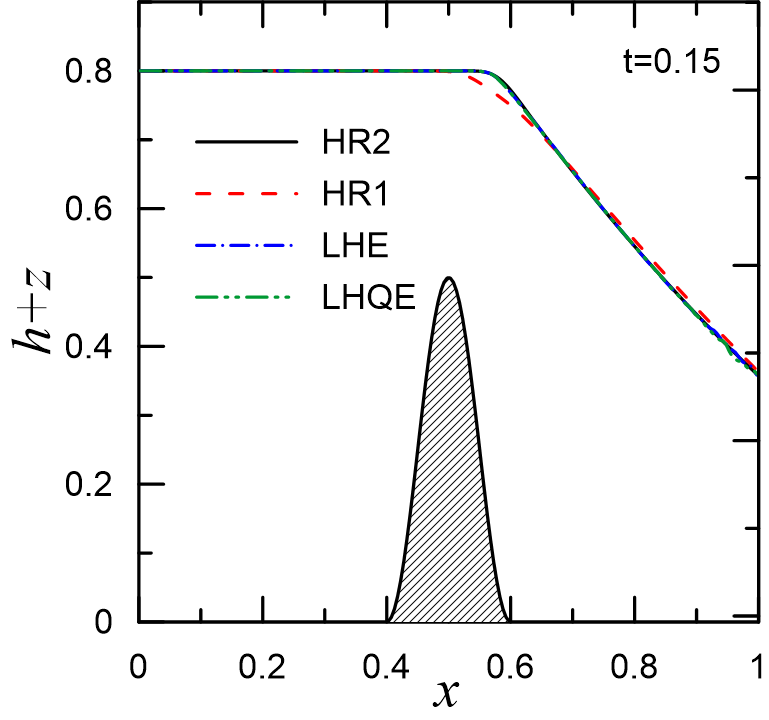}
  \includegraphics[scale=0.92]{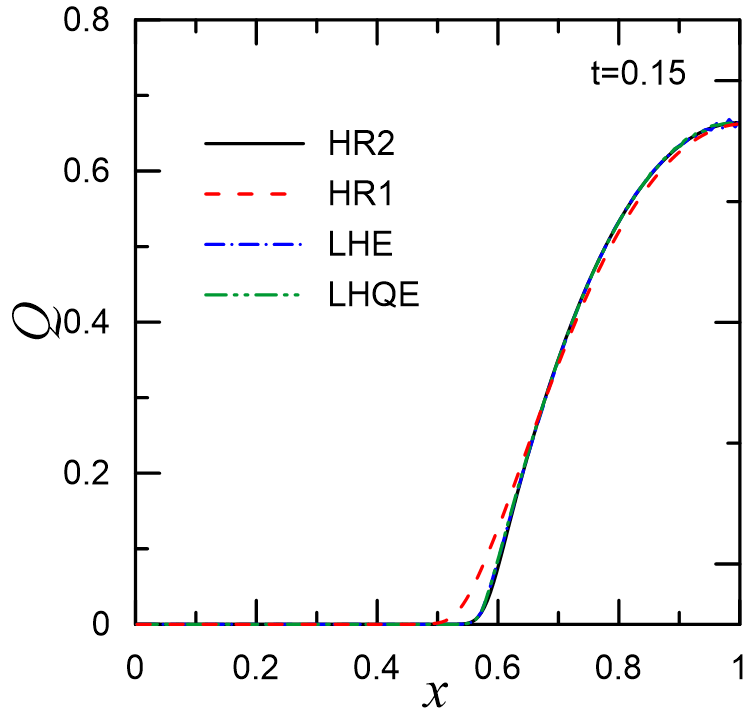}
  \includegraphics[scale=0.92]{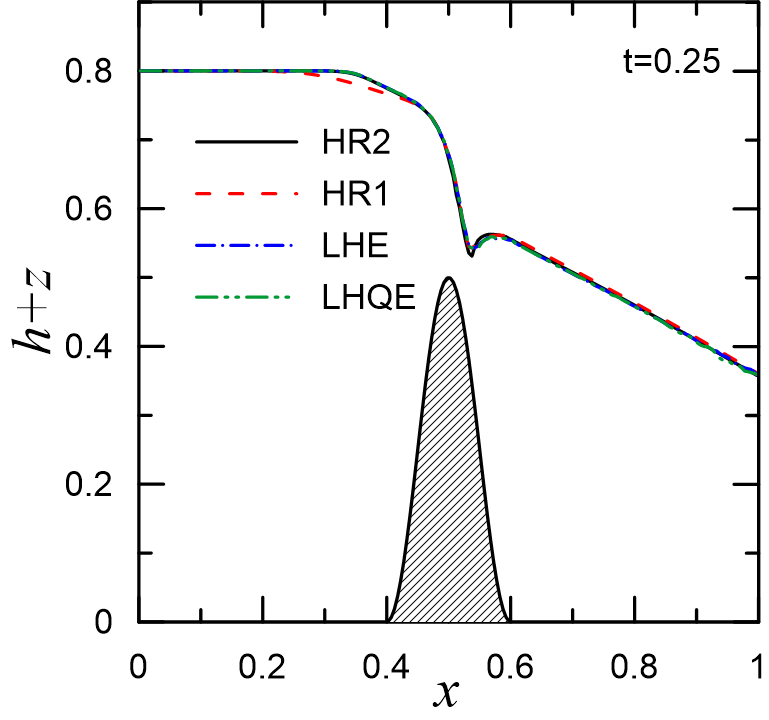}
  \includegraphics[scale=0.92]{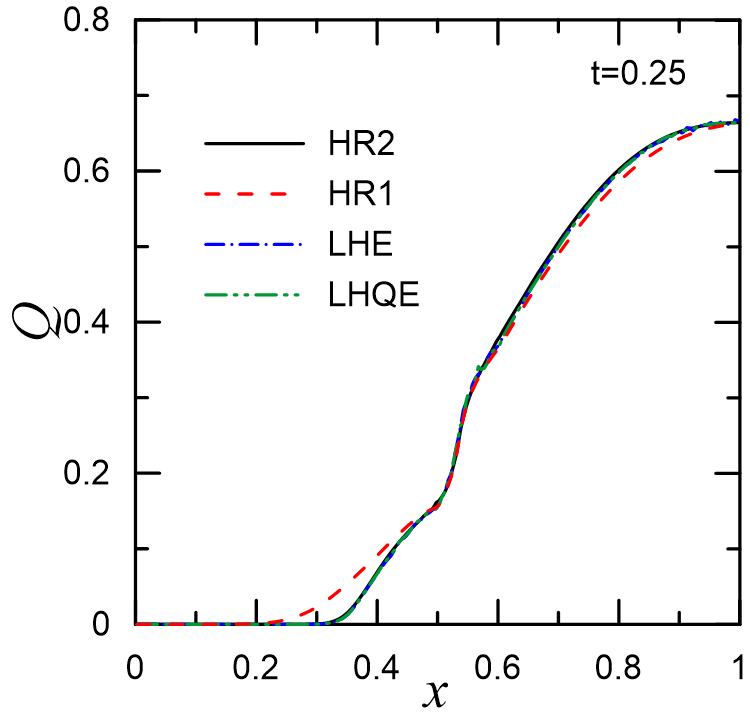}
  \includegraphics[scale=0.92]{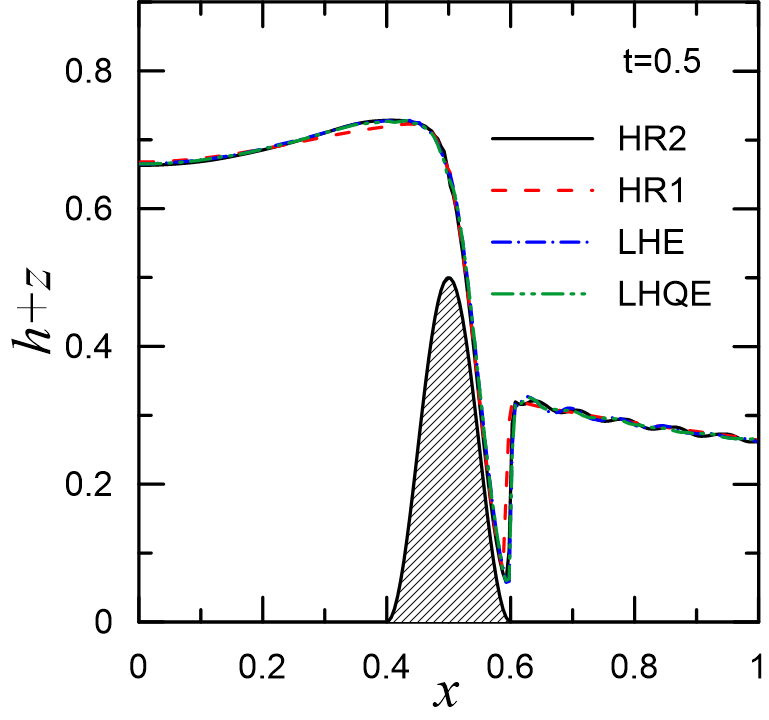}
  \includegraphics[scale=0.92]{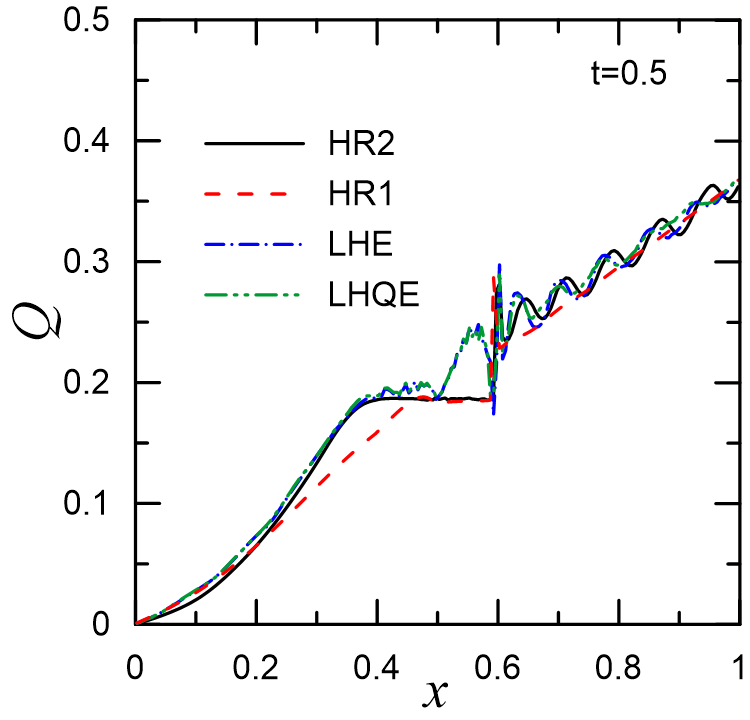}
  \includegraphics[scale=0.92]{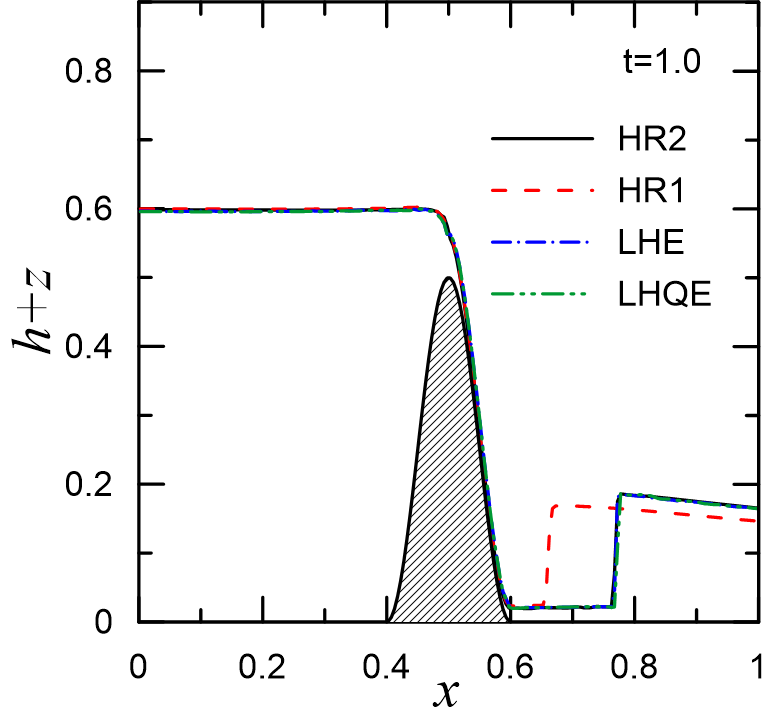}
  \includegraphics[scale=0.92]{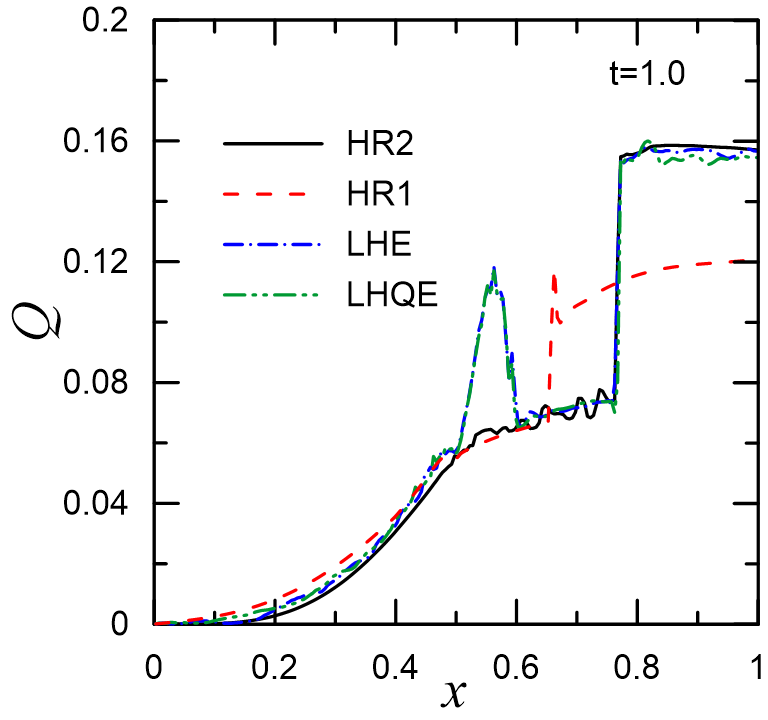}
\caption{Drainage on a non-flat bottom. Water levels and discharges at various times t=0.15,0.25,0.5,1.0 s.}
\label{fig:12}       
\end{figure*}

\subsection{Drainage on a Non-Flat Bottom}
\label{sec:55}

We consider drainage of a symmetric rectangular reservoir to a dry bed through its boundaries, leaving water in topographic depressions. Due to the symmetry, the flow is computed on half the domain, with wall boundary conditions on the left boundary, and open boundary conditions on the right boundary. 
The boundary condition on the right side of the domain allows water that was at rest to flow freely through the right boundary into the originally dry region.
The bottom topography consists of one hump 
\begin{eqnarray}
\label{eq:6.8}
  z(x) = 
  \begin{cases} 
   0.25 \left[ 1+\cos(\pi(x-0.5)/0.1) \right] &    \text{if  } \vert x - 0.5 \vert < 0.1, \\
   0 &   \text{otherwise  } .
  \end{cases} \quad
\end{eqnarray}
After drainage begins, the solution converges to a steady-state solution in which water exists only to the left of the hump.

\begin{figure*}[!h]
  \centering 
  \includegraphics[scale=0.86]{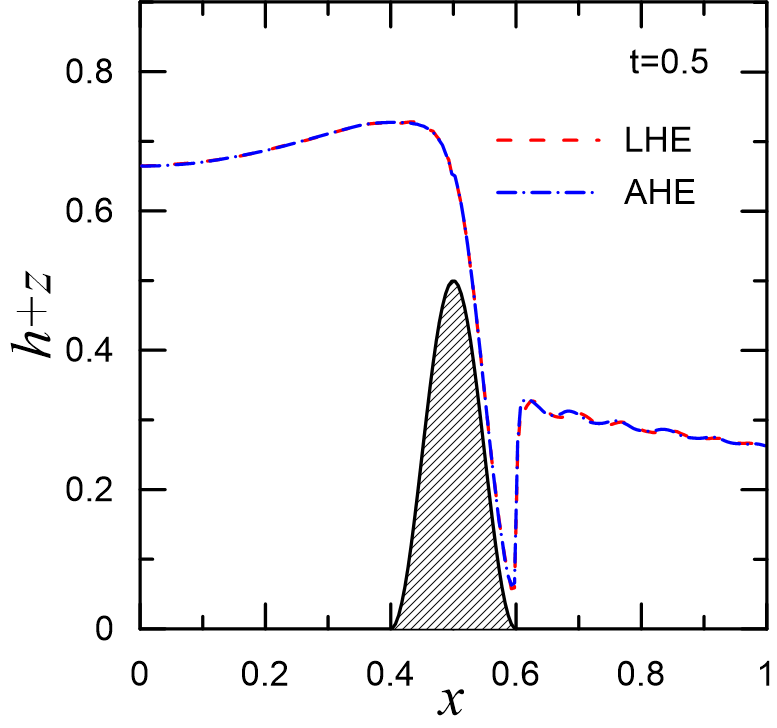}
  \includegraphics[scale=0.86]{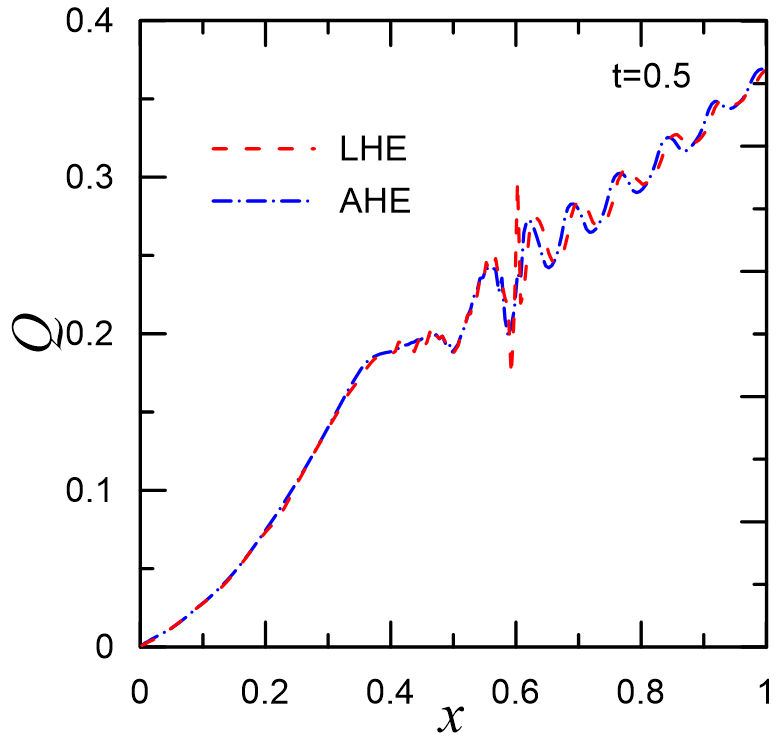}
  \includegraphics[scale=0.86]{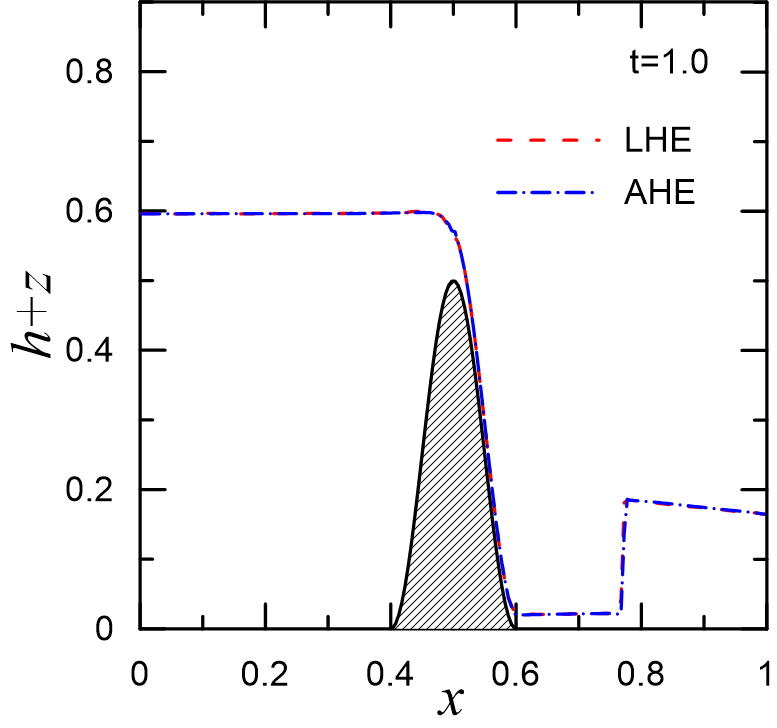}
  \includegraphics[scale=0.86]{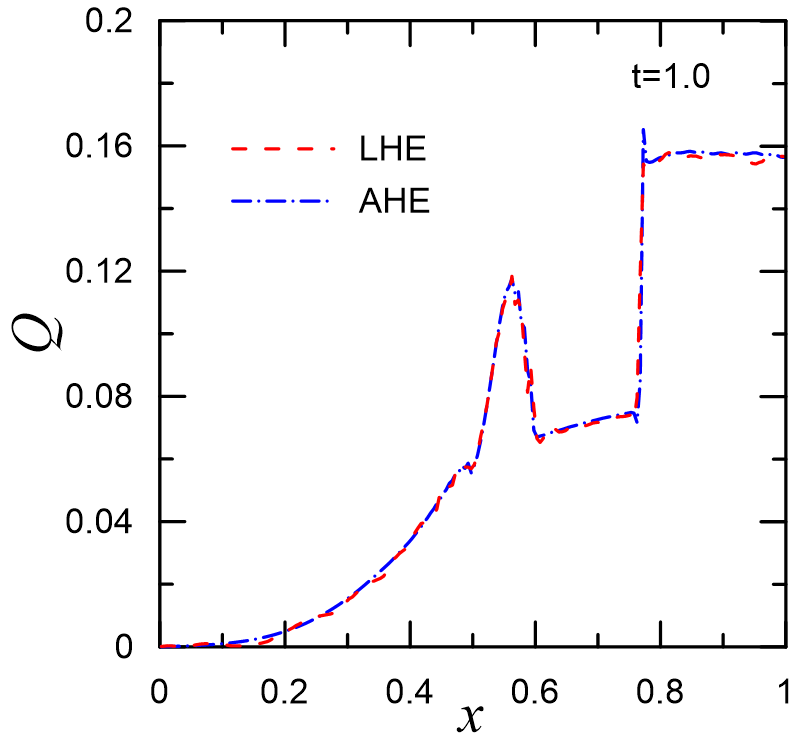}
\caption{Drainage on a non-flat bottom. Water levels and discharges at times t=0.5,1.0 s.}
\label{fig:13}       
\end{figure*}
%
  
\begin{figure*}[!h]
  \centering 
  \includegraphics[scale=0.86]{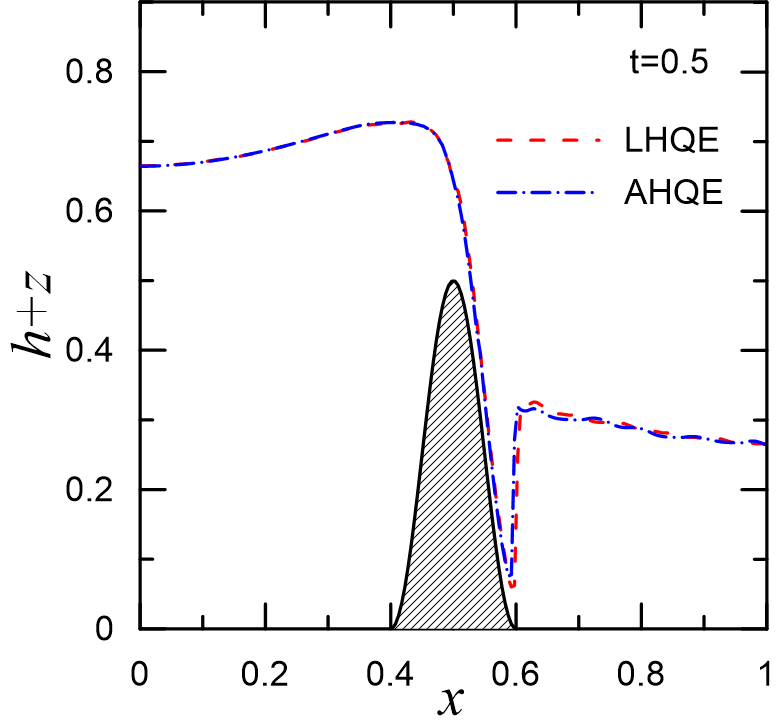}
  \includegraphics[scale=0.86]{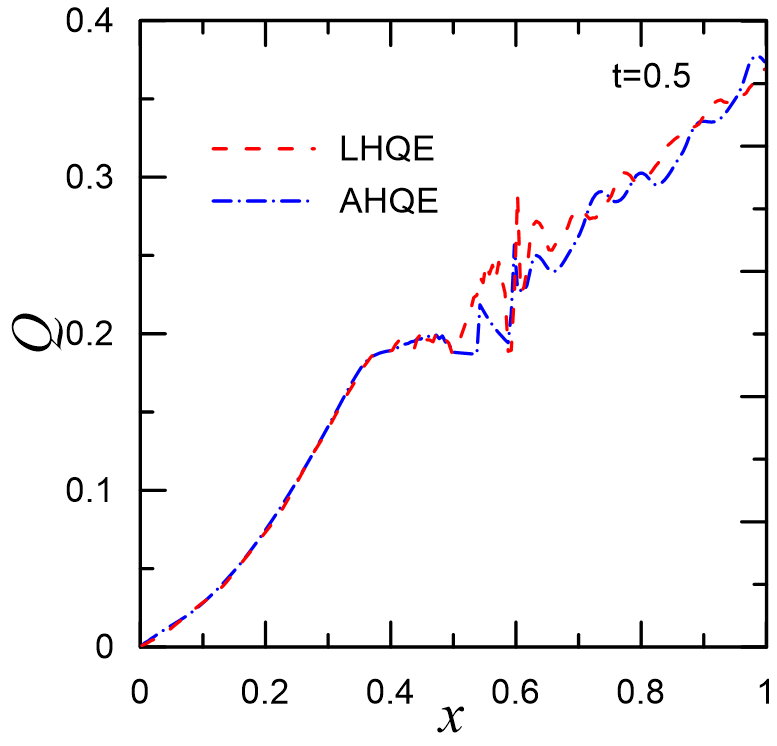}
  \includegraphics[scale=0.86]{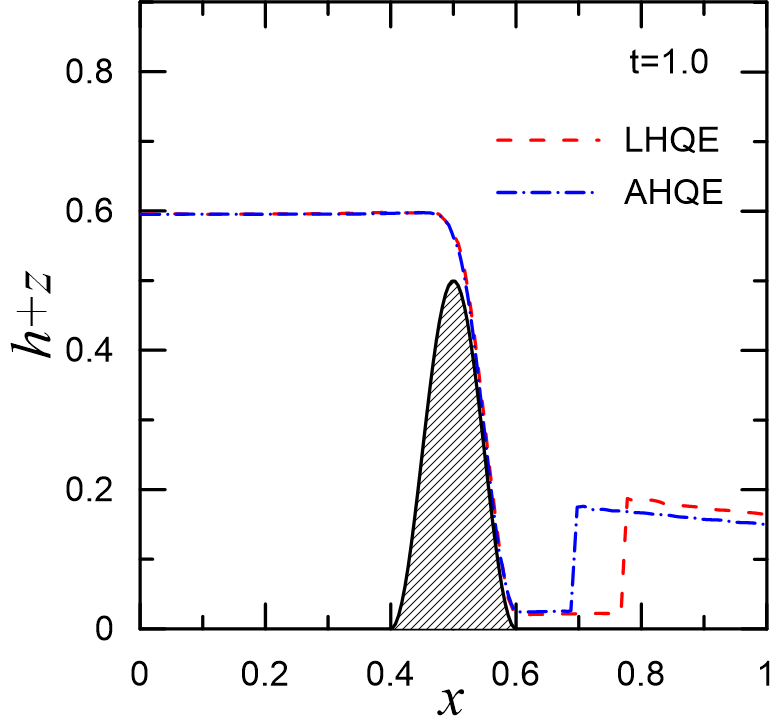}
  \includegraphics[scale=0.86]{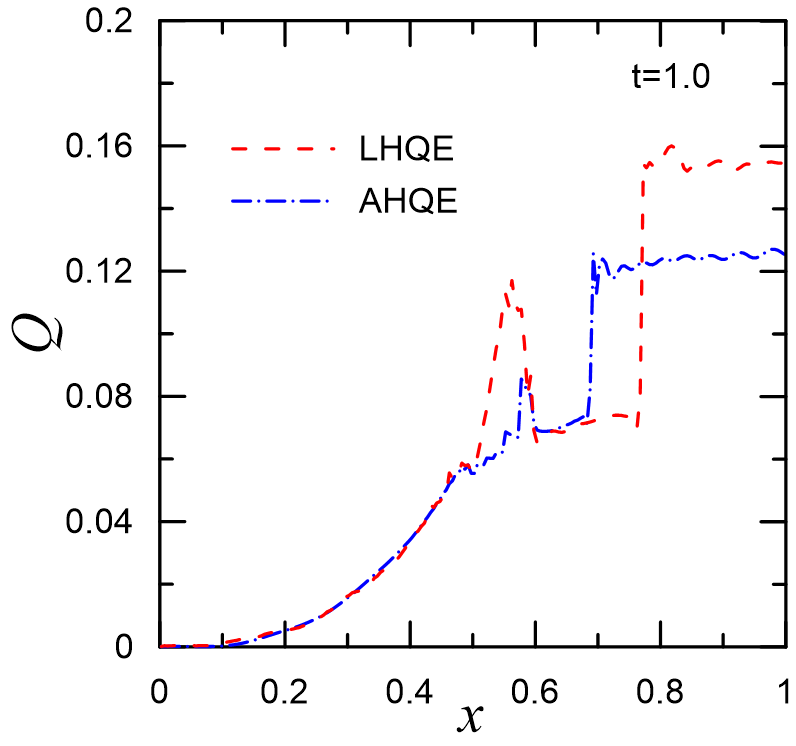}
\caption{Drainage on a non-flat bottom. Water levels and discharges at times t=0.5,1.0 s.}
\label{fig:14}       
\end{figure*}

Numerical results of water flow at different times, obtained on a uniform grid with N=200 cells, are presented in Fig.~\ref{fig:12}-\ref{fig:14}. Fig.~\ref{fig:12} shows that all numerical schemes, except HR1, give similar results for the water surface level.
The first-order HR1 scheme produces a more diffusive water level profile.
The most significant difference in the computed discharges is observed over the right side of the hump.

In Fig.~\ref{fig:13}, the numerical results obtained with the LHE and AHE schemes are in good agreement, in contrast to the results presented in Fig.~\ref{fig:14}.

Note that the AHE and AHQE schemes use flux limiters, which are approximate solutions of the corresponding linear programming problems.  The numerical results in Fig.~\ref{fig:14} show their strong dependence on the numerical diffusion of the applied difference schemes.

\begin{figure*}[!h]
  \includegraphics{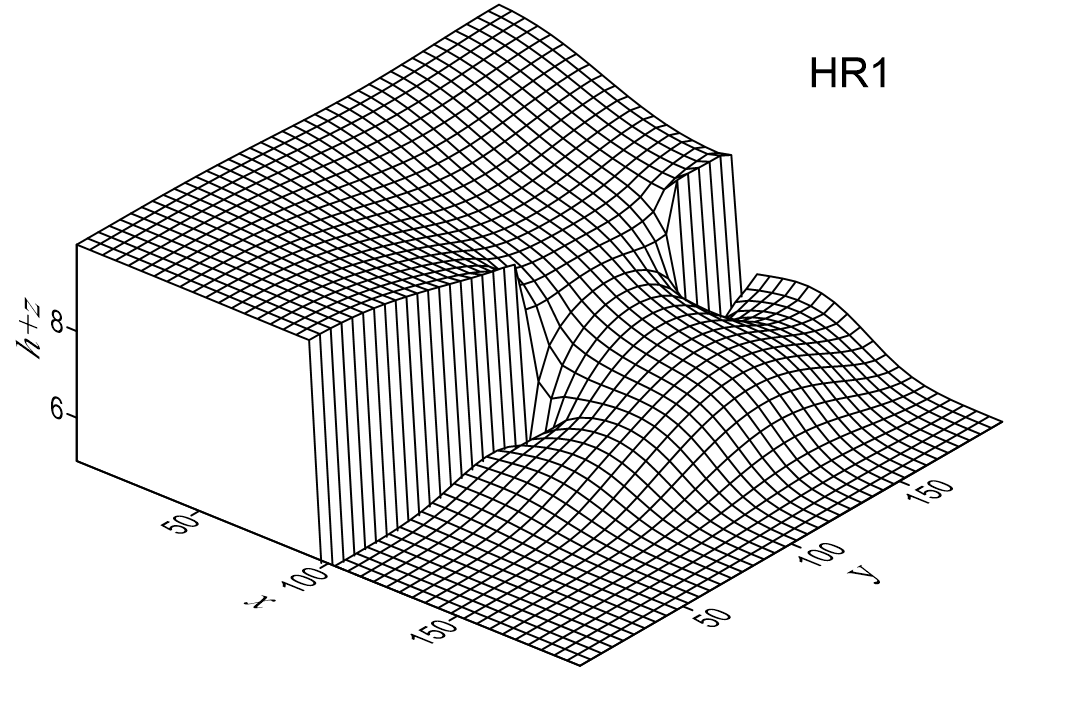}
  \includegraphics{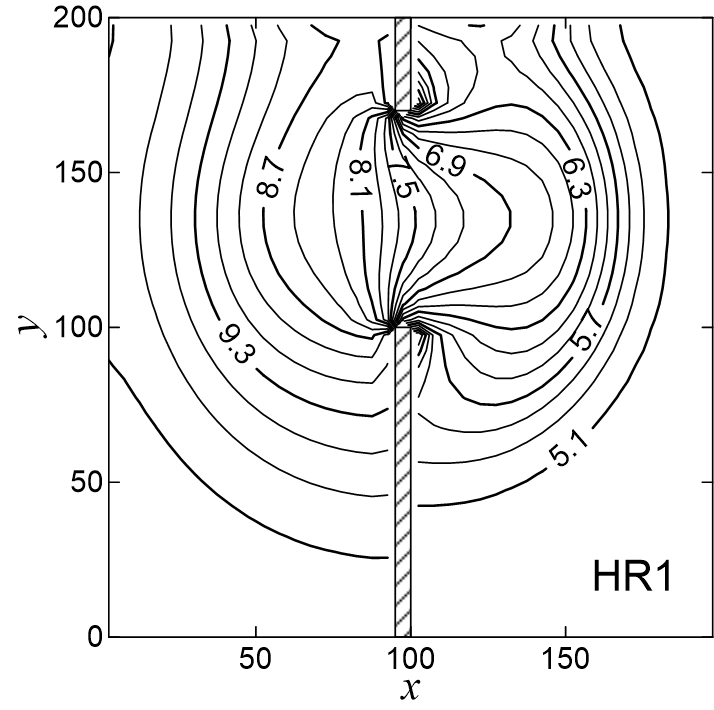}
\caption{Water surface levels and depth contours for the partial dam-break flow at t = 7.2 s computed with the HR1 scheme.}
\label{fig:15}       
\end{figure*}

\subsection{2D Partial Dam Break}
\label{sec:56}

In this section, a partial dam break problem with a nonsymmetrical breach is considered. The spatial domain is defined as a channel with 200 m in length and 200 m in width, the dam is located in the middle of the domain at a distance of 100 m. The bottom is horizontal and frictionless.
Initially, the upstream and downstream water depths are set at 10 and 5 m, respectively.
The breach is 75 m long, located 30 m from the left bank and 95 m from the right bank.

\begin{figure*}[!hbt]
  \centering 
  \includegraphics{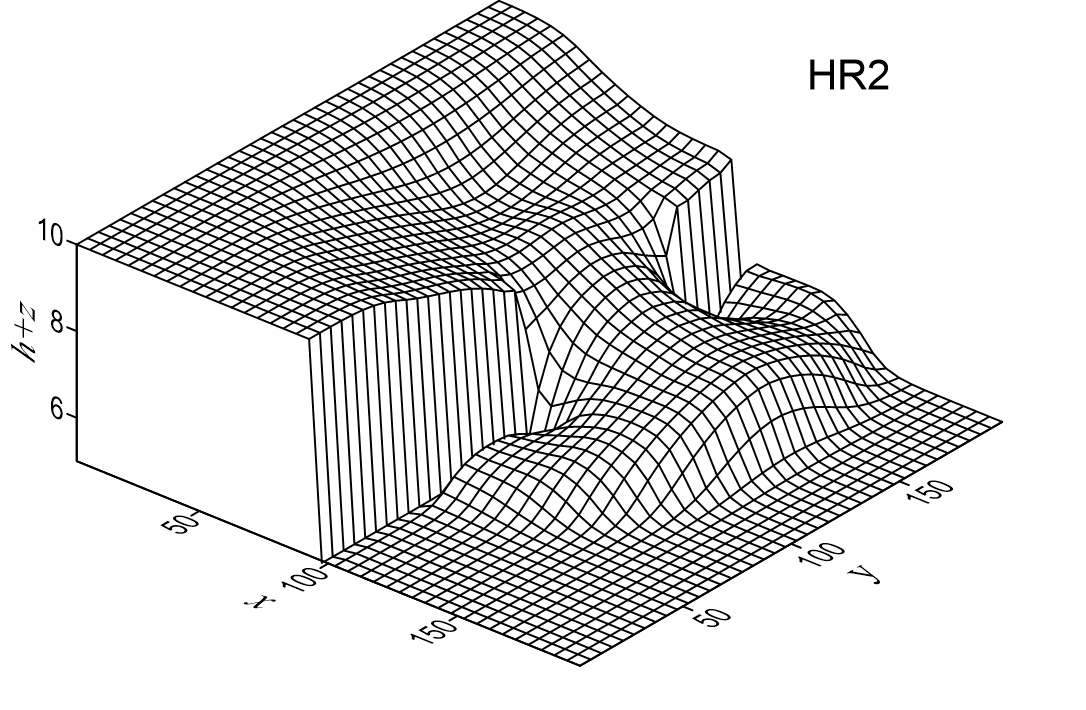}
  \includegraphics{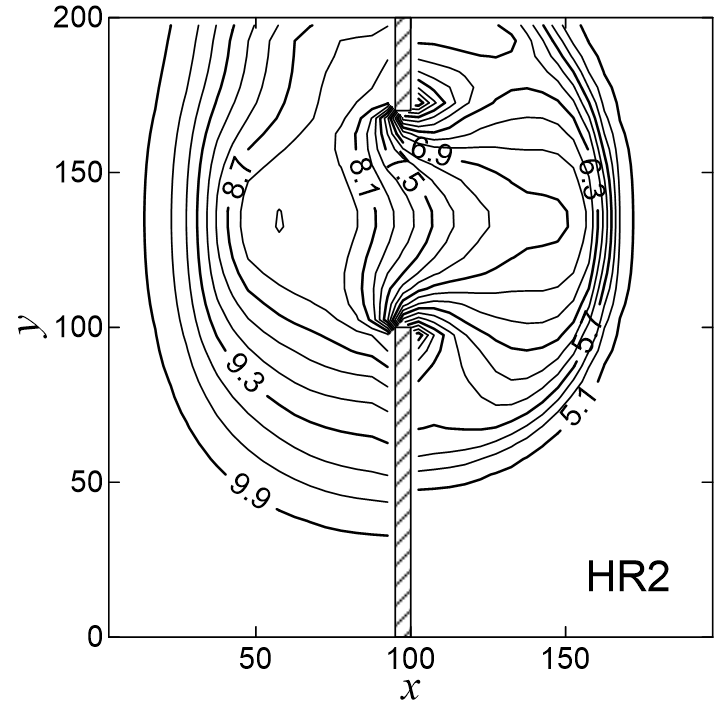}
\caption{Water surface levels and depth contours for the partial dam-break flow at t = 7.2 s computed with the HR2 scheme.}
\label{fig:16}       
\end{figure*}
%
\begin{figure*}[!hbt]
  \centering 
  \includegraphics{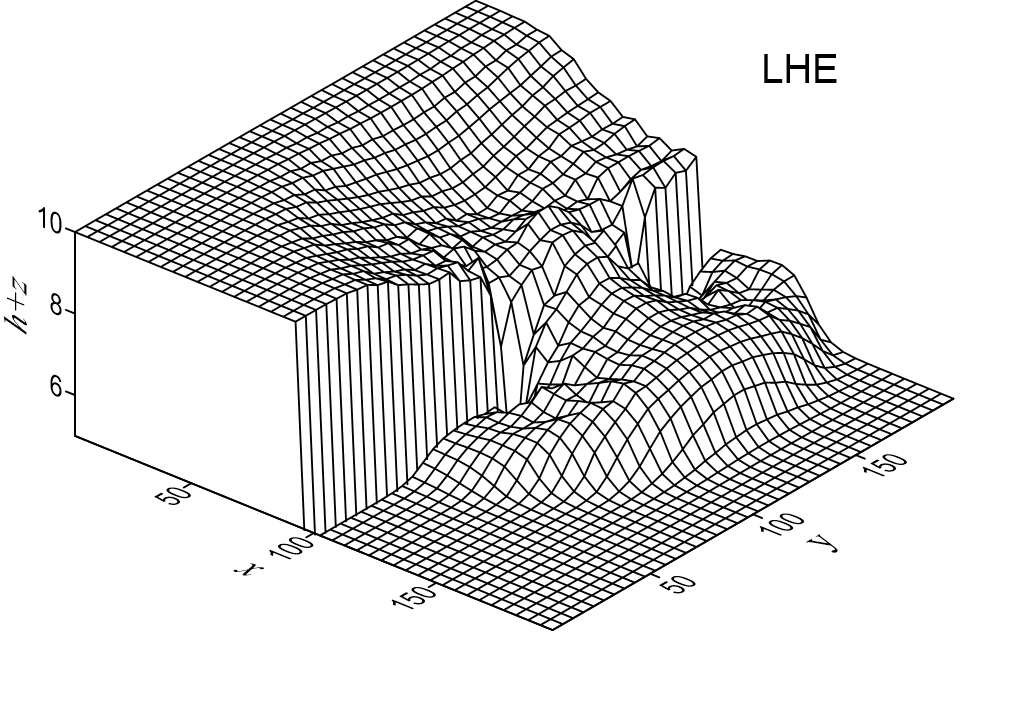}
  \includegraphics{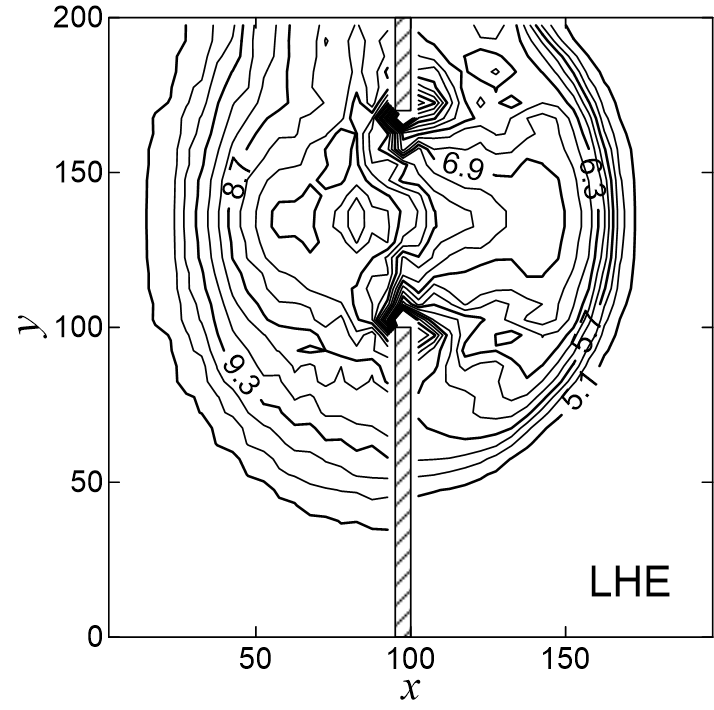}
\caption{Water surface levels and depth contours for the partial dam-break flow at t = 7.2 s computed with the LHE scheme.}
\label{fig:17}       
\end{figure*}
%
\begin{figure*}[!h]
  \centering 
  \includegraphics{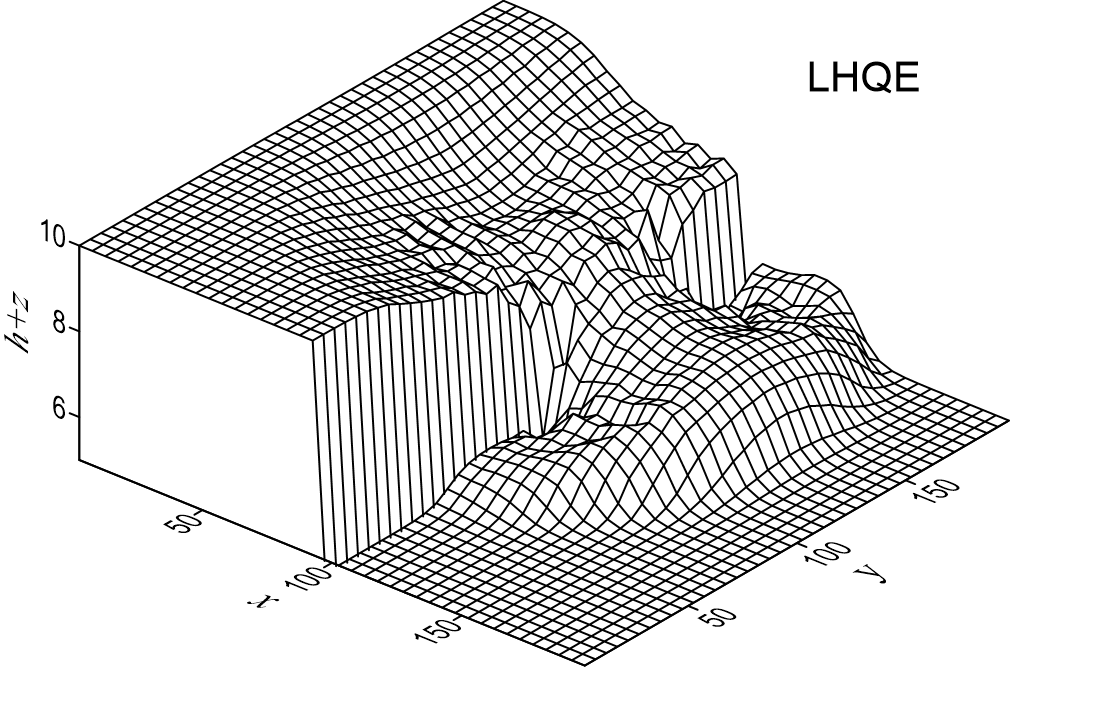}
  \includegraphics{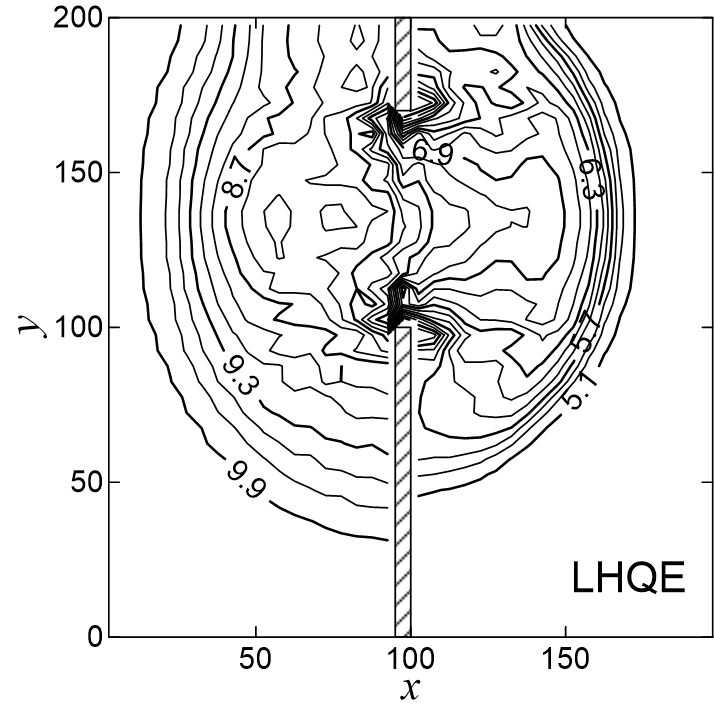}
\caption{Water surface levels and depth contours for the partial dam-break flow at t = 7.2 s computed with the LHQE scheme.}
\label{fig:18}       
\end{figure*}
%

\begin{figure*}[!hbt]
  \centering 
  \includegraphics{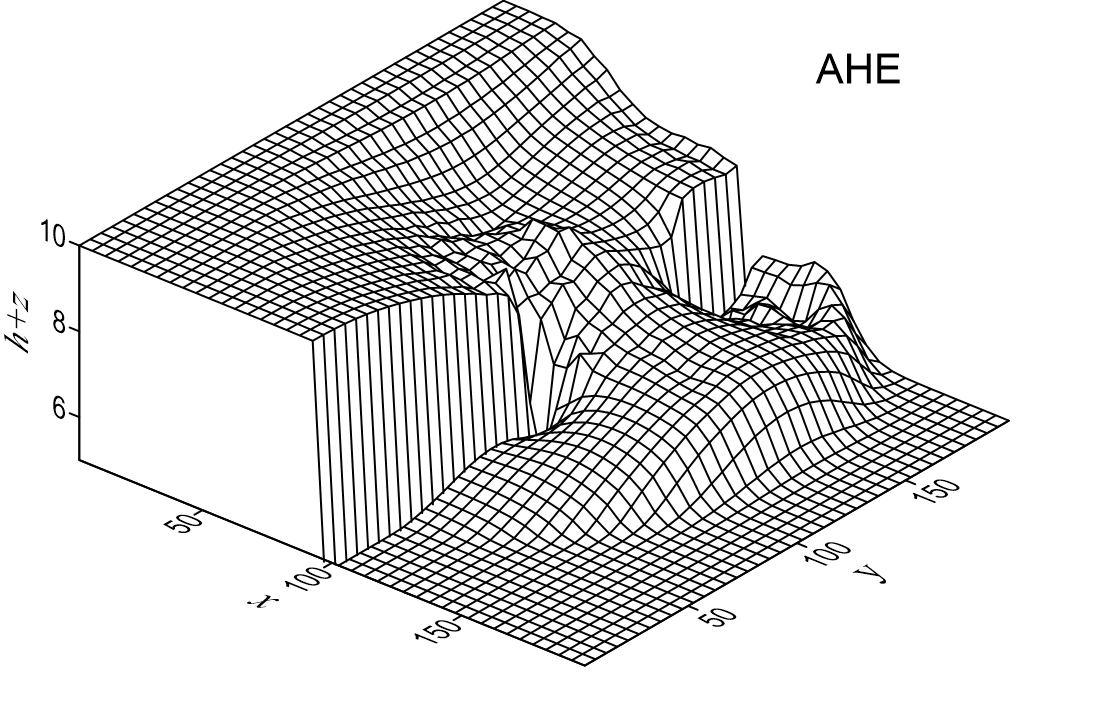}
  \includegraphics{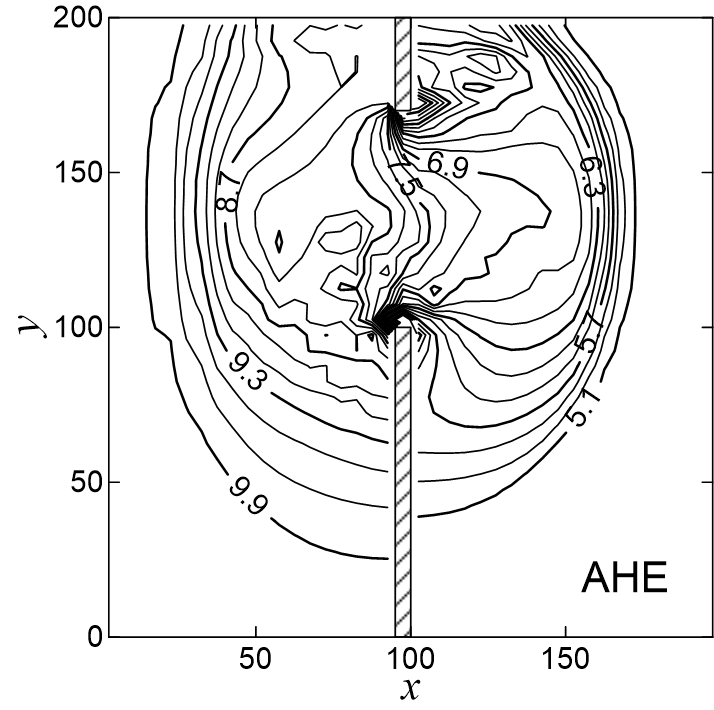}
\caption{Water surface levels and depth contours for the partial dam-break flow at t = 7.2 s computed with the AHE scheme.}
\label{fig:19}       
\end{figure*}
%
\begin{figure*}[!hbt]
  \centering 
  \includegraphics{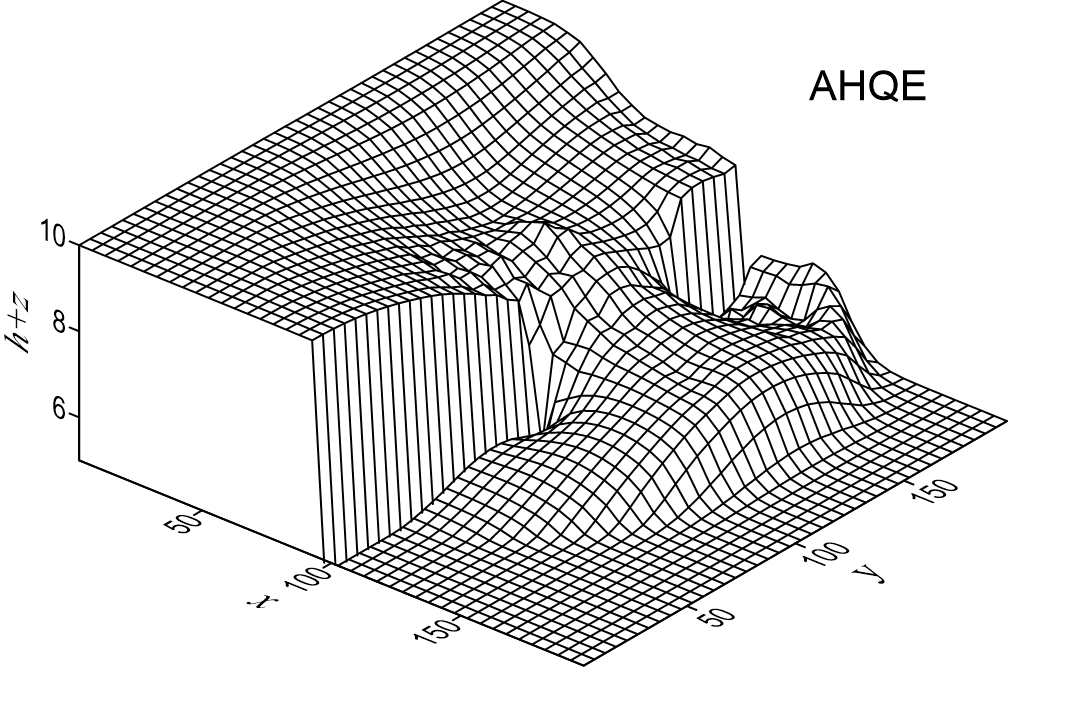}
  \includegraphics{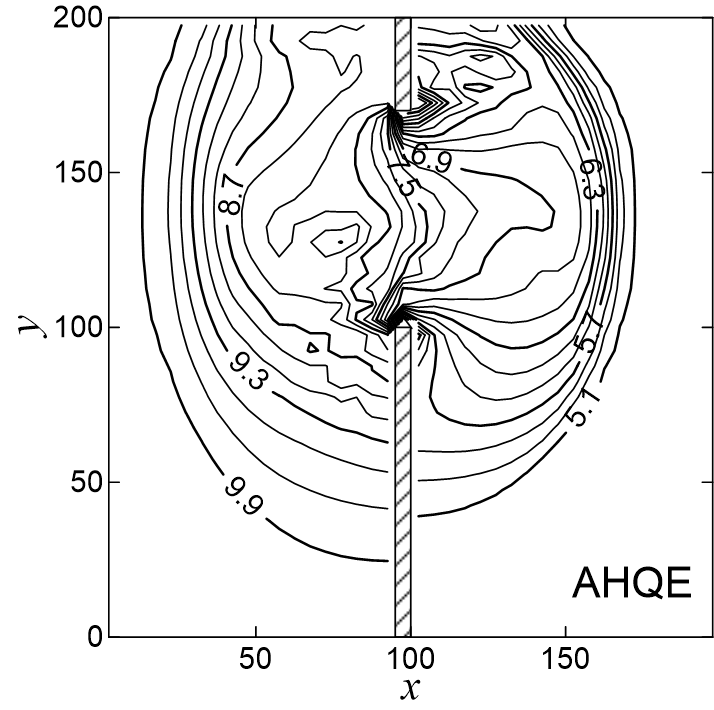}
\caption{Water surface levels and depth contours for the partial dam-break flow at t = 7.2 s computed with the AHQE scheme.}
\label{fig:20}       
\end{figure*}

The computational domain is discretized by a 40 x 40 square grid.
Fig.~\ref{fig:15}-\ref{fig:20} show a three-dimensional view of the water surface levels and water depth contours 7.2 s after the dam failure.
The numerical results obtained with the HR1 scheme are the most diffusive of the others. 
The water surface levels and water depth countours obtained by the LHE, LHQE, AHE, and AHQE schemes are similar to the numerical results of HR2 but are non-smooth. 
AHE and AHQE, whose flux limiters are approximate solutions to the corresponding optimization problems, produce smoother solutions than LHE and LHQE, but their solutions are nonsymmetric about the center of the breach.

\section{Conclusions}
\label{sec:6}

We presented the flux correction design for a hybrid scheme to obtain an entropy-stable solution of shallow water equations with variable topography. The hybrid scheme is an explicit HR scheme whose numerical flux is a convex combination of a first-order Rusanov flux and a high-order flux. We studied the conditions under which a first-order HR scheme with the Rusanov flux satisfied the fully discrete entropy inequality. The flux limiters for the hybrid scheme can be an exact or approximate solution to the corresponding optimization problem in which constraints valid for the first-order HR scheme are applied to the hybrid scheme. It is proved that in the vicinity of a numerical solution of the first-order HR scheme, there is a unique flux correction with flux limiters that are the proposed approximate solution to the optimization problem.

Numerical examples show that the hybrid HR scheme can produce oscillations in numerical results if only water surface level constraints for the optimization problem are used to compute the flux limiters. We also note that numerical results obtained with hybrid HR schemes whose flux limiters are exact and approximate solutions to the optimization problem can differ significantly.


%
 \section*{Conflict of interest}
 The author declare that he has no conflict of interest.
 \section*{Data Availability Statement}
 All data generated or analysed during this study are included in this published article.



\end{document}